\documentclass[a4paper,11pt]{amsart}

\usepackage{amsmath}
\usepackage{amsfonts}
\usepackage{amssymb,latexsym,amsmath,graphicx,epsf}

\newcommand{\dbar}{\overline{\partial}}

\newcommand{\domain}{{\mathcal D}}

\newcommand{\Om}{\Omega}
\newcommand{\ra}{\rightarrow}
\newcommand{\T}{{\mathbf{t}}}

\newcommand{\DOm}{\partial \Omega}
\newcommand{\Rec}{{\mathcal R}}
\newcommand{\rp}{{\alpha}}
\newcommand{\FM}{{\mathcal A}}

\DeclareMathOperator{\de}{\partial}
\DeclareMathOperator{\dez}{\de_z}
\DeclareMathOperator{\dbarz}{\dbar_z}
\DeclareMathOperator{\dbark}{\dbar_k}
\DeclareMathOperator{\R}{\mathbb{R}}
\DeclareMathOperator{\C}{\mathbb{C}}

\DeclareMathOperator{\bndry}{\partial\Omega}
\numberwithin{equation}{section}

\begin{document}
\title[Nonlinear Inversion from Partial EIT Data]{Nonlinear Inversion from Partial EIT Data:\\Computational Experiments}

%    author one information
% \author[short version for running head]{name for top of paper}
\author{S. J. Hamilton}
\address{University of Helsinki \\Department of Mathematics and Statistics, P.O. Box 68, FI-00014 Helsinki, Finland}
\email{{sarah.hamilton@helsinki.fi}}
\thanks{SalWe Research Program for Mind and Body (Tekes - the Finnish Funding Agency for Technology and Innovation grant 1104/10)}

%    author two information
\author{S. Siltanen}
\address{University of Helsinki \\Department of Mathematics and Statistics, P.O. Box 68, FI-00014 Helsinki, Finland}
\email{samuli.siltanen@helsinki.fi}
\thanks{Academy of Finland (Finnish Centre of Excellence in Inverse Problems Research 2012--2017, decision number 250215)}

\keywords{Inverse problem, Numerical solver, Conductivity equation, Inverse conductivity problem, Complex geometrical optics solution, Nonlinear Fourier transform, Electrical impedance tomography}

\subjclass[2010]{Primary 65N21, 35R30; Secondary 45Q05}
%    The 2010 edition of the Mathematics Subject Classification is
%    now available.  If you are citing a classification from the
%    new scheme, use the following input coding instead.
%\subjclass[2010]{Primary }

%\date{}

%\begin{document}

%\maketitle
%--------------------------------------------------------------------
\begin{abstract}
Electrical impedance tomography (EIT) is a non-invasive imaging method in which an unknown physical body is probed with electric currents applied on the boundary, and the internal conductivity distribution is recovered from the measured boundary voltage data. The reconstruction task is a nonlinear and ill-posed inverse problem, whose solution calls for special regularized algorithms, such as D-bar methods which are based on complex geometrical optics solutions (CGOs). In many applications of EIT, such as monitoring the heart and lungs of unconscious intensive care patients or locating the focus of an epileptic seizure, data acquisition on the entire boundary of the body is impractical, restricting the boundary area available for EIT measurements. An extension of the D-bar method to the case when data is collected only on a subset of the boundary is studied by computational simulation. The approach is based on solving a boundary integral equation for the traces of the CGOs using localized basis functions (Haar wavelets). The numerical evidence suggests that the D-bar method can be applied to partial-boundary data in dimension two and that the traces of the partial data CGOs approximate the full data CGO solutions on the available portion of the boundary, for the necessary small $k$ frequencies. 
\end{abstract}
\maketitle
%--------------------------------------------------------------------
\section{Introduction}
%--------------------------------------------------------------------

%.........................................................................................................................
\subsection{EIT and the inverse conductivity problem}
%.........................................................................................................................
Electrical impedance tomography (EIT) is a non-invasive imaging method where an unknown physical body is probed with electric currents, and the internal conductivity distribution is recovered from the measurement data. The reconstruction task is a nonlinear and ill-posed inverse problem, whose solution calls for special regularized algorithms, such as the D-bar method \cite{Knudsen2009}. Applications of EIT include monitoring the heart and lungs of unconscious intensive care patients, industrial process monitoring and underground prospecting. 

\begin{figure}[t] 
\begin{picture}(250,110)
\put(-25,0){\includegraphics[height=3.5cm]{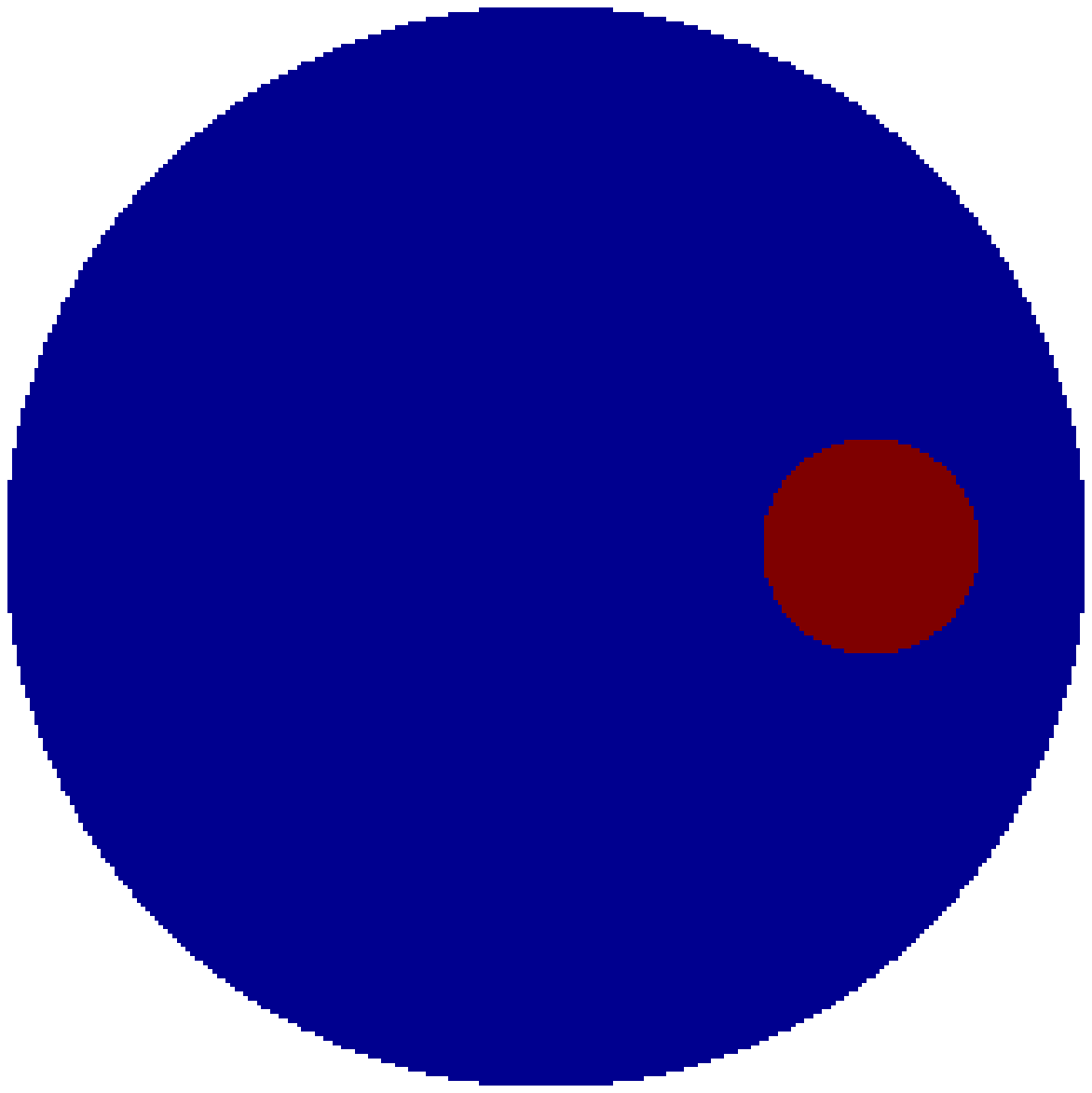}}
\put(75,0){\includegraphics[height=3.5cm]{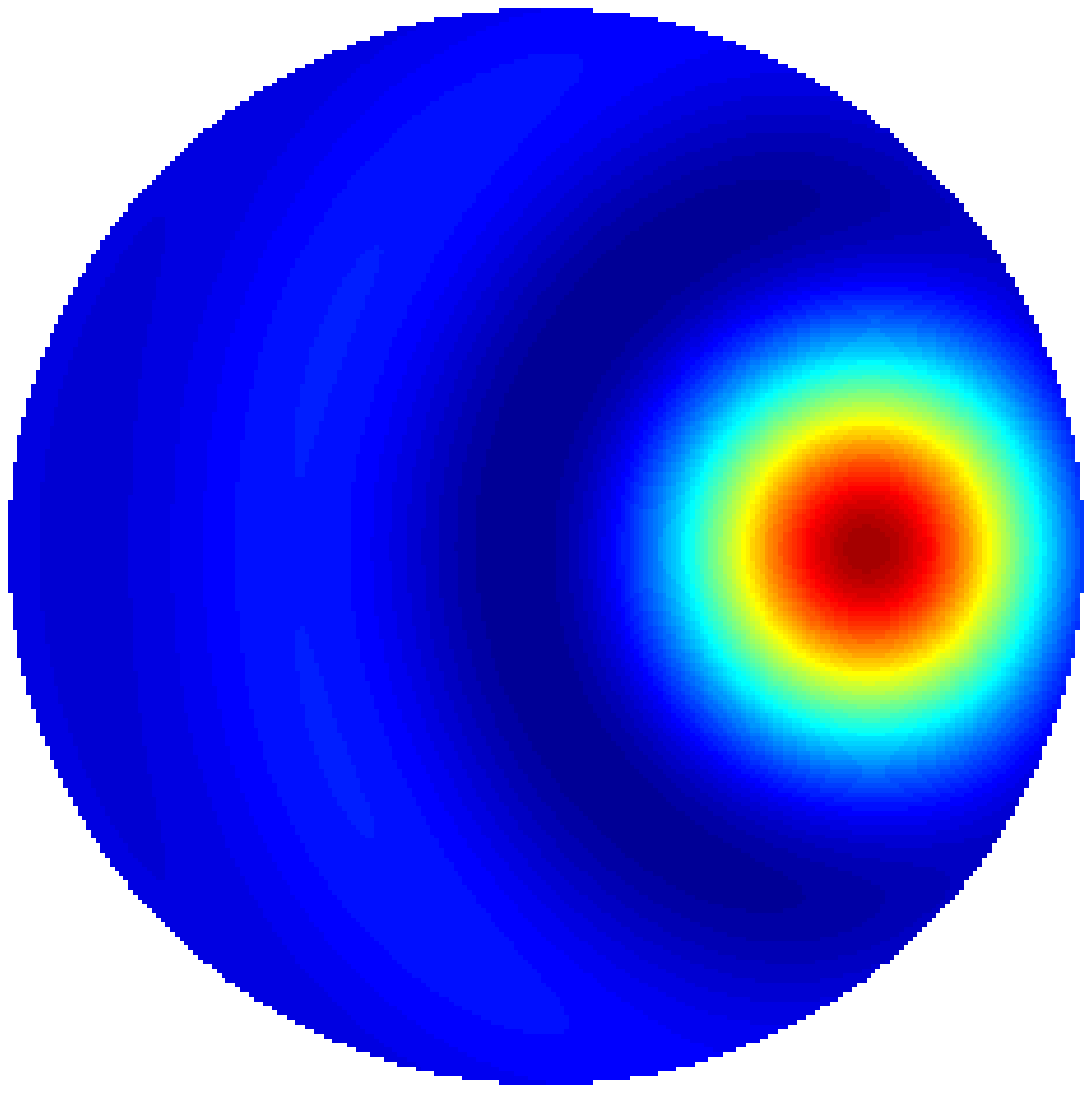}}
\put(175,0){\includegraphics[height=3.5cm]{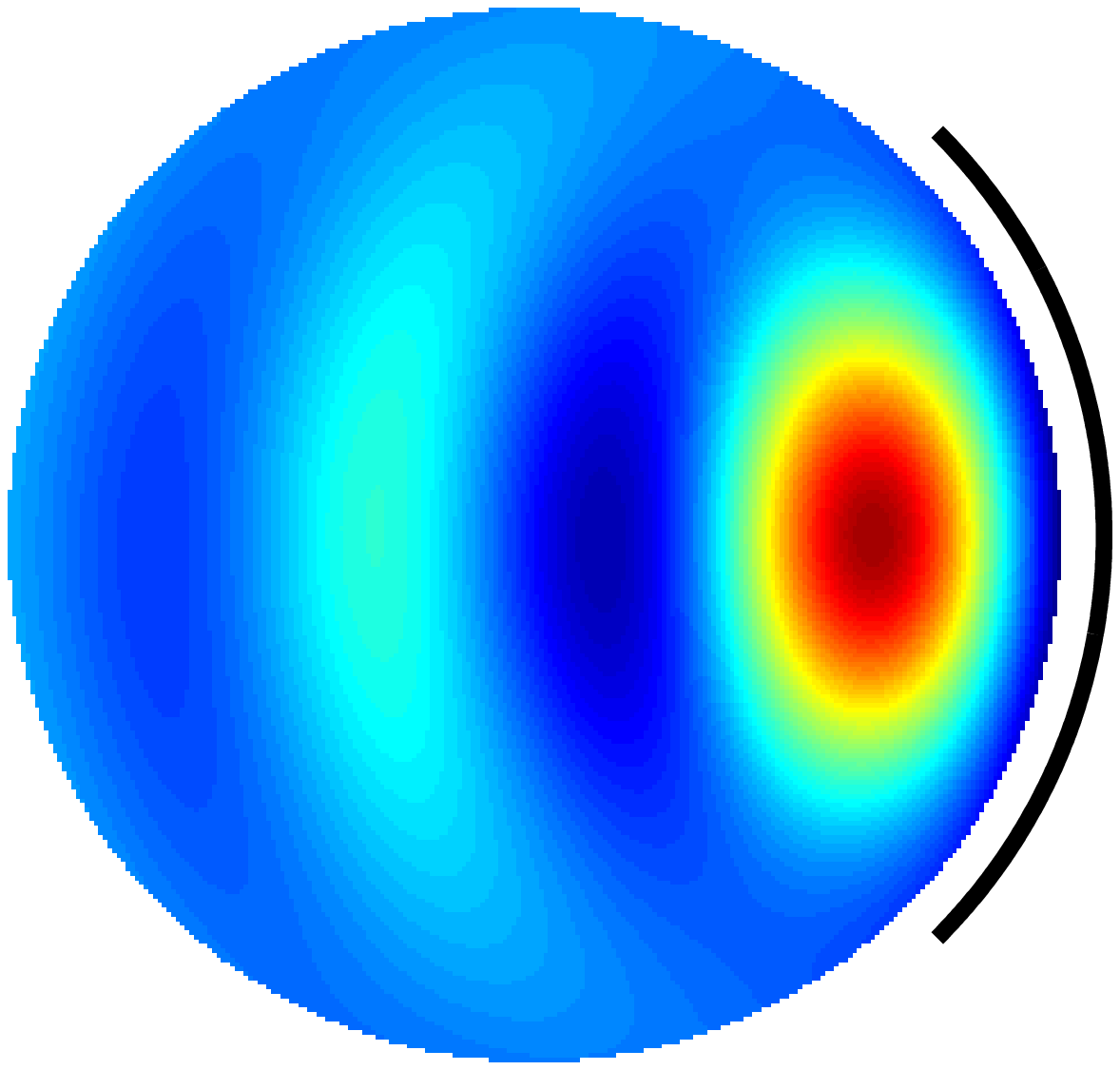}}
\put(275,20){$\Gamma$}
\end{picture}
\caption{\label{Fig:ad}Left: original conductivity. Middle: reconstruction from full-boundary data using the D-bar method. Right: reconstruction from partial-boundary data using the proposed method. The subset $\Gamma\subset\DOm$ where the measurements are available is denoted by a black line, which corresponds to 25\% of the entire boundary}
\end{figure}

Practical considerations typically restrict the boundary area available for EIT measurements: for example, it is not sensible to cover a patient completely with electrodes when imaging the heart. In this paper we study a possible extension of the D-bar method to the case when data is collected only on a subset of the boundary. See Figure \ref{Fig:ad}.
%+++++++++++++++++++++++++++++++
%+++++++++++++++++++++++++++++++
The mathematical model of EIT is the inverse conductivity problem introduced by Calder\'on \cite{Calder'on1980}. Let $\Om\subset\R^n$ be a bounded and simply connected set with a smooth boundary $\partial\Om$. Let $\sigma: \Om\to \R$ be an essentially bounded measurable function satisfying $\sigma(x)\geq c>0$ for almost every $x\in \Om$. Let $u\in H^{1}(\Om)$ be the unique solution
to
\begin{eqnarray}
  \nabla\cdot\sigma\nabla u &=& 0 \text{ in }\Om, \label{Calderon:eq1} \\
  u\big|_{\partial\Om} &=& \phi\in H^{1/2}(\partial \Om).\label{Calderon:eq2}
\end{eqnarray}
The inverse conductivity problem is to recover the conductivity $\sigma$  from the Dirichlet-to-Neumann (D-N)  map defined by
$$
  \Lambda_{\sigma} : \phi \mapsto \sigma\frac{\partial u}{\partial \nu}\Big|_{\partial \Om}.
$$
Here $\nu=\left(\nu_1,\nu_2\right)=\nu_1+i\nu_2$ is the unit outward facing normal vector to the boundary. Here $\phi$ is a voltage distribution applied on the boundary, and $\Lambda_{\sigma}  \phi$ is the resulting current flux through the boundary. Therefore, $\Lambda_\sigma$ can be seen as an ideal infinite-precision model of practical voltage-to-current measurements.

\begin{figure}[t] 
\begin{picture}(300,145)
\put(-5,0){\includegraphics[width=11cm]{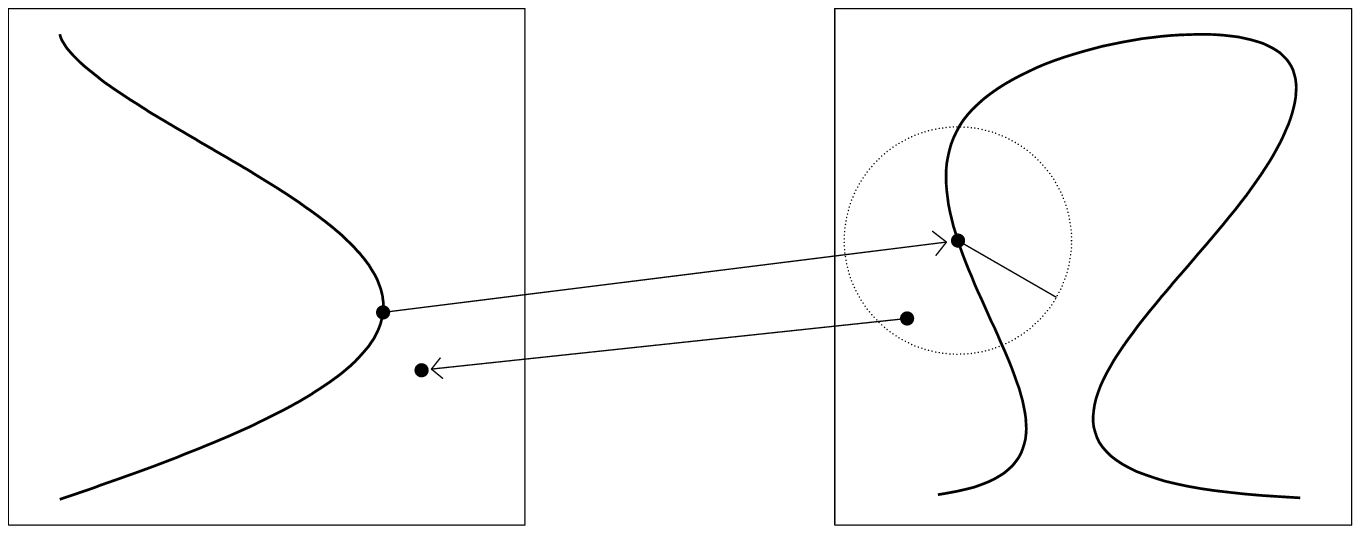}}
\put(-5,125){\small Model space}
\put(185,125){\small \small Data space}
\put(80,107){\small $L^\infty(\Om)$}
\put(195,107){\small $Y$}
\put(72,48){\small $\sigma$}
\put(148,61){\small $\FM$}
\put(232,56){\small $\delta$}
\put(148,33){\small $\Rec_\rp$}
\put(217,67){\small $\Lambda_\sigma$}
\put(265,15){\small $\FM(\domain(\FM))$}
\put(10,18){\small $\domain(\FM)$}
\put(205,45){\small $\Lambda_\sigma^\delta$}
\put(75,24){\small $\Rec_\rp(\Lambda_\sigma^\delta)$}
\end{picture}
\caption{\label{Fig:nonlinear_reg}Schematic illustration of nonlinear regularization\index{regularization!nonlinear} of the {\sc eit} problem. Here the forward map is defined as $\FM(\sigma)=\Lambda_\sigma$ with the domain of definition denoted by $\domain(\FM)\subset L^\infty(\Om)$. The conductivity $\sigma$ is approximately recovered as $\Rec_\rp(\Lambda_\sigma^\delta)$.}
\end{figure}

Calder\'on asked two main questions in his seminal article \cite{Calder'on1980}: 
\begin{itemize}
\item[(i)] Is $\sigma$ uniquely determined by $\Lambda_\sigma$?
\item[(ii)] If the answer to (i) is yes, how can one calculate $\sigma$ from $\Lambda_\sigma$? 
\end{itemize}
In practical EIT imaging only a finite-range and noisy approximate operator $\Lambda_\sigma^\delta$ is available. In general, $\Lambda_\sigma^\delta$ is not the D-N map of any conductivity. We usually only know that $\|\Lambda_\sigma^\delta-\Lambda_\sigma\|_Y\leq \delta$. Here $Y$ is an appropriate data space and $\delta>0$ can be determined from the properties of the measurement device. This leads us to a third question:
\begin{itemize}
\item[(iii)] Given $\Lambda_\sigma^\delta$ and $\delta$, how can one design a continuous map from $Y$ to $L^\infty(\Omega)$ whose output is a useful approximation to $\sigma$?
\end{itemize}
As the inverse conductivity problem is ill-posed, the forward map $\FM:\sigma\mapsto\Lambda_\sigma$ does not have a continuous inverse. Therefore, question (iii) needs to be answered by constructing a {\em regularization strategy} \cite{Engl1996}.
More precisely, a family of continuous mappings $\Rec_\rp:Y\ra L^\infty(\Om)$ must be defined, parameterized by $0<\rp<\infty$, such that
\begin{equation}\label{regstrat}
  \lim_{\rp\ra 0}\|\Rec_\rp(\Lambda_\sigma)-\sigma\|_{L^\infty(\Om)}=0,
\end{equation}
for each fixed $\sigma$. Note that (\ref{regstrat}) is closely related to question (ii) above. Furthermore, one needs to specify a choice $\rp=\rp(\delta)$ for the regularization parameter as a function of the noise level so that $\rp(\delta)\ra 0\mbox{ as }\delta\ra 0$. Finally, the reconstruction error $\|\Rec_{\rp(\delta)}(\Lambda_\sigma^\delta)-\sigma\|_{L^\infty(\Om)}$ must vanish in the zero noise limit:
for any fixed $\sigma$ we must have 
\begin{equation}\label{regstrat(ii)}
  \sup_{\Lambda_\sigma^\delta\in Y} \left\{\|\Rec_{\rp(\delta)}(\Lambda_\sigma^\delta)-\sigma\|_{L^\infty(\Om)} \,:\, \|\Lambda_\sigma^\delta-\Lambda_\sigma\|_Y\leq \delta\right\}\ra 0\mbox{ as }\delta\ra 0.
\end{equation}
For more details, see Figure \ref{Fig:nonlinear_reg} and \cite{Knudsen2009,Mueller2012}. 

%+++++++++++++++++++++++++++++++
%+++++++++++++++++++++++++++++++
%\afterpage{\clearpage}
%.........................................................................................................................
\subsection{D-bar methods for full-boundary data}
%.........................................................................................................................
From the practical viewpoint, the solution of the inverse conductivity problem is a computational algorithm that implements a regularization strategy $\Rec_{\rp}$ satisfying (\ref{regstrat}) and (\ref{regstrat(ii)}). Achieving such a goal is typically a large project involving several milestones, often corresponding to one of the following two types:
\begin{itemize}
\item[(a)] A theoretical breakthrough that outlines a computational approach
\item[(b)] Successful computational experiments that inspire further theoretical study 
\end{itemize}
Let us review the history of a specific two-dimensional D-bar method for EIT in light of (a) and (b). 

{\bf 1996(a):} Nachman showed uniqueness (i) and introduced a infinite-precision reconstruction method (ii) for twice differentiable conductivities in \cite{Nachman1996}. The proof used a nonlinear Fourier transform based on so-called \emph{complex geometrical optics} (CGO) solutions, first defined by Faddeev in 1966 \cite{Faddeev1966} and later rediscovered in 1987 by Sylvester and Uhlmann in the context of 3D EIT \cite{Sylvester1987}. Thus, \cite{Nachman1996} represents a breakthrough in the form of (a) since it is the basis of the first numerical D-bar method \cite{Siltanen2000,Siltanen2001a,Mueller2003}. 

{\bf 2004(b):} Isaacson {\em et al.} demonstrated in \cite{Isaacson2004,Isaacson2006} that the D-bar method performs well on practical data measured from laboratory phantoms and from human subjects. The mandatory regularization step was provided by low-pass filtering in the nonlinear frequency domain. The need for such filtering is evident from the structure of the experimental nonlinear Fourier transforms: they blow up outside a disc centered at the origin.

{\bf 2009(a):} The numerical evidence from practical imaging experiments \cite{Isaacson2004,Isaacson2006} inspired a rigorous regularization proof of convergence  in the form of (\ref{regstrat(ii)}), see \cite{Knudsen2009}. This gave an answer to (iii) and outlined a method for choosing the regularization parameter as the inverse of the nonlinear cutoff frequency. We outline the reconstruction method in Section \ref{sec:Method1} below.

There is an analogous history for other uniqueness proofs and related algorithms in two-dimensions. We review them briefly below without specifying explicitly the progress steps of types (a) and (b). 

Brown and Uhlmann were able to prove uniqueness for real-valued conductivities assuming only one derivative in \cite{Brown1997}. This result was complemented by constructive steps and numerical implementation by Knudsen and Tamasan \cite{Knudsen2004a,Knudsen2002,Knudsen2003}; see also \cite{Knudsen2004}. Francini \cite{Francini2000} extended the uniqueness proof to complex conductivities whose real and imaginary parts are twice differentiable, and her approach was subsequently implemented in \cite{Hamilton2012,Hamilton_Thesis_2012,HM12_NonCirc,Natalia_Thesis}. We outline this reconstruction method in Section \ref{sec:Method2} below. Both methods involve transforming \eqref{Calderon:eq1} to a first order system of $\dez$ and $\dbarz$ equations. 

Astala and P\"{a}iv\"{a}rinta answered Calder\'on's questions (i) and (ii) in their original smoothness category $\sigma\in L^\infty(\Om)$, see \cite{Astala2006a,Astala2006}. This approach has been implemented numerically as well \cite{Astala2010,Astala2011}.

Despite the above developments, some questions still remain open:
\begin{itemize}
\item Is it possible to give a regularization analysis (iii) for less smooth conductivities than twice differentiable? There is numerical evidence of type (b) available since all of the above EIT methods produce noise-robust images when applied to data arising from discontinuous conductivities and regularized by nonlinear low-pass filtering \cite{Knudsen2008a,Knudsen2007,Hamilton2012,Astala2011}.
\item Can the D-bar methodology be used in the case of partial-boundary data?  We discuss this in Section \ref{sec:partialdata} below in the two-dimensional case. 
\end{itemize}

%.........................................................................................................................
\subsection{Extension to partial-boundary data}\label{sec:partialdata}
%.........................................................................................................................
It is of high practical importance to be able to compute EIT reconstructions from data measured only on a part of the boundary. One possibility for designing such algorithms would be to take one of the recent theoretical breakthroughs, such as \cite{Knudsen2006,Kenig2007a,Nachman2010,Imanuvilov2010}, and implement it in the spirit of (a) above. However, we do not discuss such approaches in this paper. We proceed along (b) and produce novel numerical evidence suggesting that it may be possible to use the classical D-bar approach for partial data reconstructions. It is our hope that these computational results inspire further theoretical advances.

Our starting point is the assumption that only a proper subset $\Gamma\subset\DOm$ is available for measurements. We consider voltage-to-current data represented ideally by the restricted D-N map $\widetilde{\Lambda}_\sigma$, defined as follows. Let $\widetilde{\phi}\in H^{1/2}(\partial \Om)$ satisfy $\mbox{supp}(\widetilde{\phi})\subset\Gamma$ and let $u\in H^{1}(\Om)$ be the unique solution of the conductivity equation
\begin{eqnarray}
  \nabla\cdot\sigma\nabla u &=& 0 \text{ in }\Om, \label{Calderon:eq1pb} \\
  u\big|_{\partial\Om} &=& \widetilde{\phi}.\label{Calderon:eq2pb}
\end{eqnarray}
Our partial D-N map is then defined by 
\begin{equation}
  \widetilde{\Lambda}_{\sigma} : \widetilde{\phi} \mapsto \sigma\frac{\partial u}{\partial \nu}\Big|_{\Gamma}.
\end{equation}
The practical data is a finite-range and noisy approximate operator $\widetilde{\Lambda}_\sigma^\delta$ satisfying $\|\widetilde{\Lambda}_\sigma^\delta-\widetilde{\Lambda}_\sigma\|_Y\leq \delta$.

Let us briefly explain our approach in the context of the regularized D-bar method \cite{Knudsen2009} based on Nachman's uniqueness proof \cite{Nachman1996}. In the full-boundary data case, it begins by solving this Fredholm integral equation of the second kind for the (approximate) traces of the CGO solutions on $\DOm$:
\begin{equation}\label{fullBIE1.0}
  \psi(z,k) = e^{ikz} - \int_{\bndry}G_k(z-\zeta)(\Lambda_\sigma^\delta-\Lambda_1)\psi(\zeta,k)\;dS(\zeta), \quad z\in\bndry,
\end{equation}
where $G_k$ is the Faddeev Green's function \cite{Faddeev1966}, here defined in the sense of tempered distributions,
\begin{equation}\label{def:Gk}
   G_k(z) := e^{ikz}g_k(z),\quad g_k(z) := \frac{1}{(2\pi)^2}  \int_{\R^2}\frac{e^{iz\cdot\xi}}{|\xi|^2+2k(\xi_1+i\xi_2)}d\xi.
\end{equation}

In the case of partial-boundary data, we solve the following equation for the unknown functions $\omega(\,\cdot\,,k):\Gamma\ra\C$:
\begin{equation}\label{partialBIE1.0}
  \omega(z,k) = e^{ikz}- \int_{\Gamma}G_k(z-\zeta)(\widetilde{\Lambda}_\sigma^\delta-\widetilde{\Lambda}_1)\omega(\zeta,k)\;dS(\zeta), \quad z\in\Gamma\subset\bndry.
\end{equation}
Now the hypothesis is that
\begin{equation}\label{mainapprox}
  \psi(z,k)|_\Gamma \approx \omega(z,k), \quad z\in\Gamma\subset\bndry,\quad\text{for some $k\in\C$}.
\end{equation}
If (\ref{mainapprox}) holds, it opens up a variety of extensions of D-bar methods to partial-boundary data applications.  

%.........................................................................................................................
\subsection{Focus of this paper}
\begin{figure}[t]
\hspace{1em}
\begin{picture}(300,220)
\put(-10,90){\includegraphics[width=2.4cm]{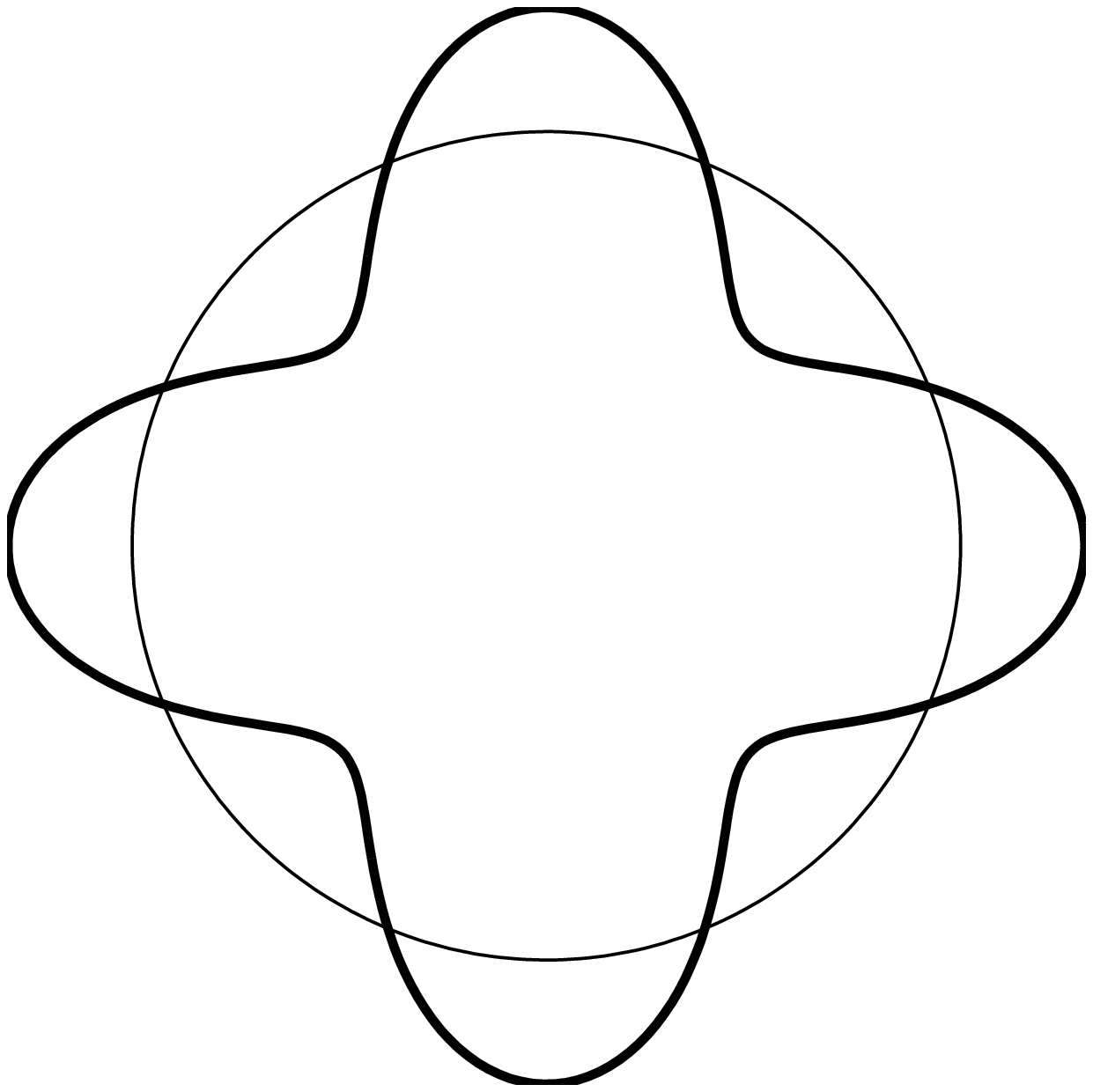}}
\put(-10,0){\includegraphics[width=2.4cm]{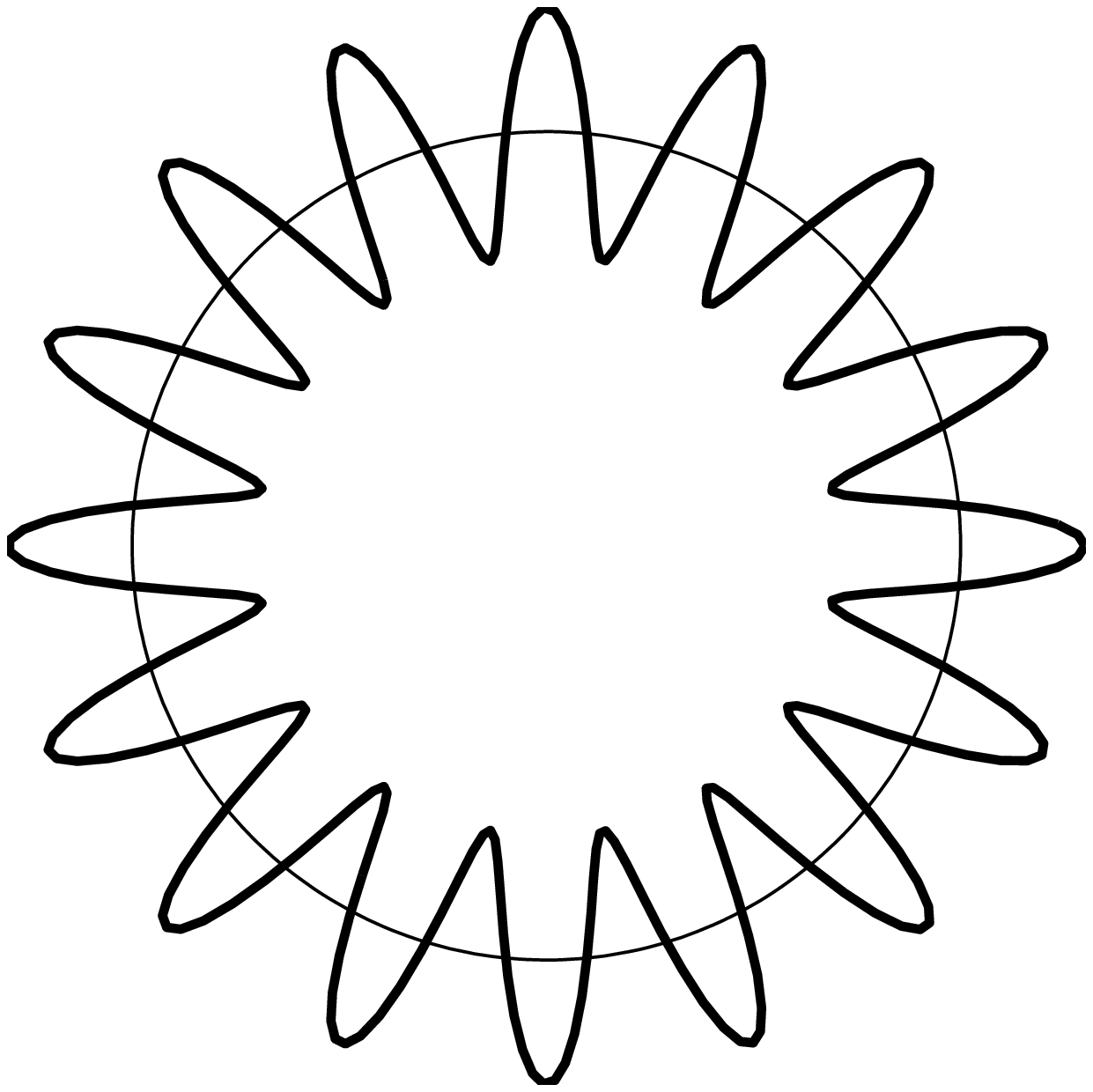}}
\put(75,90){\includegraphics[width=2.4cm]{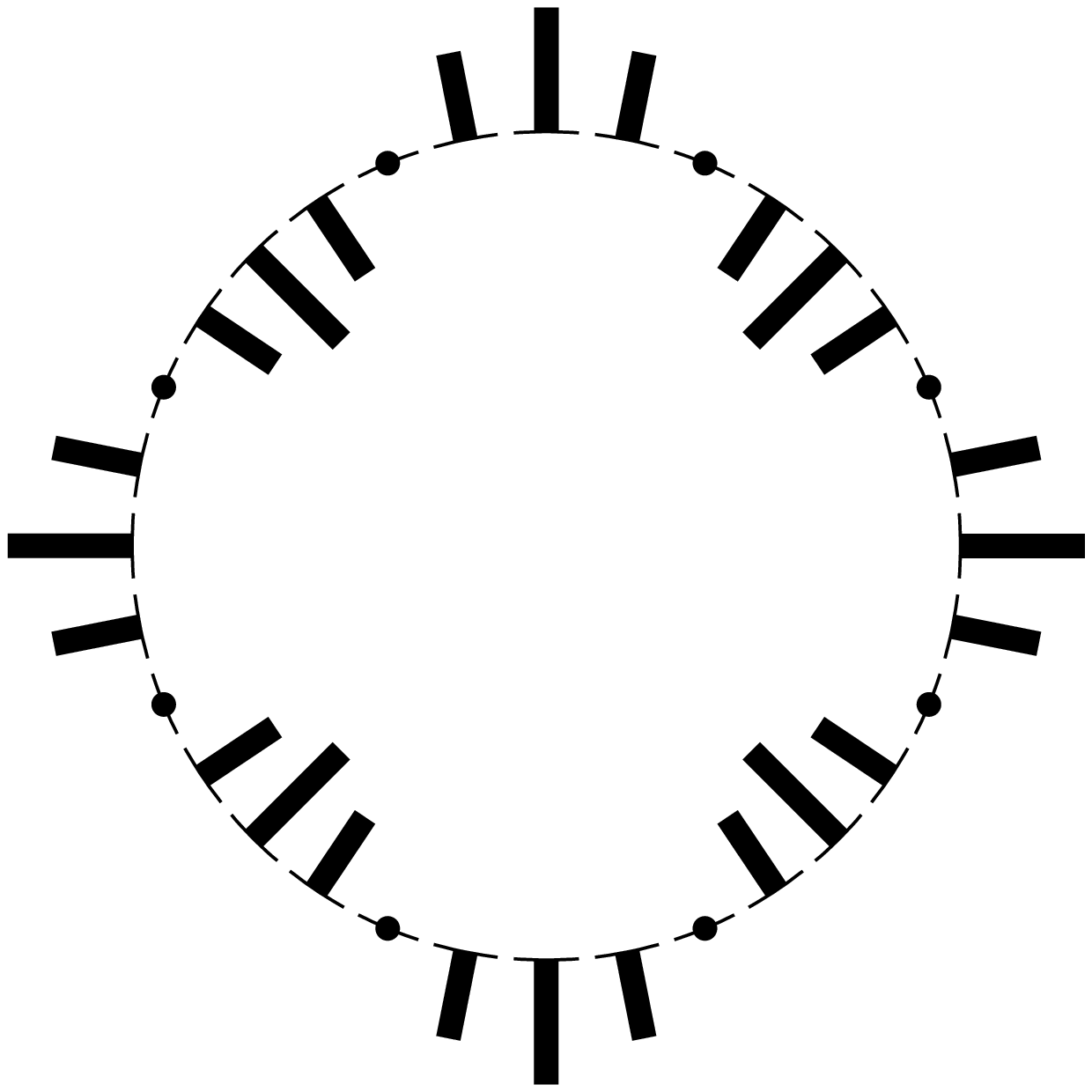}}
\put(75,0){\includegraphics[width=2.4cm]{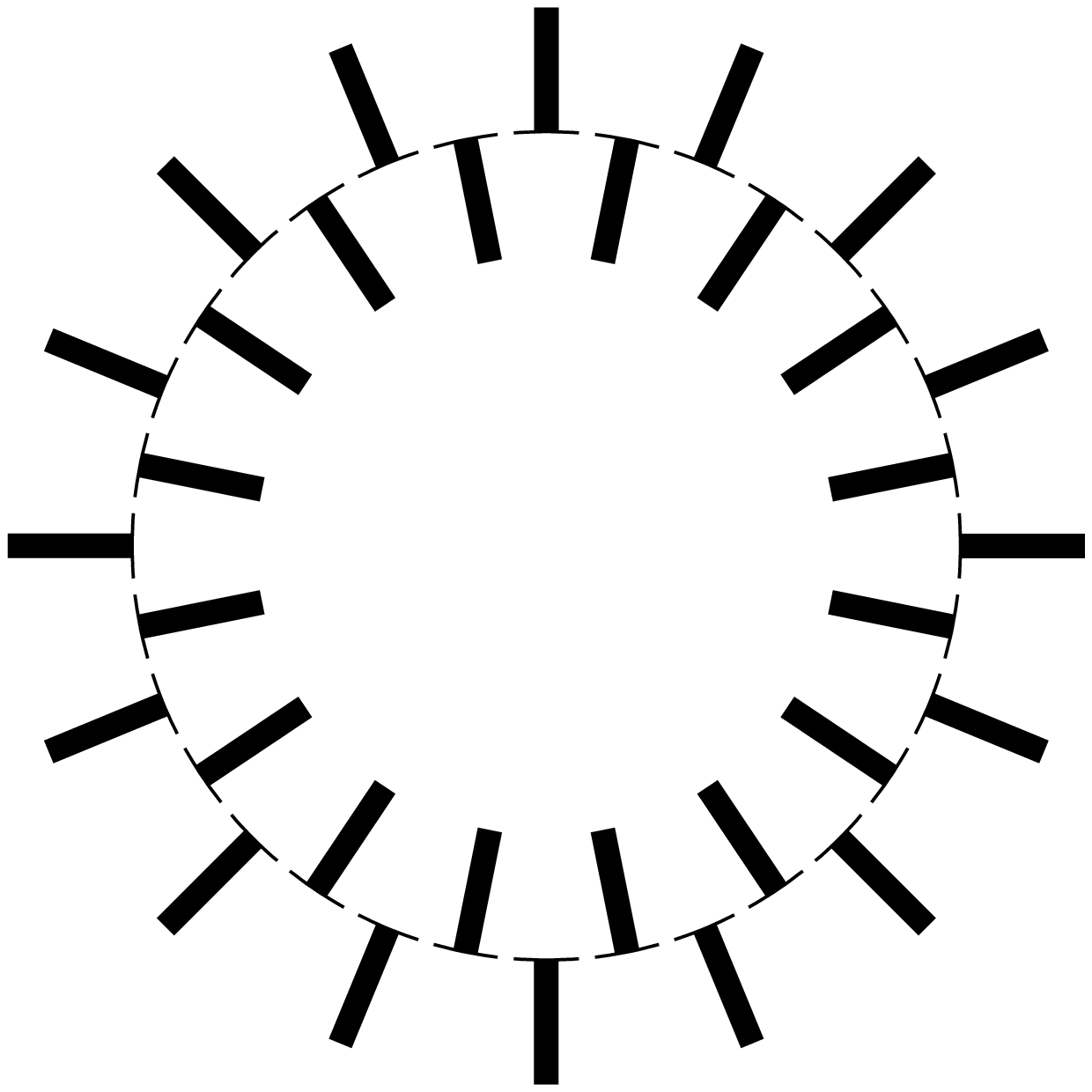}}
\put(165,90){\includegraphics[width=1.9cm]{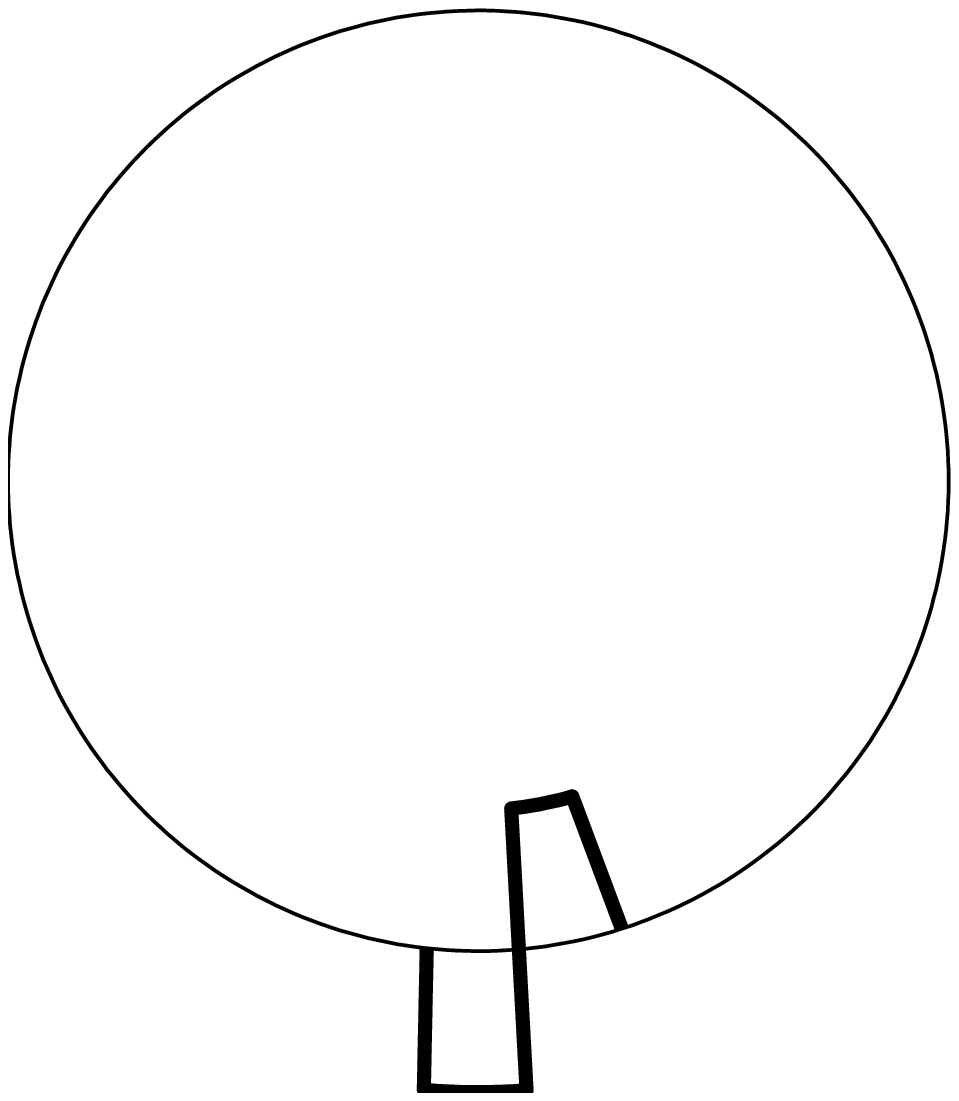}}
\put(165,0){\includegraphics[width=1.9cm]{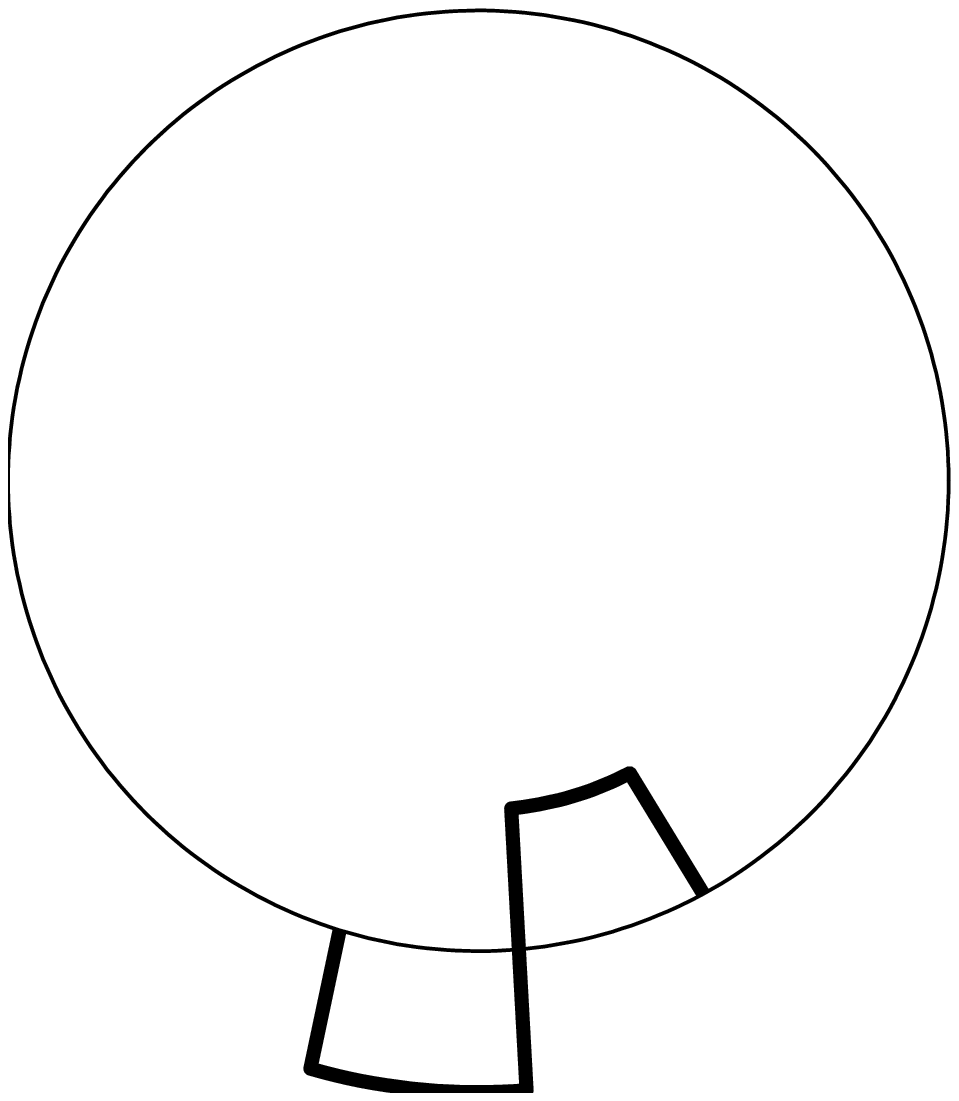}}
\put(245,90){\includegraphics[width=2cm]{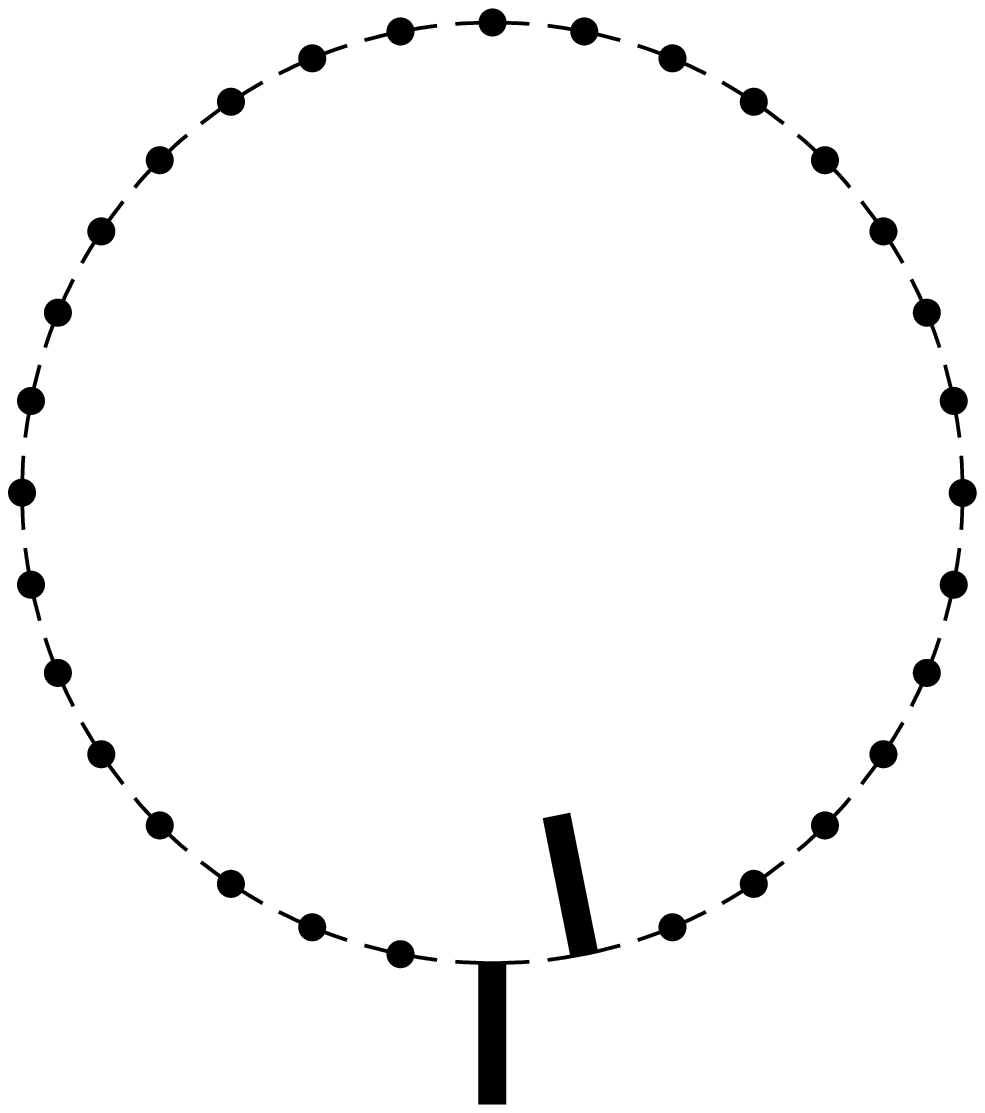}}
\put(245,0){\includegraphics[width=2cm]{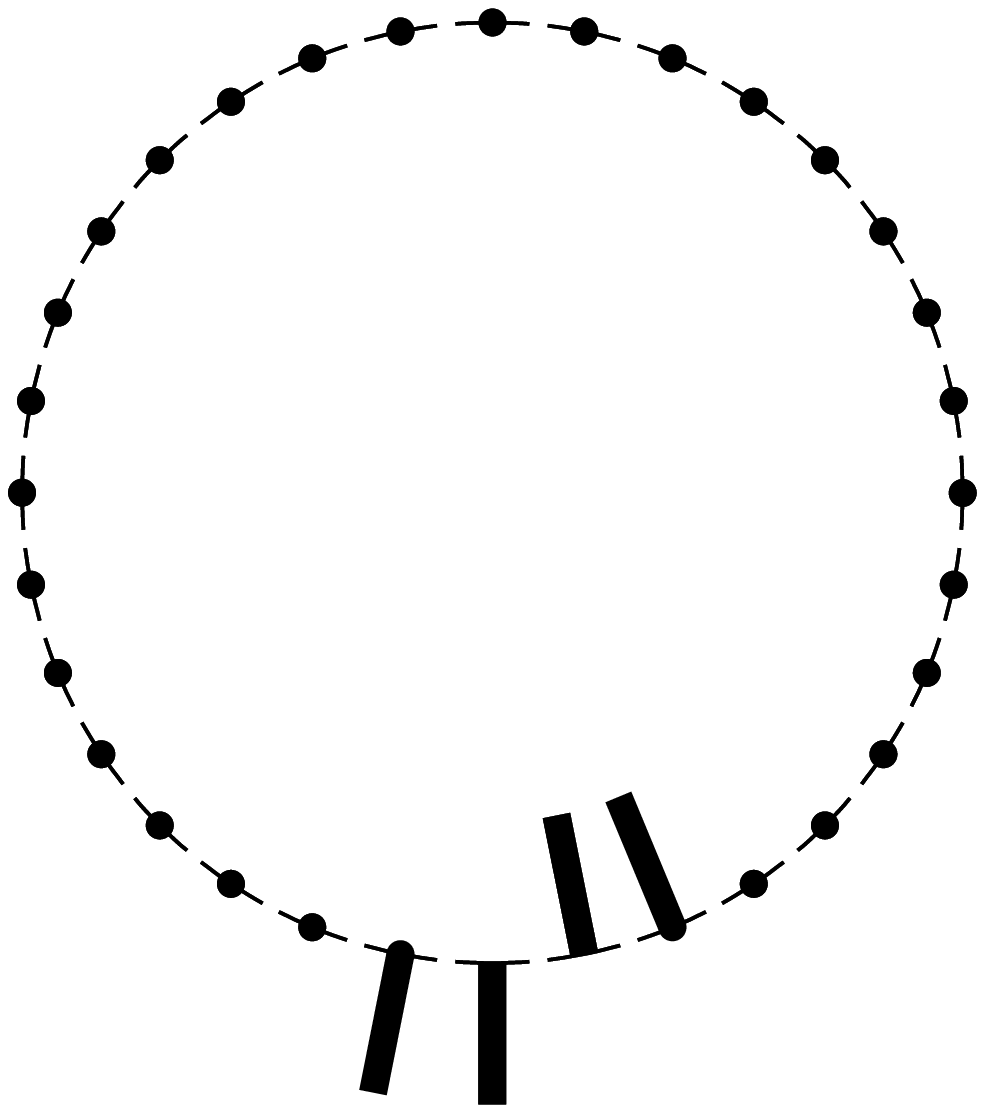}}
\put(0,205){\small Trigonometric basis functions and }
\put(0,193){\small corresponding electrode inputs:}
\put(0,181){\small for full-boundary data}
\put(168,205){\small Haar wavelet functions and }
\put(168,193){\small corresponding electrode inputs:}
\put(168,181){\small for both full and partial data}
 \put(150,-2){\line(0,1){220}}
\end{picture}
\caption{\label{Fig:trig_and_Haar}Illustration of various basis functions used for solving boundary integral equations. Also shown are voltage patterns applied using 32 electrodes, approximating the basis functions. The right-hand side functions and patterns are localized and therefore may be more suitable for working with partial-boundary data.}
\end{figure}

%.........................................................................................................................
How does one solve (\ref{partialBIE1.0}) numerically? Computational solution methods for boundary integral equations (BIEs) of type (\ref{fullBIE1.0}), corresponding to the continuum model, have most often been based on representing the unknown CGO solutions in terms of (generalized) trigonometric bases, where the basis functions are essentially supported on the entire boundary \cite{Knudsen2009,Astala2011,Mueller2012}. This approach is not directly applicable to partial data problems. In this work we present new numerical experiments on the unit disc (\emph{without loss of generality}) involving the solution of the above-mentioned BIEs using localized basis functions supported only on a subset of the boundary, in this case the Haar wavelets which are naturally applicable to the partial (as well as full) boundary continuum model. See Figure \ref{Fig:trig_and_Haar}.

Let us stress that at present there is no proof available for the solvability of equation \eqref{partialBIE1.0}. However, we did not encounter any problems when numerically solving \eqref{partialBIE1.0}, suggesting that it may be possible to prove unique solvability under appropriate assumptions.

We demonstrate that it is possible to recover the traces of the CGO solutions approximately on the part of the boundary available for measurements. In other words, the approximation in \eqref{mainapprox} is quite good in the $C^2$ and discontinuous conductivity examples we consider. In addition, we show below that these partial traces lead to interesting and useful reconstructions of practically relevant discontinuous conductivities.

Our new results may be useful in extending  three-dimensional D-bar reconstructions, such as \cite{Cornean2006,Boverman2008a,Bikowski2010a,Delbary2011}, to partial-boundary data.

We mention that numerical reconstructions using restricted information about the conductivity have been published in cases of partial-boundary data, see \cite{Mueller1999,Ide2007,Ide2010,Uhlmann2008}. The present work differs from those in that we aim to recover the full unknown conductivity function instead of inclusions in a known background. Also, there is an alternative methodology for partial-data EIT based on resistor networks \cite{Mamonov2010, Borcea2010,Borcea2010b}; our work again represents a very different approach. Finally, we mention that there is a large body of work on iterative solution methods for EIT in cases of partial data; those studies are fundamentally different from our direct (non-iterative) approach.

%+++++++++++++++++++++++++++++++
%+++++++++++++++++++++++++++++++

%--------------------------------------------------------------------
\section{Method 1: The method based on the Schr\"odinger equation}\label{sec:Method1} 
%--------------------------------------------------------------------
As mentioned above, the proof by \cite{Nachman1996} transforms the conductivity equation \eqref{Calderon:eq1pb} to the Schr\"{o}dinger equation 
\[\left(-\Delta + q(z)\right)\psi(z,k) = 0,\quad z\in\Omega\text{ and } k\in \C,\]
via the change of variables $\psi(z,k) = \sqrt{\sigma}u(z,k)$, where
\[q(z) = \frac{\Delta\sqrt{\sigma(z)}}{\sqrt{\sigma(z)}},\]
denotes the Schr\"{o}dinger potential and $\R^2$ is associated with $\C$ via $z=(x,y)=x+iy$.  Without loss of generality, the conductivity $\sigma$ is assumed to be $1$ near $\bndry$ and then extended to 1 in all of $\C$.  Existence and uniqueness are then studied for the well known  Schr\"{o}dinger equation with CGO solutions $\psi$ that are asymptotic to $e^{ikz}$ and $kz=(k_1+ik_2)(x+iy)$.  

The alternative Lippmann-Schwinger formulation 
\begin{equation}\label{eq-LS}
\mu(z,k) = 1-g_k\ast\left(q\mu\right), \quad z,\; k\in \C,
\end{equation}
of the Schr\"{o}dinger equation uses the related CGO solutions $\mu(z,k) = e^{-ikz}\psi(z,k)$ where $g_k$ is related to the Faddeev Green's function and defined in \eqref{def:Gk}. 

The reconstruction method of  Nachman \cite{Nachman1996} from infinite-precision data consists of the following two steps:
\[\Lambda_\sigma \overset{1}{\longrightarrow} \T(k)  \overset{2}{\longrightarrow} \sigma.\]
\begin{itemize}
\item[{\bf Step 1:}] {\bf From boundary measurements $\Lambda_\sigma$ to the scattering transform $\T$.}\\
For each fixed  $k\in\C$, solve in $H^{1/2}(\DOm)$ the integral equation
\begin{equation}\label{eq-Psi-INT} 
  \psi(z,k) = e^{ikz} - \int_{\bndry}G_k(z-\zeta)(\Lambda_\sigma-\Lambda_1)\psi(\zeta,k)\;dS(\zeta), \quad z\in\bndry,
\end{equation}
for the CGO solutions $\psi$ where the D-N map of the homogeneous conductivity $1$ is denoted by $\Lambda_1.$ Then, substitute $\psi$ into the formula for the nonlinear scattering transform $\T:\C\ra\C$:
\begin{align}\label{Tdef}
  \T(k)=\int_{\bndry}e^{i\bar{k}\bar{z}} (\Lambda_\sigma-\Lambda_1)\psi(z,k)\;dS(z),
\end{align}
where $dS$ denotes arclength measure on $\DOm$.
\vspace{1em}
  \item[{\bf Step 2:}] {\bf From the scattering transform $\T$ to the conductivity $\sigma$.}\\ Denote $e(z,k):=\exp(i(kz+\overline{k}\overline{z}))$. For each fixed $z\in \Om$, solve the integral equation
    \begin{align}\label{eq:int}
      \mu(z,k)=  1+\frac{1}{(2\pi)^2}\int_{\R^2}
      \frac{\T(k^\prime)}{(k-k^\prime)\bar{k}^\prime}e(-z,k^\prime)
      \overline{\mu(z,k^\prime)}dk^\prime_1 dk^\prime_2,
    \end{align}
then the conductivity is recovered by $\sigma(z)=\mu(z,0)^2$.
\end{itemize}
The integral equation \eqref{eq:int} was obtained from a corresponding partial differential equation, a so-called D-bar
equation,  which involves the
derivative with respect to the complex variable $\overline k$. This is where the D-bar method gets its name.
%--------------------------------------------------------------------
\section{Method 2: The method based on a $\dez$ and $\dbarz$ system}\label{sec:Method2}
%--------------------------------------------------------------------
As the existence/uniqueness result by Francini \cite{Francini2000} holds for complex admittivities $\gamma=\sigma+i\omega\epsilon$, where $\omega$ is the frequency of the applied current and $\epsilon$ denotes the electrical permittivity, and is an extension of that by Brown and Uhlmann \cite{Brown1997} we will formulate the problem in the complex case.   Let $u_1(z,k)$ and $u_2(z,k)$ be two CGO solutions to \eqref{Calderon:eq1pb} with asymptotic behavior $\frac{e^{ikz}}{ik}$ and $\frac{e^{ik\bar{z}}}{-ik}$, respectively.  Introduce a matrix $\Psi(z,k)$ of CGO solutions related to $u_1$ and $u_2$ by
\[\left(\begin{array}{c}
\Psi_{11}\\
\Psi_{21}\\
\end{array}\right) =\gamma^{1/2}\left(\begin{array}{c}
\dez u_1\\
\dbarz u_1\\
\end{array}\right), \quad\quad \left(\begin{array}{c}
\Psi_{12}\\
\Psi_{22}\\
\end{array}\right) =\gamma^{1/2}\left(\begin{array}{c}
\dez u_2\\
\dbarz u_2\\
\end{array}\right),\]
for  $z\in\Omega$ and  $k\in\C$.  The transformed system is then
\begin{equation}\label{eq-Dpsi-Qpsi}
D\Psi = Q\Psi,
\end{equation}
where $D$ is a matrix of $\dez$ and $\dbarz$ partial derivatives and $Q$ represents a matrix potential 
\begin{equation}
\label{matrixim}
\begin{split}
Q(z) = \begin{pmatrix} 0 & -\frac{1}{2}\partial_z\log\,\gamma(z) \\{-\frac{1}{2}\bar{\partial}_z\log \,\gamma(z)} & 0\end{pmatrix},
\qquad D = \begin{pmatrix} \dbarz &  0\\0 & \dez\end{pmatrix}.
\end{split}
\end{equation}
The admittivity $\gamma$ is assumed to be 1 near $\bndry$ and is extended to 1 in all of $\C$.  Existence and uniqueness of solutions are then studied for \eqref{eq-Dpsi-Qpsi} for $z\in\C$ instead of \eqref{Calderon:eq1pb}.  As it is more practical to work with CGOs with finite asymptotic behavior, we often make use of the related matrix of CGO solutions $M(z,k)\sim\left(\begin{array}{cc}
1&0\\
0&1
\end{array}\right)$ defined by

\begin{equation}
\label{eq-M-Psi}
M (z,k)=\Psi(z,k)\begin{pmatrix} e^{-izk} & 0\\0 & e^{i\bar{z}k}\end{pmatrix}=\begin{pmatrix} e^{-izk}\Psi_{11}(z,k) & e^{i\bar{z}k}\Psi_{12}(z,k) \\ e^{-izk}\Psi_{21}(z,k) & e^{i\bar{z}k}\Psi_{22}(z,k) \end{pmatrix}.
\end{equation}

Similarly to Method 1, the full data direct reconstruction algorithm \cite{Hamilton2012} also involves solving Fredholm integral equations of the second kind for CGO solutions using D-N data, evaluating a nonlinear scattering transform $S(k)$, solving a $\dbark$ equation, and using the recovered CGO solutions at $k=0$ to reconstruct the conductivity and permittivity.  The method can be summarized in the following steps:

\vspace{0.5em}
\[\Lambda_\gamma \overset{1}{\longrightarrow} S(k)  \overset{2}{\longrightarrow} M(z,0)   \overset{3}{\longrightarrow} \gamma.\]

\vspace{0.5em}
\begin{itemize}
\item[{\bf Step 1:}] {\bf From boundary measurements $\Lambda_\gamma$ to the scattering transform $S$.}\\
For fixed $k\in\C\setminus0$, solve Fredholm integral equations of the second kind on $\bndry$ for the traces of the CGO solutions $u_1(z,k)$ and $u_2(z,k)$:
\small{\begin{eqnarray}
\hspace{4em}u_1(z,k) &=&  \frac{e^{ik{z}}}{ik}- \int_{\bndry}G_k(z-\zeta)(\Lambda_\gamma - \Lambda_1) u_1(\zeta,k) dS(\zeta), \quad z\in\bndry \label{bie_u_1}\\
\hspace{4em}u_2(z,k) &=&  \frac{e^{-ik\bar{z}}}{-ik} - \int_{\bndry}G_k(-\bar{z}+\overline{\zeta})(\Lambda_\gamma - \Lambda_1) u_2(\zeta,k) dS(\zeta), \quad z\in\bndry. \label{bie_u_2}
\end{eqnarray}}\normalsize
Use the traces of $u_1$ and $u_2$ to compute the off diagonal entries of the {\sc CGO} solutions $\Psi(z,k)$ for $z\in\bndry$ from the BIEs
\begin{eqnarray}
\hspace{2em}\Psi_{12}(z,k) &=& \int_{\bndry}\frac{e^{i \bar{k}(z-\zeta)}}{4\pi(z-\zeta)}\left[\Lambda_\gamma-\Lambda_1\right]u_2(\zeta,k)\;dS(\zeta), \quad z\in\bndry \label{bie_psi_12}\\
\hspace{2em}\Psi_{21}(z,k) &=& \int_{\bndry}\overline{\left[\frac{e^{i k(z-\zeta)}}{4\pi(z-\zeta)}\right]}\left[\Lambda_\gamma-\Lambda_1\right]u_1(\zeta,k)\;dS(\zeta), \quad z\in\bndry, \label{bie_psi_21}
\end{eqnarray}
and compute the off-diagonal entries of the scattering matrix $S(k)$
\begin{eqnarray}
S_{12}(k) &=& \frac{i}{2\pi}\int_{\bndry}e^{-i\bar{k}z}\Psi_{12}(z,k) (\nu_1+i\nu_2)dS(z), \quad k\in\C \label{bie_S_12}\\
S_{21}(k) &=& -\frac{i}{2\pi}\int_{\bndry}e^{i\bar{k}\bar{z}}\Psi_{21}(z,k) (\nu_1-i\nu_2)dS(z), \quad k\in\C. \label{bie_S_21}
\end{eqnarray}
Interpolate the scattering data $S(k)$ to include $k=0$.
\vspace{1em}
\item[{\bf Step 2:}] {\bf From the scattering transform $S(k)$ to CGO solutions $M(z,0)$.}\\
Solve the $\dbar_k$ equation \eqref{eq-dbar-k} for the matrix $M(z,k)$
\begin{equation}\label{eq-dbar-k}
\dbar_k M(z,k)=M(z,\bar{k})\left(\begin{array}{cc}
e(z,\bar{k}) & 0\\
0 & e(z,-k)\\
\end{array}\right) S(k).
\end{equation}
\vspace{1em}
\item[{\bf Step 3:}] {\bf From CGO solutions $M(z,0)$ to the Admittivity $\gamma=\sigma + i\omega\epsilon$.}\\
Reconstruct the matrix potential $Q$ from 
\begin{equation}\label{eq-M-to-Q}
Q_{12}(z)=\frac{\dbarz M_+(z,0)}{M_-(z,0)},\quad\quad
Q_{21}(z)=\frac{\dez M_-(z,0)}{M_+(z,0)},
\end{equation}
where
\begin{eqnarray}
M_+(z,k)&=& M_{11}(z,k)+e^{-i(kz+\bar{k}\bar{z})}M_{12}(z,k)\label{eq-M+-thm}\\
M_-(z,k)&=& M_{22}(z,k)+e^{i(kz+\bar{k}\bar{z})}M_{21}(z,k)\label{eq-M--thm},
\end{eqnarray}
and use either $Q_{12}$ or $Q_{21}$ to recover $\gamma$
\begin{equation}\label{eq-Q-to-GAM}
\gamma(z)=\exp\left\{-\frac{2}{\pi}\int_{\Omega}\frac{Q_{12}(\zeta)}{\bar{z}-\bar{\zeta}}\;d\mu(\zeta) \right\}=\exp\left\{-\frac{2}{\pi}\int_{\Omega}\frac{Q_{21}(\zeta)}{z-\zeta}\;d\mu(\zeta)\right\},
\end{equation}
where the integration takes place over $\Omega$ rather than all of $\C$ due to the compact support of the matrix potential $Q$.
\end{itemize}
%--------------------------------------------------------------------
\section{Computation of Partial Boundary Data CGO Solutions}
%--------------------------------------------------------------------
In this work we use \emph{localized basis functions} in place of \emph{global basis functions}.  As mentioned above, the most commonly used global basis functions for the continuum model are the exponential trigonometric basis functions $e^{in\theta}$.  When electrode models (such as the gap, shunt, or complete electrode model) are used, a trigonometric basis of sines and cosines is often used instead \cite{DeAngelo_Mueller2010,Isaacson2004}.  As mentioned above, the common thread of these global basis functions is that their support is essentially the entire boundary of the domain.  By contrast, localized basis functions are supported on a subset of the boundary.  Examples of localized basis patterns include the skip patterns and adjacent patterns (see e.g., \cite{Hamilton_Thesis_2012,Ethan,HM12_NonCirc}) as well as the Haar wavelets.  In this work we use the Haar wavelets as they are localized basis functions that can be naturally used in both the continuum and electrode model cases.

As the boundary integral equations in Methods 1 and 2 are very similar, we will describe, without loss of generality, the computation in detail for Method 1.  In order to solve the boundary integral equation \eqref{eq-Psi-INT} 
\[\psi(z,k)= e^{ikz}-\int_{\bndry} G_k\left(z-\zeta\right)\left(\Lambda_\sigma-\Lambda_1\right)\psi(\zeta,k)\;dS(\zeta), \quad z\in\bndry,
\]
for the traces of the CGO solutions $\psi(z,k)$, we will need the Dirichlet-to-Neumann (D-N) map, and thus we must first discuss the applied voltage patterns, in this case, the Haar wavelets.
%.........................................................................................................................
\subsection{Description of Haar Wavelets}\label{sec-haar}
%.........................................................................................................................
Let $\Gamma$ denote a subset of the boundary $\bndry$ and let  $|\Gamma| = \mathcal{L}$ denote the length of the subset $\Gamma$.  
The first wavelet is the \emph{scaling function} which we will denote $\phi_1$, and is defined as:
\begin{equation}\label{eq-haar-phi-constant}
\begin{array}{cccc}
\phi_1(z) &=& h_1& z\in\Gamma\subseteq\bndry\\
h_1&=&\sqrt{\frac{1}{\mathcal{L}}}.
\end{array}
\end{equation}
If $z\in\bndry\setminus\Gamma$, the scaling function $\phi_1(z)$ is set to zero, as will be the case for the subsequent Haar wavelets.  

For ease of notation, let $d(z)$ be the distance, along the subset $\Gamma$ of the boundary, a point $z$ is from the beginning point $z_0$ on $\Gamma$ (corresponding to the smallest $\theta$ value in the traditional counter-clockwise orientation) and $z_{\mathcal{L}}$ the ending point (corresponding to the largest $\theta$ value).  Thus, $d=0$ at $z_0$ and $d=\mathcal{L}$ at $z_\mathcal{L}$.

The second wavelet is the so-called \emph{mother wavelet} which we will denote $\phi_2$ and is defined as:
\begin{equation}\label{eq-haar-phi-mother}
\phi_2(z) = \begin{cases}
h_1,& \left\{z\in\Gamma\middle| 0\leq d(z)<\frac{\mathcal{L}}{2}\right\}\\
-h_1,& \left\{z\in\Gamma\middle| \frac{\mathcal{L}}{2}\leq d(z)\leq\mathcal{L}\right\}.\\
\end{cases}
\end{equation}
As the basis functions need to be orthonormal, we require $\langle \phi_m, \phi_n\rangle=\delta_{m,n}$, which the above wavelets satisfy by construction.

The third and fourth Haar wavelets $\phi_3$ and $\phi_4$ are copies of the  mother wavelet $\phi_2$, squished into 1/2 the length of the support of $\phi_2$ as follows:
\begin{equation}\label{eq-haar-phi-3-4}
\begin{array}{cc}
\phi_3(z) =& \begin{cases}
h_2,& \left\{z\in\Gamma\middle| 0\leq d(z)<\frac{\mathcal{L}}{4}\right\}\\
-h_2,& \left\{z\in\Gamma\middle| \frac{\mathcal{L}}{4}\leq d(z)\leq\frac{\mathcal{L}}{2}\right\},\\
\end{cases}\\
&\\
\phi_4(z) =& \begin{cases}
h_2,& \left\{z\in\Gamma\middle| \frac{\mathcal{L}}{2}\leq d(z)<\frac{3\mathcal{L}}{4}\right\}\\
-h_2,& \left\{z\in\Gamma\middle| \frac{3\mathcal{L}}{4}\leq d(z)\leq\mathcal{L}\right\},\\
\end{cases}\\
&\\
&h_2= \sqrt{\frac{2}{\mathcal{L}}}.
\end{array}
\end{equation}
Notice that these new wavelets satisfy $\langle \phi_m, \phi_n\rangle=\delta_{m,n}$ for $m, n=1,\ldots,4$.

An exact formula for the $j$-th height function $h_j$ is
\[h_j = \sqrt{\frac{2^{j-1}}{\mathcal{L}}}, \quad j\geq1,\]
corresponding to Haar wavelets with support width $w_j$
\[w_j= \frac{\mathcal{L}}{2^{j-1}},\quad j\geq2.\]
%.........................................................................................................................
\subsection{Formation of the D-N map Using Haar Wavelets}
%.........................................................................................................................
As the main goal is to use only partial boundary data, thus applying and measuring data only on subset of the domain, it is more natural to apply voltages (rather than currents) and form the Dirichlet-to-Neumann (D-N) map directly.  Although in practice currents are frequently applied, and thus the Neumann-to-Dirichlet (N-D) map is formed first (which is done to dampen noise), that approach requires the inversion of the N-D map, which for partial data poses new questions.  As a preliminary approach, we proceed with Dirichlet data.

The conductivity equation \eqref{Calderon:eq1pb} can be solved using the Finite Element Method.  For each Haar wavelet, the Dirichlet boundary value problem is solved and the resulting solution $u$ in $\Omega$ is used to determine the current flux (Neumann data) at the boundary.  This allows the determination of the Neumann data  corresponding to the prescribed Dirichlet data and formation the D-N map:
\[\Lambda_\sigma f= \sigma\left.\frac{\partial u}{\partial\nu}\right|_{\bndry}.\]
Note that for the cases in this document, the conductivity on and near the boundary is 1.

The discrete matrix approximation to the D-N map is formed using the following formula for the $(m,n)$-th entry
\begin{equation}\label{eq-L-mat}
\Lambda_\sigma^M(m,n) := \langle \Lambda_\sigma \phi_m,\phi_n\rangle=\langle \sigma\nabla\phi_m\cdot\nu, \phi_n\rangle=\langle \nabla\phi_m\cdot\nu, \phi_n\rangle,
\end{equation}
where $\phi_j$ are the Haar wavelets described in Section~\ref{sec-haar} that now serve as the Dirichlet data and $\langle\cdot,\cdot\rangle$ denotes the $L^2$ inner product.  As $\nu$ denotes the outward facing unit normal and $\Omega$  is the unit disc, at the boundary point $z=e^{i\theta}=\cos(\theta)+i\sin(\theta)$, we have $\nu=\left(\cos(\theta),\sin(\theta)\right)$.  

%.........................................................................................................................
\subsection{Solution of the Full Data BIE}
%.........................................................................................................................
After forming the D-N map using the Haar wavelets, localized basis functions, we proceed  to solving the boundary integral equation for the traces of CGO solutions $\psi(z,k)$.  This involves the solution of a Fredholm integral equation of the second kind.  Following the approach of \cite{DeAngelo_Mueller2010,Hamilton2012,HM12_NonCirc,Hamilton_Thesis_2012}, we expand the exponential $e^{ikz}$ and CGOs $\psi(z,k)$ in the Haar wavelet patterns $\left\{\Phi_j\right\}_{j=1}^J$.  Let $z_\ell$ be an evaluation point on $\bndry$ and $J$ denote the number of linearly independent Haar wavelet functions used.  Then the values of the CGO solution $\psi$ and complex exponential $e^{ikz}$, for a given complex number $k$, at position $z_\ell$ on $\bndry$ are given by
\begin{equation} \label{eq-Psi-expansion}
\psi(z_\ell,k) \approx \sum_{j=1}^{J} b_j(k)\Phi^j_\ell, 
\end{equation}
and 
\begin{equation}\label{eq-exp-expansion}
e^{ikz_\ell} \approx \sum^{J}_{j = 1}c_j(k)\Phi^j_\ell,
\end{equation}
where $\Phi$ denotes the normalized Haar wavelets such that their $\ell^2$ norm is 1, i.e. they are related via
\[\Phi_j = \sqrt{\frac{\mathcal{L}}{L}}\phi_j.\]
Let ${\bf b}(k)$ denote the column vector ${\bf b}(k)=[b_1(k),\ldots,b_{J}(k)]^T$, and define ${\bf c}(k)$ analogously where $T$ denotes the standard matrix non-conjugate transpose. 

Let $E_{\ell^\prime}$ denote the $\ell^\prime$-th subdivision of the boundary $\bndry$ ($\ell^\prime=1,\ldots,L$) centered at the center of the $\ell^\prime$-th boundary element $z_{\ell^\prime}$ with length $2\pi/L$. Splitting the integral over $\bndry$ into a sum of integrals over the subsections $E_{\ell^\prime}$
\begin{eqnarray*}
\psi(z_\ell,k)
&\approx& e^{ikz_\ell}-\sum_{\ell^\prime=1}^L\int_{E_{\ell^\prime}} G_k\left(z_\ell-\zeta\right)\delta\Lambda_\sigma \psi(\zeta_{\ell^\prime},k)\;dS(\zeta)\\
&=& e^{ikz_\ell} -\sum_{\ell^\prime=1}^L\int_{E_{\ell^\prime}} G_k\left(z_\ell-\zeta\right)\;dS(\zeta)\;\left[\delta\Lambda_\sigma\psi(\zeta_{\ell^\prime},k)\right],
\end{eqnarray*}
where for ease of notation
\[\delta\Lambda_\sigma = \Lambda_\sigma - \Lambda_1.\]
Using the expansions for $\psi(z_\ell,k)$ and $e^{ikz_\ell}$, \eqref{eq-Psi-expansion} and \eqref{eq-exp-expansion} respectively, we have
\begin{eqnarray*}
\sum_{j=1}^{J}b_{j}(k)\Phi^j_\ell
&\approx&\sum_{j=1}^{J}c_{j}(k)\Phi^j_\ell- \sum_{\ell^\prime=1}^L\int_{E_{\ell^\prime}} G_k\left(z_\ell-\zeta\right)\;dS(\zeta)\;\left[\delta\Lambda_\sigma\sum_{j=1}^{J}b_{j}(k)\Phi^j_{\ell^\prime}\right]\\
&=&\sum_{j=1}^{J}c_{j}(k)\Phi^j_\ell- \sum_{\ell^\prime=1}^L\int_{E_{\ell^\prime}} G_k\left(z_\ell-\zeta\right)\;dS(\zeta)\sum_{j=1}^{J} b_{j}(k)f_j\left(\zeta_{\ell^\prime}\right),
\end{eqnarray*}
where $f_j\left(\zeta_{\ell^\prime}\right)$ denotes the action of the discretized $\delta\Lambda_\sigma^M$ matrix on the $j$-th normalized Haar wavelet basis function evaluated at $\zeta_{\ell^\prime}$.  Define the matrix approximation to the Faddeev Green's function as 
\begin{equation}\label{eq-Gk-cases}
{\mathbf{G}}_k(\ell,\ell^\prime) = \begin{cases}
G_k\left(z_\ell,\zeta_{\ell^\prime}\right) & \ell\neq\ell^\prime\\
0 & \ell=\ell^\prime,\\
%\frac{L}{\mathcal{L}}\int_{E_{\ell^\prime}} G_k(z_\ell-\zeta)\;dS(\zeta) & \ell=\ell^\prime,\\
\end{cases}
\end{equation}
removing the singularity at $G_k(0)$.  Then
\begin{equation}\label{eq-25-DeAngelo-u1}
\sum_{j=1}^{J}b_{j}(k)\Phi^j_\ell
\approx\sum_{j=1}^{J}c_{j}(k)\Phi^j_\ell- \frac{2\pi}{L}\sum_{j=1}^{J}b_{j}(k)\sum_{\ell^\prime=1}^L {\mathbf{G}}_k(\ell,\ell^\prime)  f_j\left(\zeta_{\ell^\prime}\right).
\end{equation}
Following \cite{DeAngelo_Mueller2010}
\begin{equation}\label{eq-DN-int-approx}
f_p(\zeta_{\ell^\prime})\approx \left(\Phi\delta\Lambda_\sigma^M\right)(\ell^\prime,j),
\end{equation}
i.e., the $(\ell^\prime,j)$ entry in the matrix resulting from multiplication of the  matrix  of normalized basis functions $\Phi$ and the discretized difference in D-N maps $\delta\Lambda_\sigma^M$.  Using the properties of matrix multiplication, equation~\eqref{eq-25-DeAngelo-u1} can be rewritten as
\[\sum_{j=1}^{J}b_{j}(k)\Phi^j_\ell =\sum_{j=1}^{J}c_{j}(k)\Phi^j_\ell - \frac{2\pi}{L}\sum_{j=1}^{J}b_{j}(k)\left({\mathbf{G}}_k \Phi\delta\Lambda_\sigma^M\right)(\ell,j),\]
or equivalently,
\[\Phi \mathbf{b}=\Phi\mathbf{c}-\frac{2\pi}{L}{\mathbf{G}}_k \Phi\delta\Lambda_\sigma^M\mathbf{b},\]
a matrix equation for the unknown coefficients $\mathbf{b}$ which are needed in the normalized Haar wavelet basis expansion of $\psi(z,k)$.

Using the orthonormality of the normalized Haar wavelet basis functions in the matrix $\Phi$, we multiply both sides of the equation by $\Phi^T$, and then solve
\begin{equation}\label{eq-28-DeAngelo-u1}
(I+A)\mathbf{b}=\mathbf{c},
\end{equation}
where
\begin{equation}\label{eq-DeAngelo-A-matrix}
A=\frac{2\pi}{L}\Phi^T{\mathbf{G}}_k \Phi\delta\Lambda_\sigma^M.
\end{equation}
To reiterate, for each desired value of $k\in\C$,  expand $e^{ikz}$ for $z\in\bndry$ in the normalized Haar wavelets $\Phi$ to define the vector of coefficients $\mathbf{c}$, and solve the system \eqref{eq-28-DeAngelo-u1} using GMRES for the unknown coefficients $\mathbf{b}$.  These coefficients are then used to reconstruct $\psi(z,k)$ for the specified value of $k$ via \eqref{eq-Psi-expansion}.  
%.........................................................................................................................
\subsection{Solution of the Partial Data BIE}
%.........................................................................................................................
We now proceed to the problem of interest, namely, the solution of the boundary integral equation \eqref{eq-Psi-INT} when only part of the boundary is accessible for data acquisition.    The solution method is nearly identical to the full data Haar wavelet case presented above.

Now $\Gamma$ is a proper subset of the boundary $\bndry$ and the Haar wavelets and D-N map are formed as above. The D-N map now corresponds to data taken only on the proper subset $\Gamma$ since the applied voltage is 0 off $\Gamma$.  

Let $\tilde{z}$ denote the boundary values $z$ restricted to $\Gamma$ and $\tilde{\psi}$ the corresponding partial data CGO solutions.  We then expand $e^{ik\tilde{z}}$ and $\tilde{\psi}$ as before and solve the resulting system for each desired value of $k\in\C$:
\[(I+\tilde{A})\mathbf{\tilde{b}}=\mathbf{\tilde{c}},\]
where
\begin{equation}\label{eq-DeAngelo-A-matrix-G0-PARTIAL}
\tilde{A}=\frac{\mathcal{L}}{L} \tilde{\Phi}^T\tilde{\mathbf{G}}_k\tilde{\Phi}\delta\tilde{\Lambda}_\sigma^M,
\end{equation}
and $\mathcal{L} = |\Gamma|$.
%--------------------------------------------------------------------
\section{Numerical Reconstruction of Conductivities from Partial Data Using Method~2}\label{sec-method2-pdata}
%--------------------------------------------------------------------
As stated above, the boundary integral equations \eqref{bie_u_1} and \eqref{bie_u_2} are nearly identical to the BIE \eqref{eq-Psi-INT} described above for Method~1.  Thus, their traces on the subset $\Gamma$ of the boundary can be recovered by solving analogous formulas.  The coefficients $\mathbf{\tilde{b}^1}$ for $u_1$ are determined by solving
\[(I+\tilde{A})\mathbf{\tilde{b}^1}=\mathbf{\tilde{c}^1},\]
where $\mathbf{\tilde{c}^1}$ are the coefficients in the Haar expansion of $\frac{e^{ik\tilde{z}}}{ik}$.  Similarly, the coefficients $\mathbf{\tilde{b}^2}$ for $u_2$ are determined by solving
\[(I+\tilde{A_2})\mathbf{\tilde{b}^2}=\mathbf{\tilde{c}^2},\]
where $\mathbf{\tilde{c}^2}$ are the coefficients in the Haar expansion of $\frac{e^{-ik\overline{\tilde{z}}}}{-ik}$ and $\tilde{A}_2$ now contains the matrix approximation of $G_k(-\bar{z}+\bar{\zeta})$ instead of $G_k(z-\zeta)$.

A natural question is whether these CGO solutions, which match very well on $\Gamma$ with their full data counterparts (see Section~\ref{sec-Test1-numerics}), can be used to produce informative reconstructions of the conductivity (and/or permittivity) near the region of the accessible boundary.  Our aim was to understand the extent of the impact of the partial data CGO solutions on the remainder of a D-bar algorithm.  Therefore, as an initial test, we left the remainder of the algorithm for Method~2 intact which  means computing the intermediate CGO solutions $\Psi_{12}$ and $\Psi_{21}$ using the partial D-N map and the partial data CGO solutions $u_1$ and $u_2$, computing the scattering transforms $S_{12}$ and $S_{21}$ over $\Gamma$, and proceeding with Steps~2-3 as before.  The steps of the proposed partial data algorithm are included here for the reader's convenience.

\vspace{0.5em}
\[\tilde{\Lambda}_\gamma \overset{1}{\longrightarrow} \tilde{S}(k)  \overset{2}{\longrightarrow} \tilde{M}(z,0)   \overset{3}{\longrightarrow} \tilde{\gamma}.\]

\vspace{0.5em}
\begin{itemize}
\item[{\bf Step 1:}] {\bf From \emph{partial} boundary measurements $\tilde{\Lambda}_\gamma$ to the approx. scattering transform $\tilde{S}(k)$.}\\
For fixed $k\in\C\setminus0$ such that $|k|<R$ a fixed radius depending on the measured D-N map, solve Fredholm integral equations of the second kind on $\Gamma\subset\bndry$ for the approximate traces of the CGO solutions $\tilde{u}_1(z,k)$ and $\tilde{u}_2(z,k)$ on $\Gamma\subset\bndry$
\begin{eqnarray}
\hspace{4em}\tilde{u}_1(z,k) &=&  \left.\frac{e^{ik{z}}}{ik}\right\vert_{\Gamma} - \int_{\Gamma}G_k(z-\zeta)(\tilde{\Lambda}_\gamma - \tilde{\Lambda}_1) \tilde{u}_1(\zeta,k) dS(\zeta), \quad z\in\Gamma\label{bie_u_1pb}\\
\hspace{4em}\tilde{u}_2(z,k)&=& \left. \frac{e^{-ik\bar{z}}}{-ik}\right\vert_{\Gamma} - \int_{\Gamma}G_k(-\bar{z}+\overline{\zeta})(\tilde{\Lambda}_\gamma - \tilde{\Lambda}_1) \tilde{u}_2(\zeta,k) dS(\zeta), \quad z\in\Gamma. \label{bie_u_2pb}
\end{eqnarray}
Use the approximate traces of $u_1$ and $u_2$ (namely, $\tilde{u}_1$ and $\tilde{u}_2$) to compute the approximate off diagonal entries of the {\sc CGO} solutions $\widetilde{\Psi}(z,k)$ for $z\in\Gamma\subset\bndry$ from the BIEs
\begin{eqnarray}
\hspace{2em}\widetilde{\Psi}_{12}(z,k) &=& \int_{\Gamma}\frac{e^{i \bar{k}(z-\zeta)}}{4\pi(z-\zeta)}\left[\tilde{\Lambda}_\gamma- \tilde{\Lambda}_1\right]\tilde{u}_2(\zeta,k)\;dS(\zeta), \quad z\in\Gamma \label{bie_psi_12pb}\\
\hspace{2em}\widetilde{\Psi}_{21}(z,k) &=& \int_{\Gamma}\overline{\left[\frac{e^{i k(z-\zeta)}}{4\pi(z-\zeta)}\right]}\left[\tilde{\Lambda}_\gamma- \tilde{\Lambda}_1\right]\tilde{u}_1(\zeta,k)\;dS(\zeta), \quad z\in\Gamma, \label{bie_psi_21pb}
\end{eqnarray}
and compute the off-diagonal entries of the scattering matrix $\tilde{S}(k)$ integrating over $\Gamma\subset\bndry$
\begin{eqnarray}
\tilde{S}_{12}(k) &=& \frac{i}{2\pi}\int_{\Gamma}e^{-i\bar{k}z}\widetilde{\Psi}_{12}(z,k) (\nu_1+i\nu_2)dS(z) \label{bie_S_12pb}\\
\tilde{S}_{21}(k) &=& -\frac{i}{2\pi}\int_{\Gamma}e^{i\bar{k}\bar{z}}\widetilde{\Psi}_{21}(z,k) (\nu_1-i\nu_2)dS(z). \label{bie_S_21pb}
\end{eqnarray}
Interpolate the approximate scattering data $\tilde{S}(k)$ to include $k=0$.
\vspace{1em}
\item[{\bf Step 2:}] {\bf From the approx. scattering transform $\tilde{S}(k)$ to approx. CGO solutions $\widetilde{M}(z,0)$.}\\
Solve the $\dbar_k$ equation \eqref{eq-dbar-k} for the matrix of approx. CGO solutions $\widetilde{M}(z,k)$
\begin{equation}\label{eq-dbar-k-pb}
\dbar_k \widetilde{M}(z,k)=\widetilde{M}(z,\bar{k})\left(\begin{array}{cc}
e(z,\bar{k}) & 0\\
0 & e(z,-k)\\
\end{array}\right) \tilde{S}(k).
\end{equation}
\vspace{1em}
\item[{\bf Step 3:}] {\bf From the approximate CGO solutions $\widetilde{M}(z,0)$ to approx. admittivity $\tilde{\gamma}=\tilde{\sigma} + i\omega\tilde{\epsilon}$.}\\
Reconstruct the approximate matrix potential $\tilde{Q}$ from 
\begin{equation}\label{eq-M-to-Q-pb}
\tilde{Q}_{12}(z)=\frac{\dbarz \widetilde{M}_+(z,0)}{\widetilde{M}_-(z,0)},\quad\quad
\tilde{Q}_{21}(z)=\frac{\dez \widetilde{M}_-(z,0)}{\widetilde{M}_+(z,0)},
\end{equation}
where
\begin{eqnarray}
\widetilde{M}_+(z,k)&=& \widetilde{M}_{11}(z,k)+e^{-i(kz+\bar{k}\bar{z})}\widetilde{M}_{12}(z,k)\label{eq-M+-thm-pb}\\
\widetilde{M}_-(z,k)&=& \widetilde{M}_{22}(z,k)+e^{i(kz+\bar{k}\bar{z})}\widetilde{M}_{21}(z,k)\label{eq-M--thm-pb},
\end{eqnarray}
and use either $\tilde{Q}_{12}$ or $\tilde{Q}_{21}$ to recover the approximation $\tilde{\gamma}$
\begin{equation}\label{eq-Q-to-GAM-pb}
\tilde{\gamma}(z)=\exp\left\{-\frac{2}{\pi}\int_{\Omega}\frac{\tilde{Q}_{12}(\zeta)}{\bar{z}-\bar{\zeta}}\;d\mu(\zeta) \right\}=\exp\left\{-\frac{2}{\pi}\int_{\Omega}\frac{\tilde{Q}_{21}(\zeta)}{z-\zeta}\;d\mu(\zeta)\right\},
\end{equation}
where the integration takes place over $\Omega$ rather than all of $\C$ due to the compact support of the matrix potential $\tilde{Q}$.
\end{itemize}

For the numerical details regarding how to implement Steps~2-3 see \cite{Hamilton2012,Hamilton_Thesis_2012,HM12_NonCirc}.
%--------------------------------------------------------------------
\section{Computational Experiments}
%--------------------------------------------------------------------
We considered two test problems.  Test~1 aims to determine how a partial data D-N map affects the values of the traces of the CGO solutions on the accessible portion of the boundary.  Test~2 aims to determine the effect of the partial data CGO solutions on a D-bar algorithm.
%.........................................................................................................................
\subsection{Test~1: Partial Data Traces of CGO Solutions}\label{sec-Test1-numerics}
\begin{figure}[h!]
\centering
\includegraphics[width=2.5in]{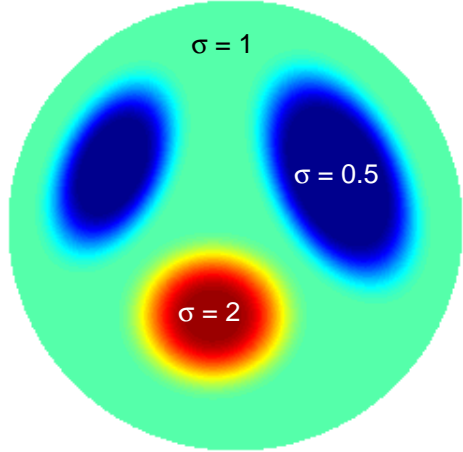}
\caption{The $C^2$ conductivity used in Test 1.}\label{fig-sigma-phantom}
\end{figure}

%.........................................................................................................................
For the first test problem we considered the $C^2$ smooth conductivity given in Figure~\ref{fig-sigma-phantom}.  The conductivity equation was first solved using the Finite Element method with 256 Haar wavelets with essential support on the entire boundary, serving as 256 different Dirichlet boundary conditions.  We considered the $3/4$, $1/2$, and $1/4$ data problems with 192, 128, and 64 Haar wavelets respectively.  Each of the partial data cases is centered around $z=1$, i.e. $\theta=0$.  

Using Method~1, we solved the full data matrix formulation \eqref{eq-28-DeAngelo-u1} of the boundary integral equation for the traces of the CGO solutions $\psi$.  The partial data traces of the CGO solutions were recovered by solving \eqref{eq-DeAngelo-A-matrix-G0-PARTIAL}.  In order to evaluate how well the reconstructed traces compare to the true traces we also solved the Lippmann-Schwinger equation \eqref{eq-LS} using the twice-differentiable conductivity in Figure~\ref{fig-sigma-phantom}.

\begin{figure}
\hspace{2em}
\begin{picture}(300,400)
\put(15,380){\textbf{Real parts of $\psi$}}
\put(170,380){\textbf{Imaginary parts of $\psi$}}
\put(-30,340){True}
\put(15,320){\includegraphics[width=4.75cm]{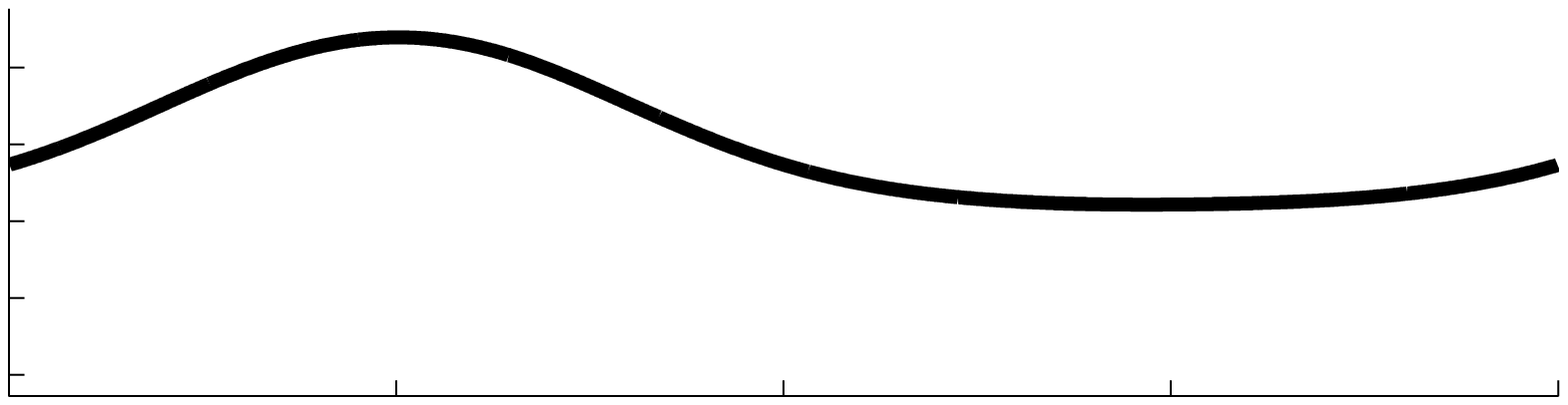}}
\put(170,320){\includegraphics[width=4.75cm]{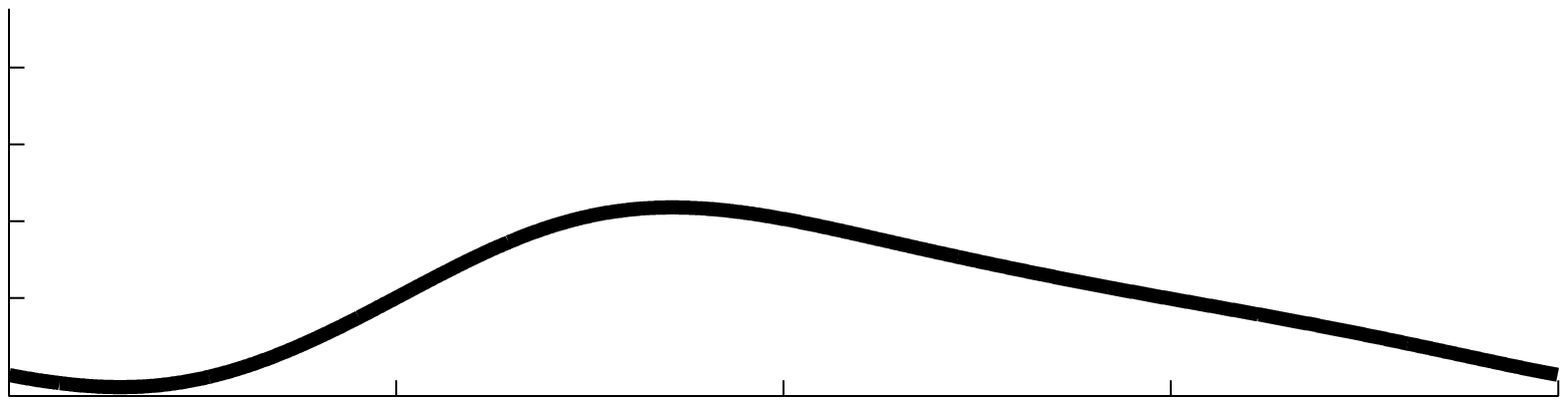}}
\put(-30,280){Full data}
\put(15,260){\includegraphics[width=4.75cm]{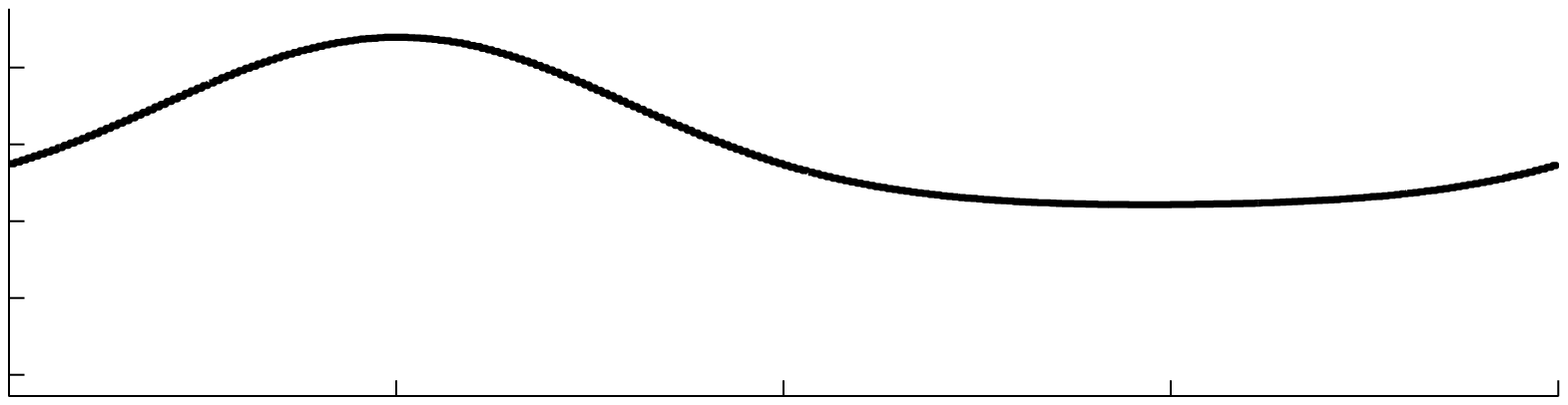}}
\put(170,260){\includegraphics[width=4.75cm]{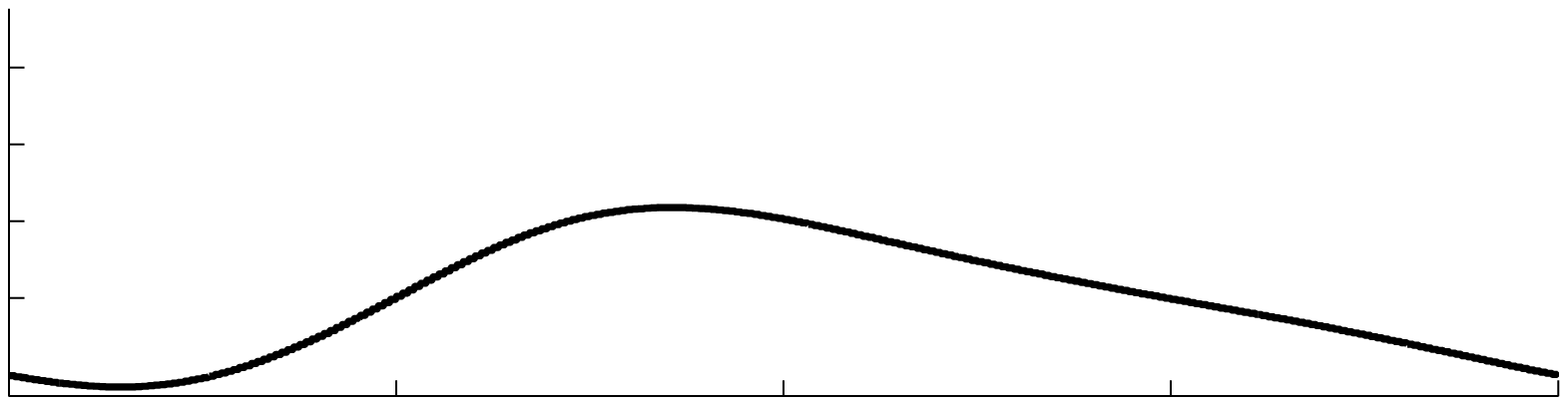}}
\put(-30,220){$\frac 34$ data}
\put(15,200){\includegraphics[width=4.75cm]{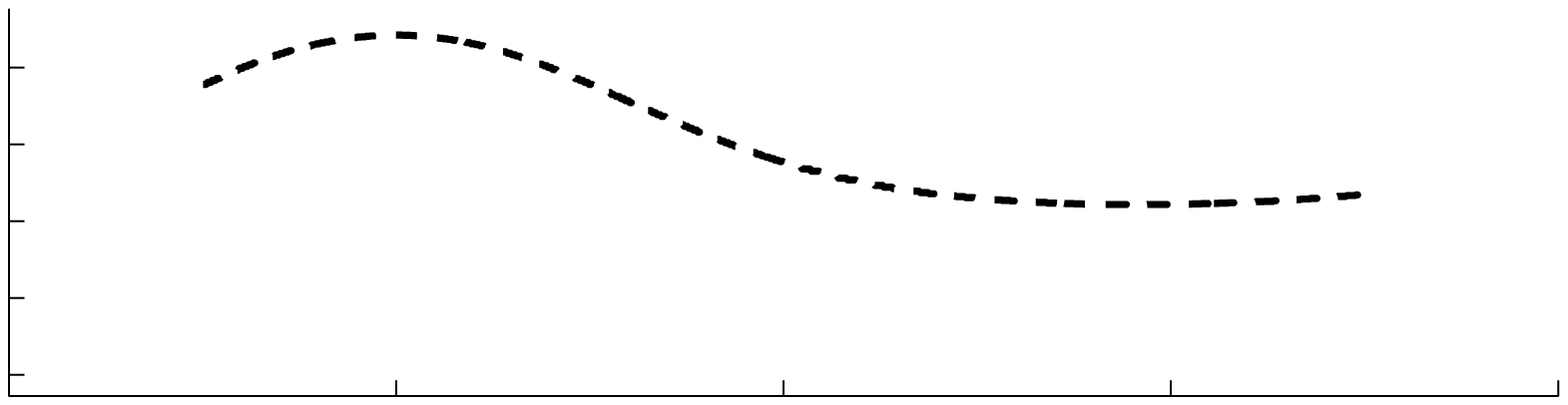}}
\put(170,200){\includegraphics[width=4.75cm]{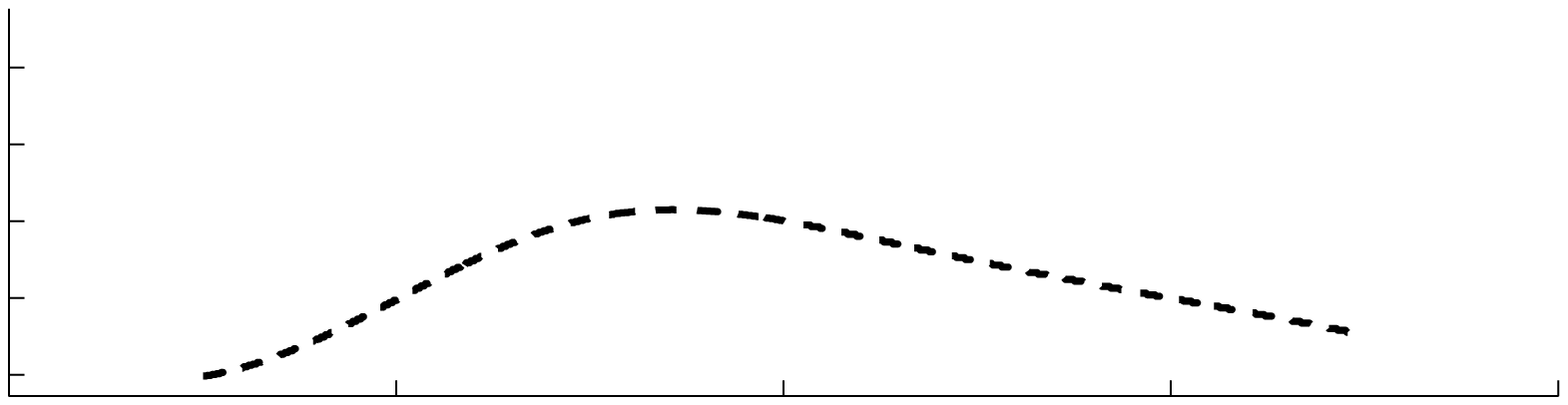}}
\put(-30,160){$\frac 12$ data}
\put(15,140){\includegraphics[width=4.75cm]{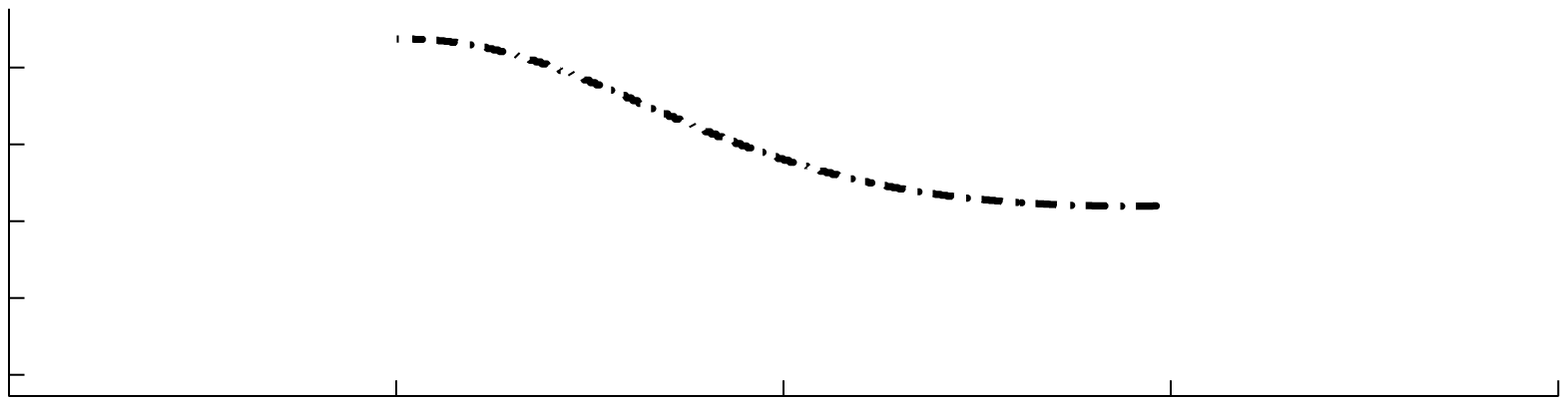}}
\put(170,140){\includegraphics[width=4.75cm]{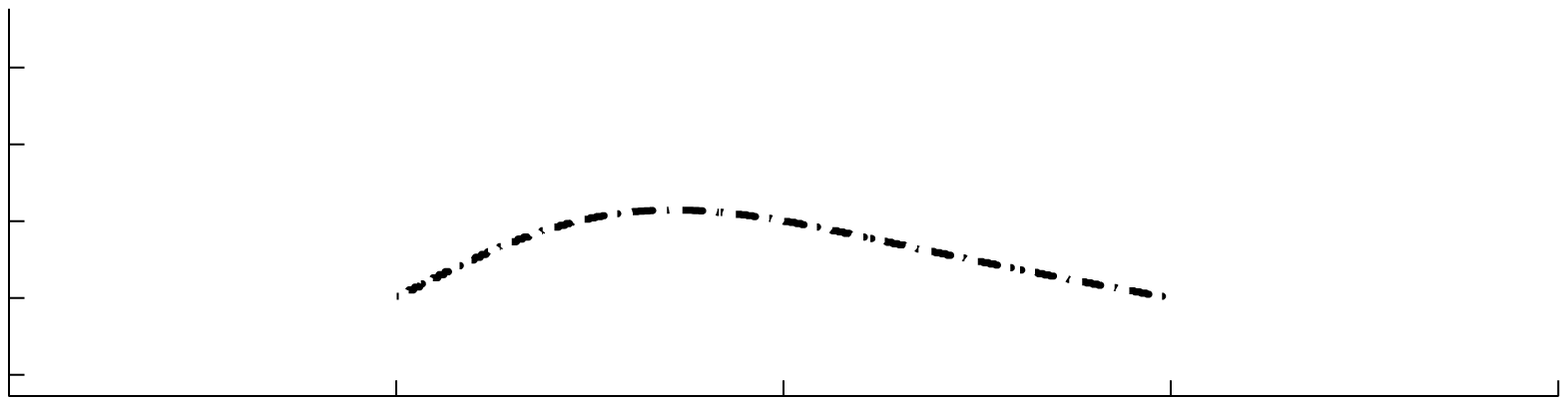}}
\put(-30,100){$\frac 14$ data}
\put(15,80){\includegraphics[width=4.75cm]{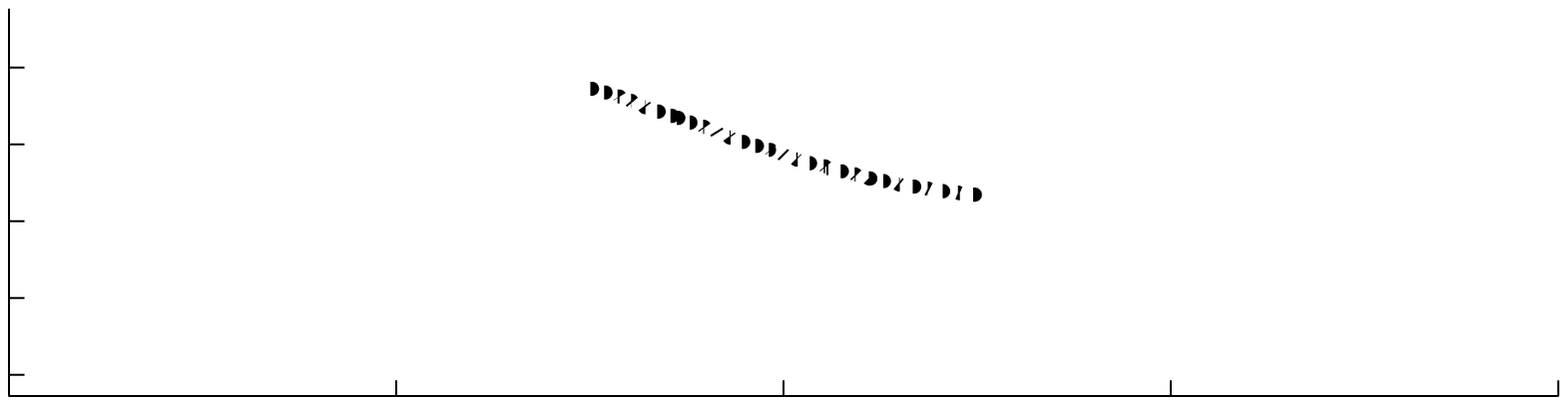}}
\put(170,80){\includegraphics[width=4.75cm]{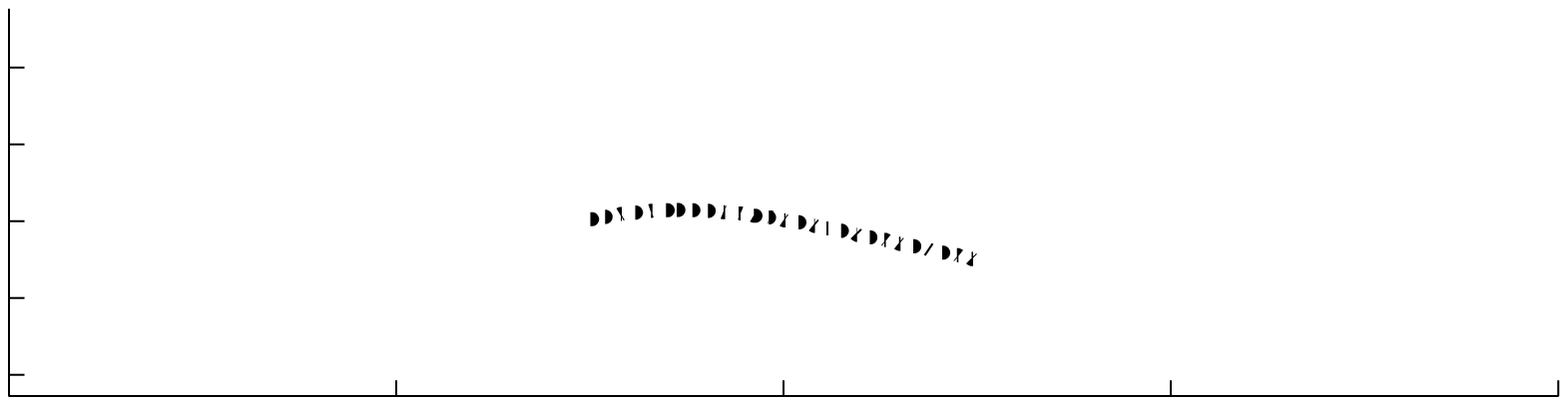}}
\put(-30,40){All}
\put(15,20){\includegraphics[width=4.75cm]{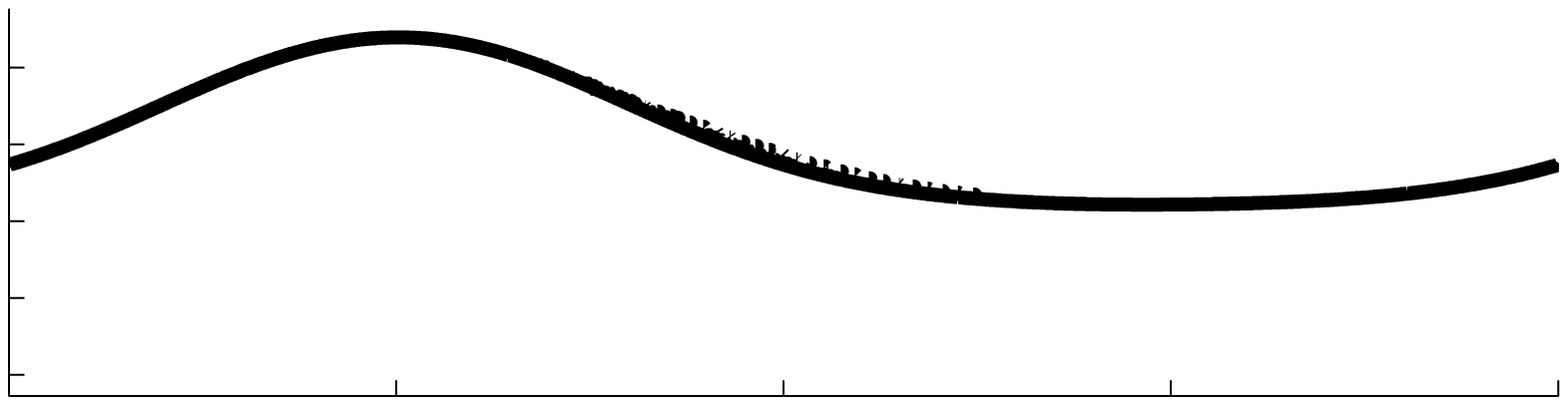}}
\put(170,20){\includegraphics[width=4.75cm]{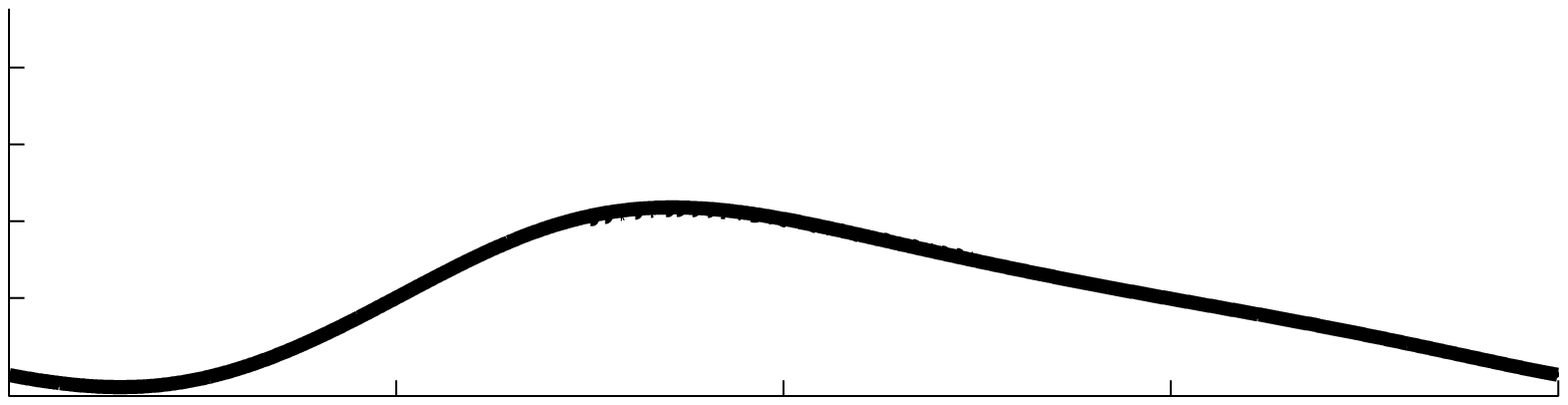}}
\put(10,10){\small$-\pi$}
\put(38,10){\small$-\pi/2$}
\put(81,10){\small$0$}
\put(110,10){\small$\pi/2$}
\put(147,10){\small$\pi$}

\put(164,10){\small$-\pi$}
\put(192,10){\small$-\pi/2$}
\put(236,10){\small$0$}
\put(265,10){\small$\pi/2$}
\put(302,10){\small$\pi$}
\end{picture}
\caption{\label{fig:tracesA}Traces of the CGO solutions $\psi$ corresponding to the $C^2$ conductivity in Figure~\ref{fig-sigma-phantom}. Here $k=0.5$.}
\end{figure}

\begin{figure}
\hspace{2em}
\begin{picture}(300,400)
\put(15,380){\textbf{Real parts of $\psi$}}
\put(170,380){\textbf{Imaginary parts of $\psi$}}
\put(-30,340){True}
\put(15,320){\includegraphics[width=4.75cm]{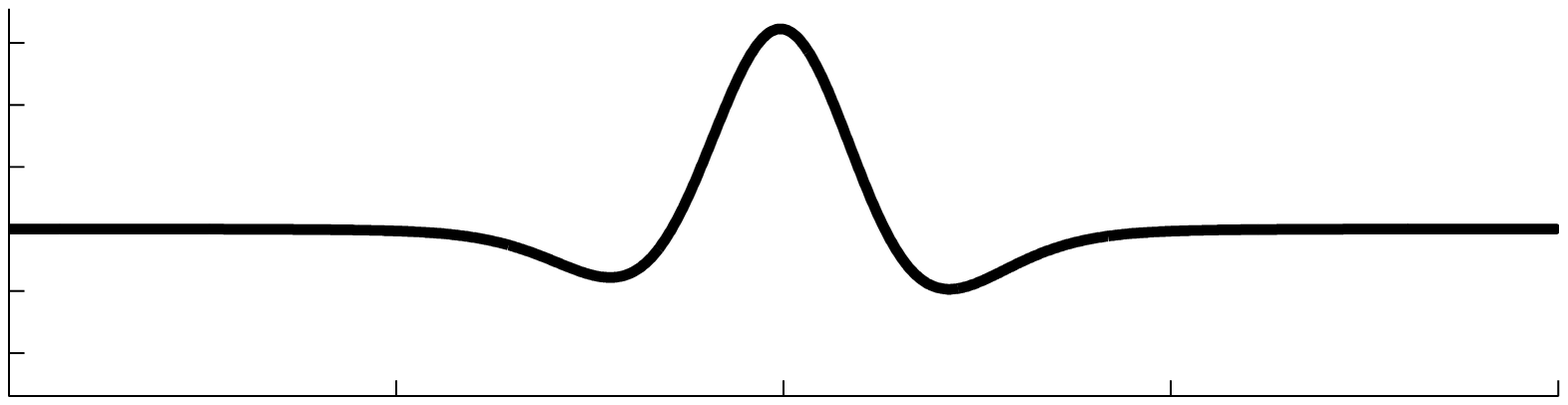}}
\put(170,320){\includegraphics[width=4.75cm]{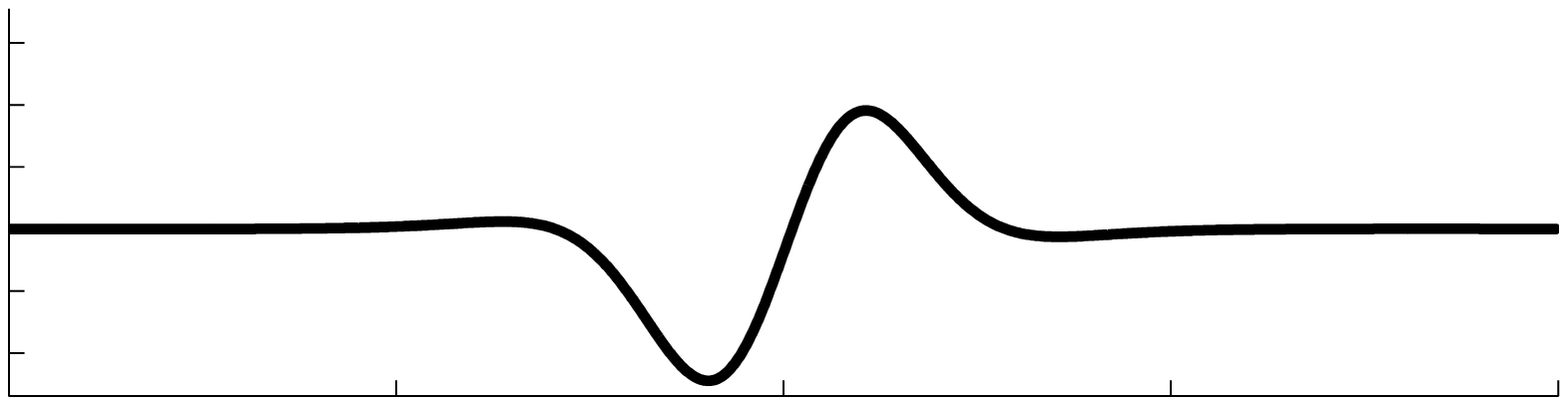}}
\put(-30,280){Full data}
\put(15,260){\includegraphics[width=4.75cm]{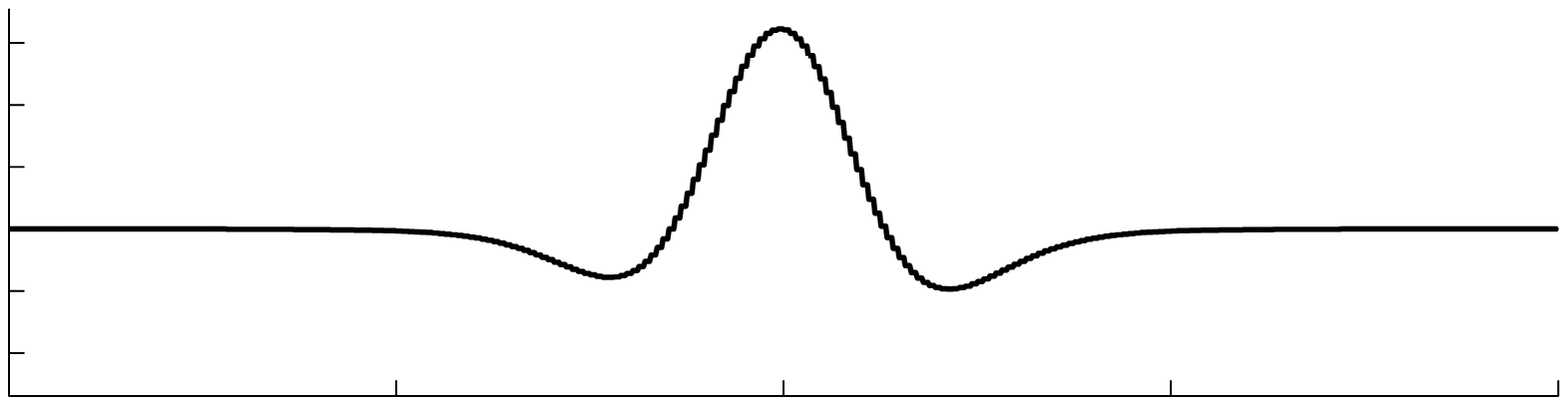}}
\put(170,260){\includegraphics[width=4.75cm]{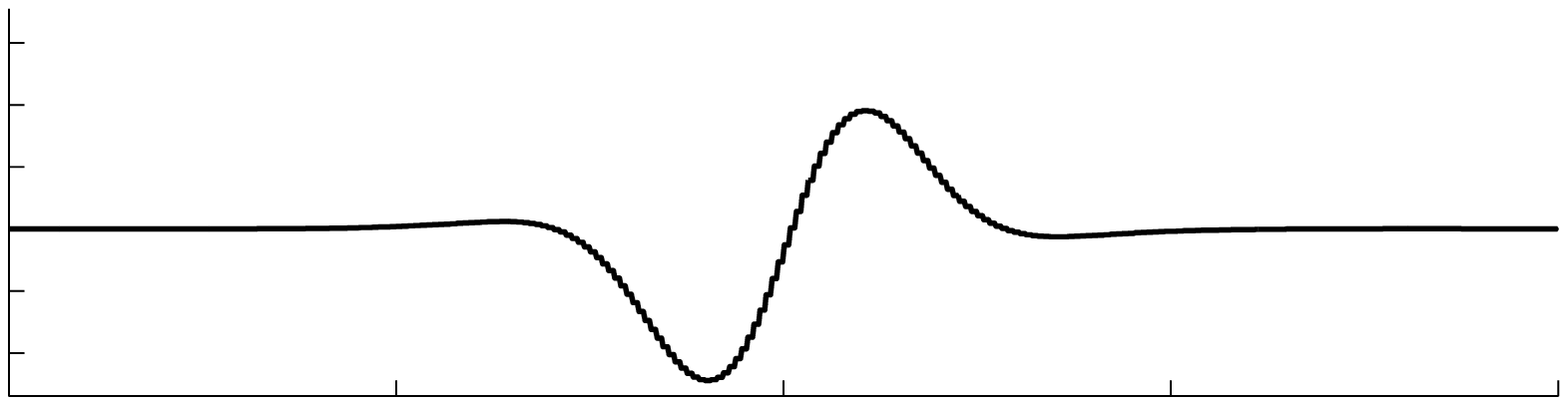}}
\put(-30,220){$\frac 34$ data}
\put(15,200){\includegraphics[width=4.75cm]{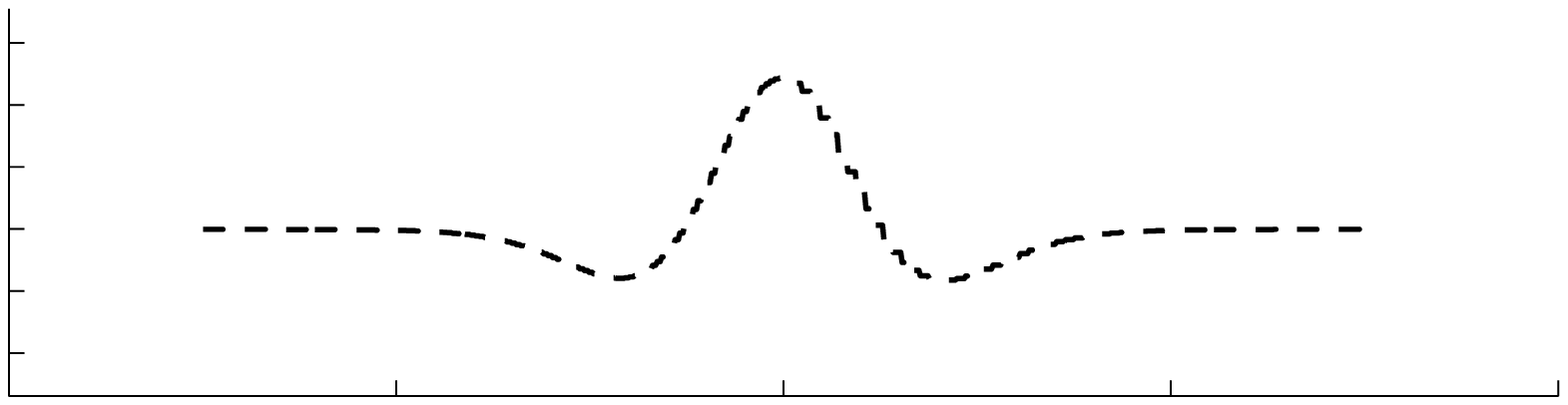}}
\put(170,200){\includegraphics[width=4.75cm]{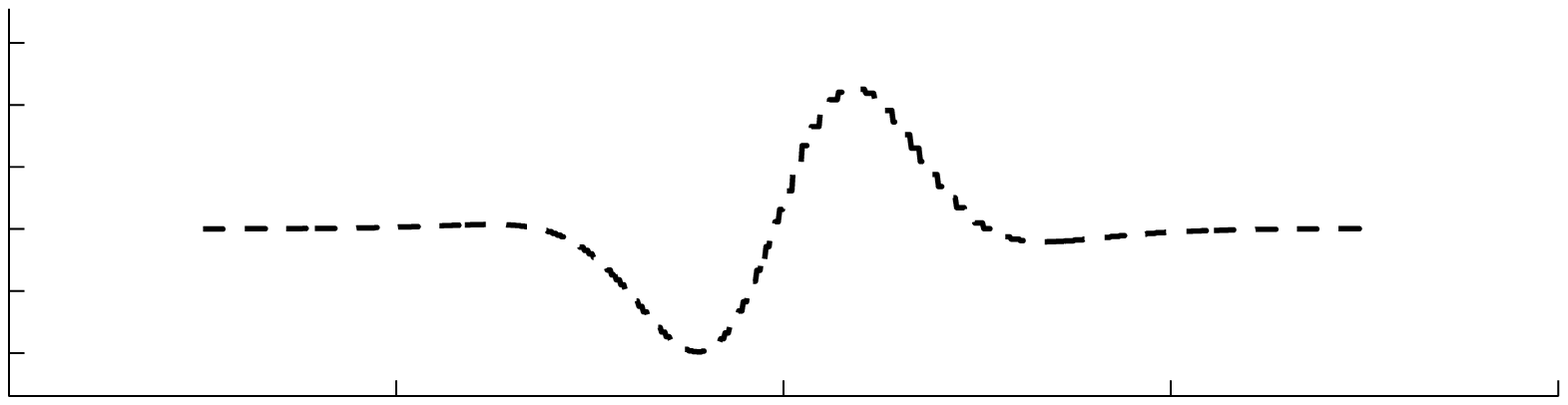}}
\put(-30,160){$\frac 12$ data}
\put(15,140){\includegraphics[width=4.75cm]{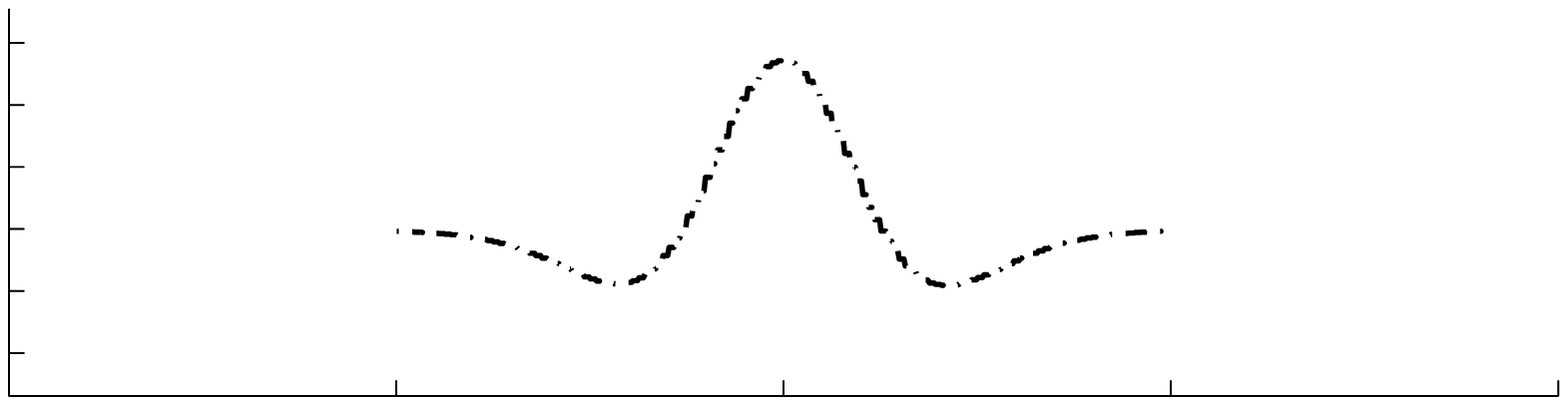}}
\put(170,140){\includegraphics[width=4.75cm]{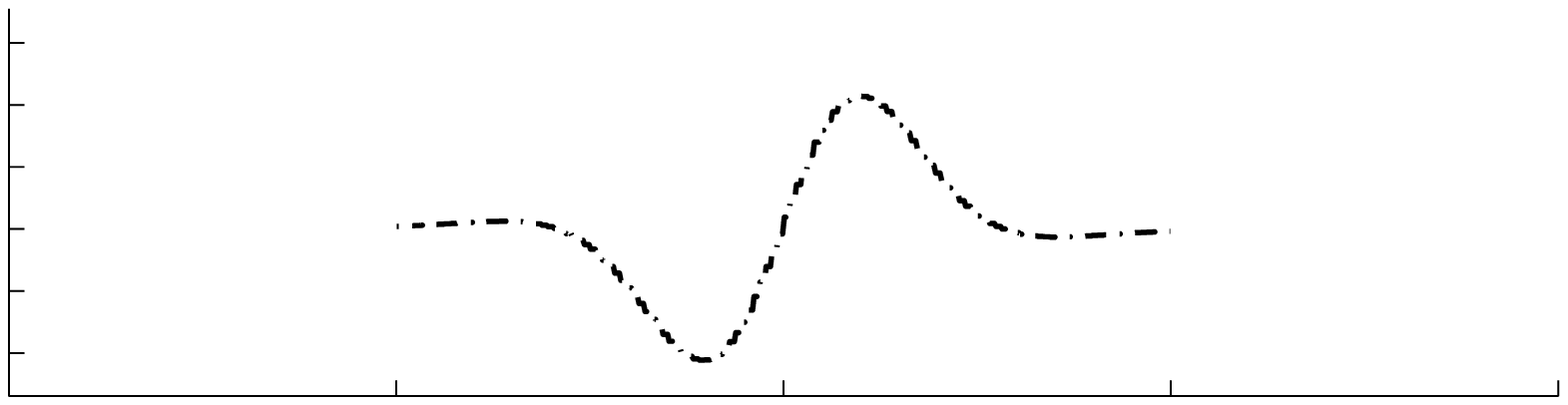}}
\put(-30,100){$\frac 14$ data}
\put(15,80){\includegraphics[width=4.75cm]{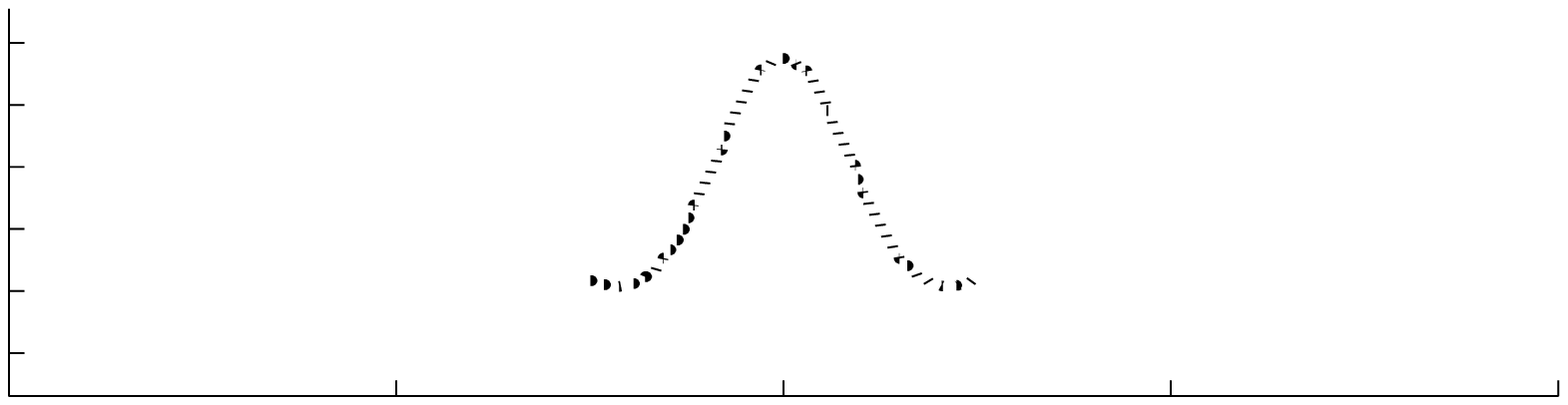}}
\put(170,80){\includegraphics[width=4.75cm]{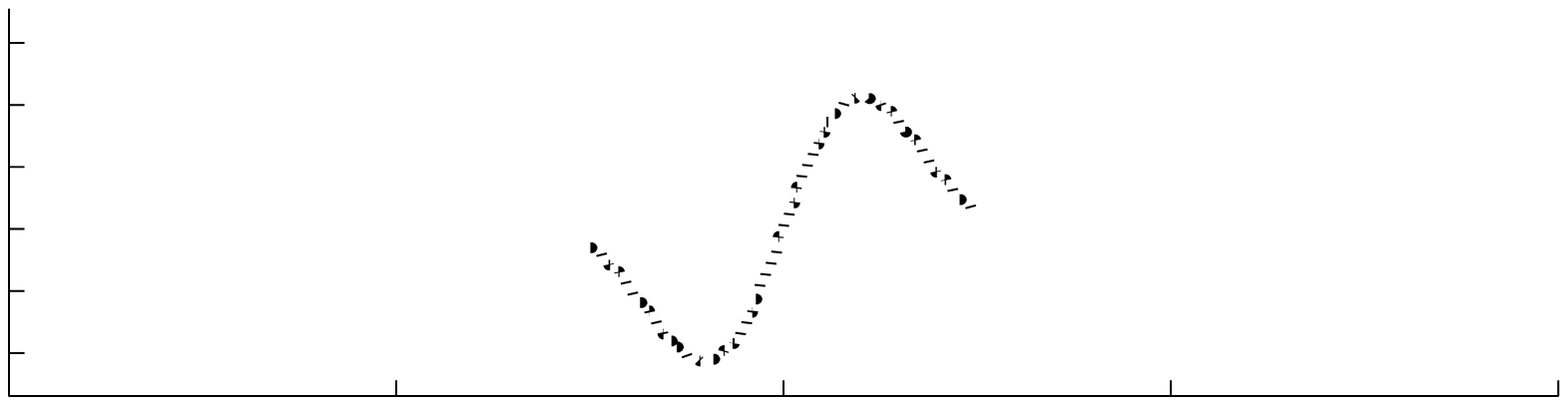}}
\put(-30,40){All}
\put(15,20){\includegraphics[width=4.75cm]{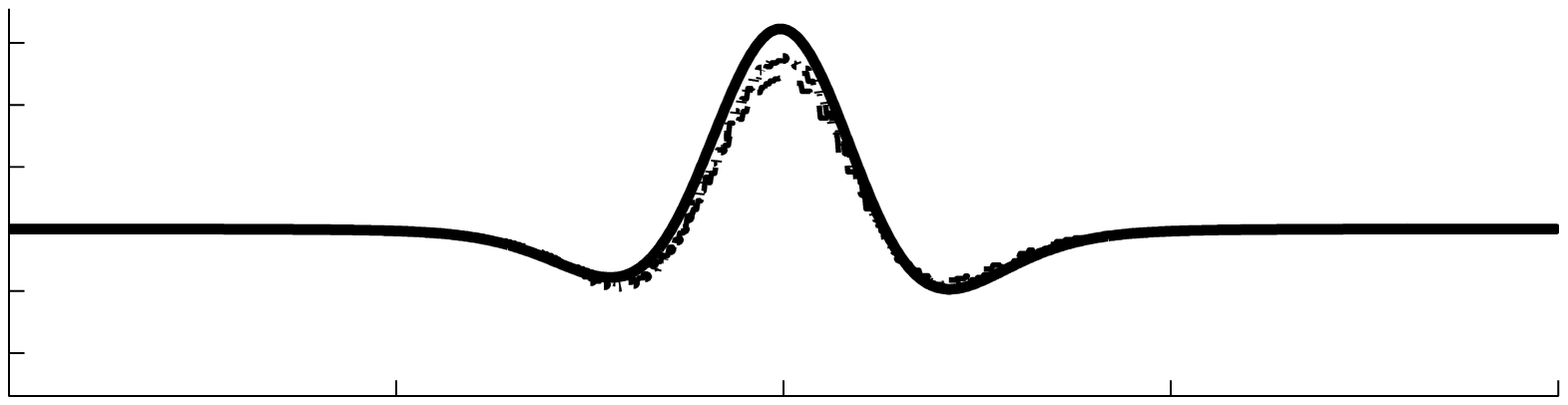}}
\put(170,20){\includegraphics[width=4.75cm]{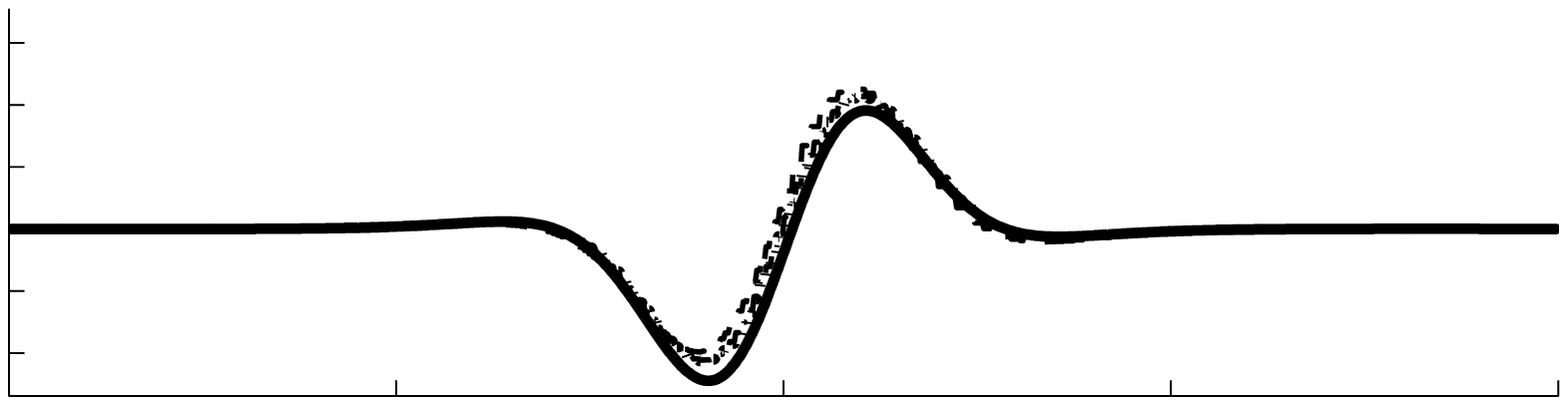}}
\put(10,10){\small$-\pi$}
\put(38,10){\small$-\pi/2$}
\put(81,10){\small$0$}
\put(110,10){\small$\pi/2$}
\put(147,10){\small$\pi$}

\put(164,10){\small$-\pi$}
\put(192,10){\small$-\pi/2$}
\put(236,10){\small$0$}
\put(265,10){\small$\pi/2$}
\put(302,10){\small$\pi$}
\end{picture}
\caption{\label{fig:tracesC}Traces of the CGO solutions  $\psi$ corresponding to the $C^2$ conductivity in Figure~\ref{fig-sigma-phantom}. Here $k=-4i$.}
\end{figure}

Figures~\ref{fig:tracesA} and \ref{fig:tracesC} show the recovered traces of the CGO solutions $\psi(z,k)$ for $k=0.5$ and $-4i$ respectively, plotted against the true traces produced via the Lippmann-Schwinger computation.  Preliminary results suggest that for small magnitude $k$, the partial data CGO solutions agree with the full data (and true) solutions on the accessible part of the domain.  As the magnitude of the frequency parameter $k$ increases, the partial data CGO solutions begin to drift slightly from the full data solutions.  However, in the nonlinear CGO approaches typically used in EIT imaging, only low frequency CGO solutions are used and therefore these results are very promising.
%+++++++++++++++++++++++++++++++
%+++++++++++++++++++++++++++++++
%%+++++++++++++++++++++++++++++++

%%+++++++++++++++++++++++++++++++
%+++++++++++++++++++++++++++++++

%.........................................................................................................................
\subsection{Test~2: D-Bar Reconstructions of Conductivities using Partial Data CGO Solutions}
\begin{figure}[h!]
\centering
\includegraphics[width=2.5in]{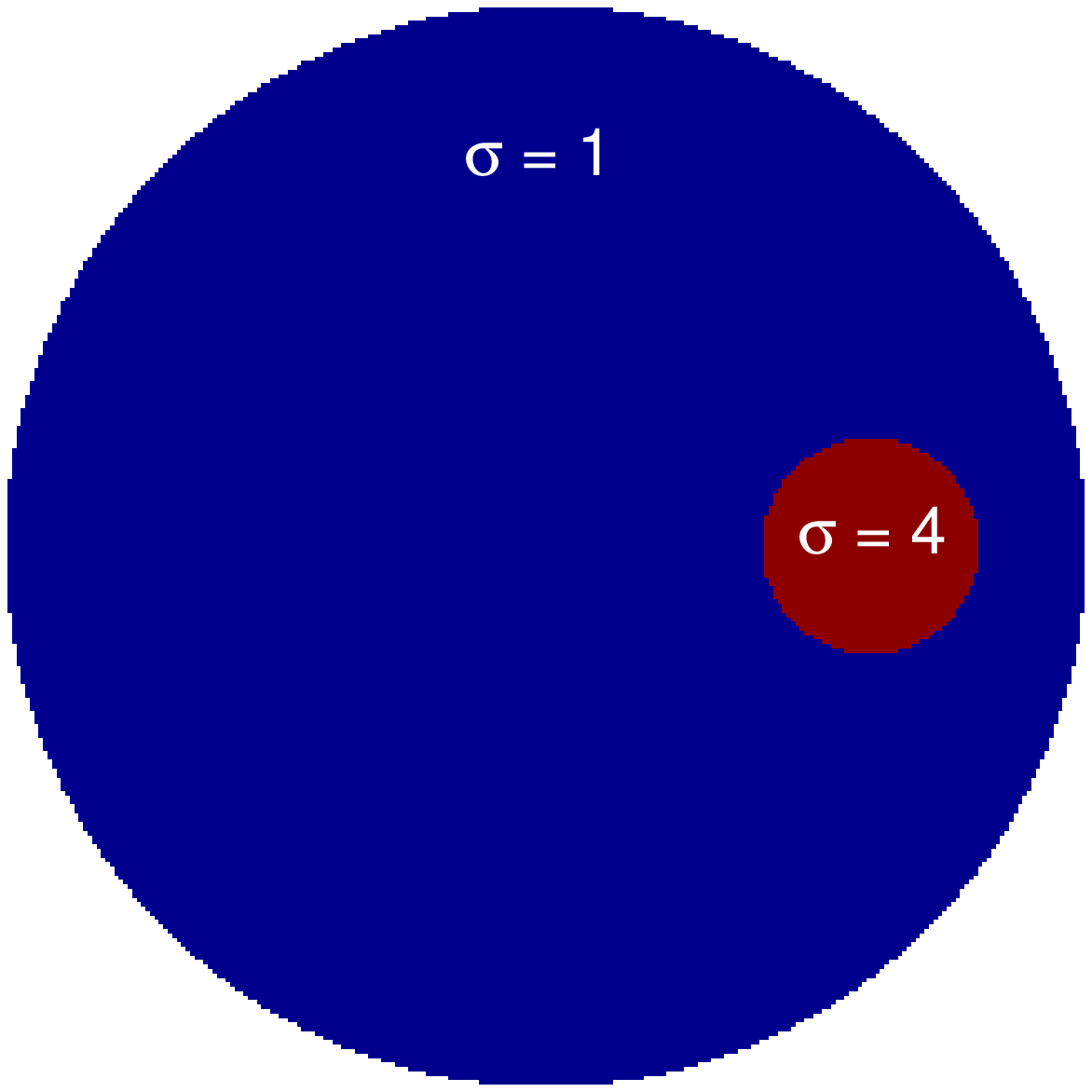}
\caption{The discontinuous conductivity used in Test 2.}\label{fig-Test-sig-2}
\end{figure}

%.........................................................................................................................
In Test~1, we saw that the traces of the CGO solutions computed using the partial D-N data corresponding to an accessible subset $\Gamma$ of the boundary $\bndry$ match the true, as well as full data, traces very well on $\Gamma$ for low frequencies $k$ and $C^2$ smooth conductivities.  Real life situations often involve cases where the conductivities are not smooth but instead only bounded.  In order to determine if the traces of the partial data CGO solutions still provide useful information when the smoothness is relaxed, we consider the discontinuous test conductivity shown in Figure~\ref{fig-Test-sig-2}.
%+++++++++++++++++++++++++++++++
%+++++++++++++++++++++++++++++++
This test phantom could represent a saline filled tank containing an object of higher conductivity.

\begin{figure}
\hspace{2em}
\begin{picture}(300,360)
\put(15,320){\textbf{Real parts of $\mathbf{u_1}$}}
\put(170,320){\textbf{Imaginary parts $\mathbf{u_1}$}}
\put(-30,280){Full data}
\put(15,260){\includegraphics[width=4.75cm]{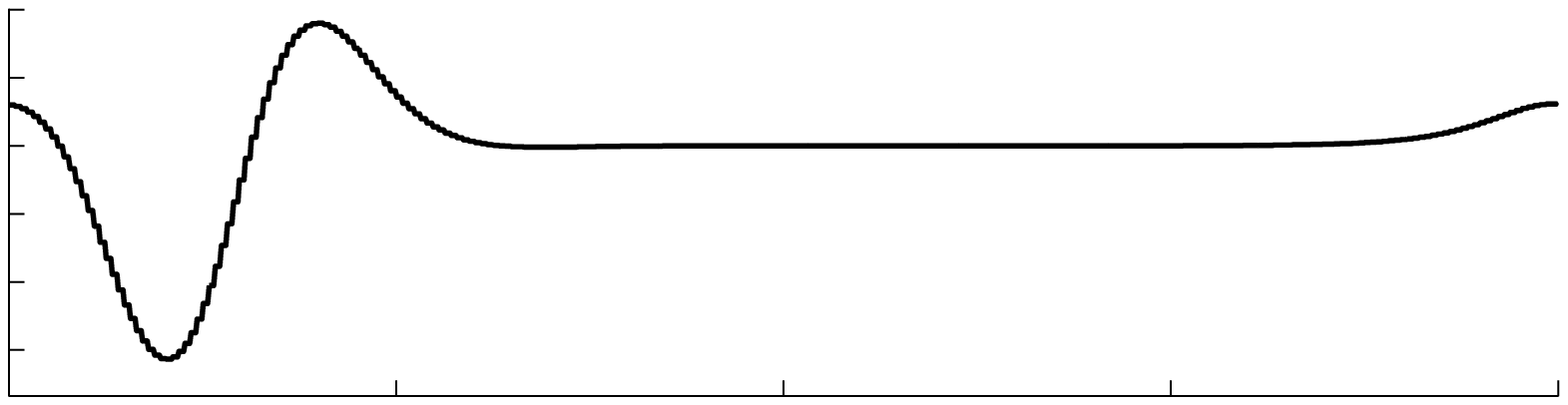}}
\put(170,260){\includegraphics[width=4.75cm]{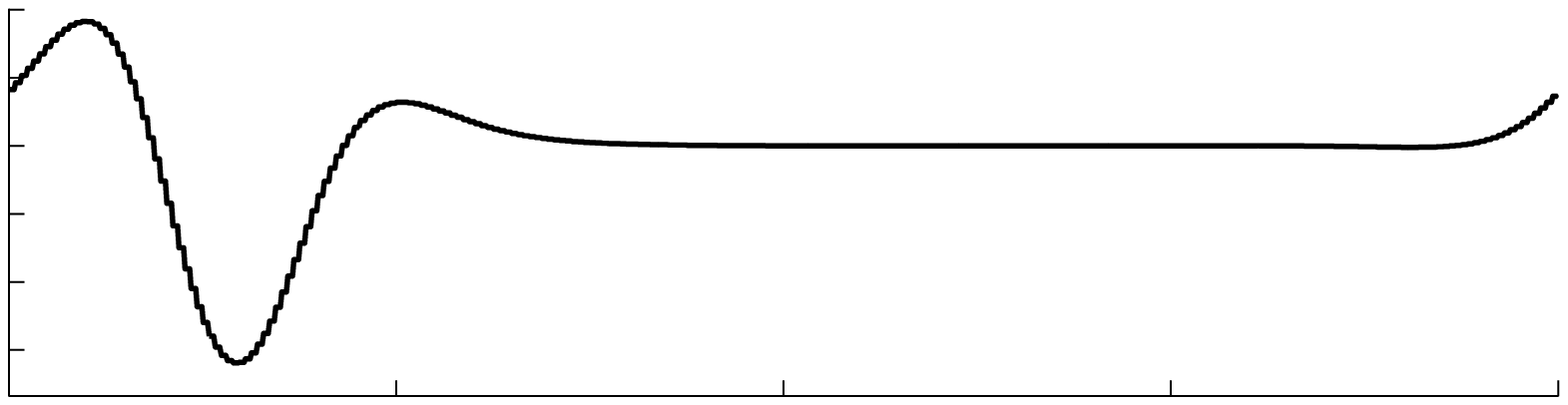}}
\put(-30,220){$\frac 34$ data}
\put(15,200){\includegraphics[width=4.75cm]{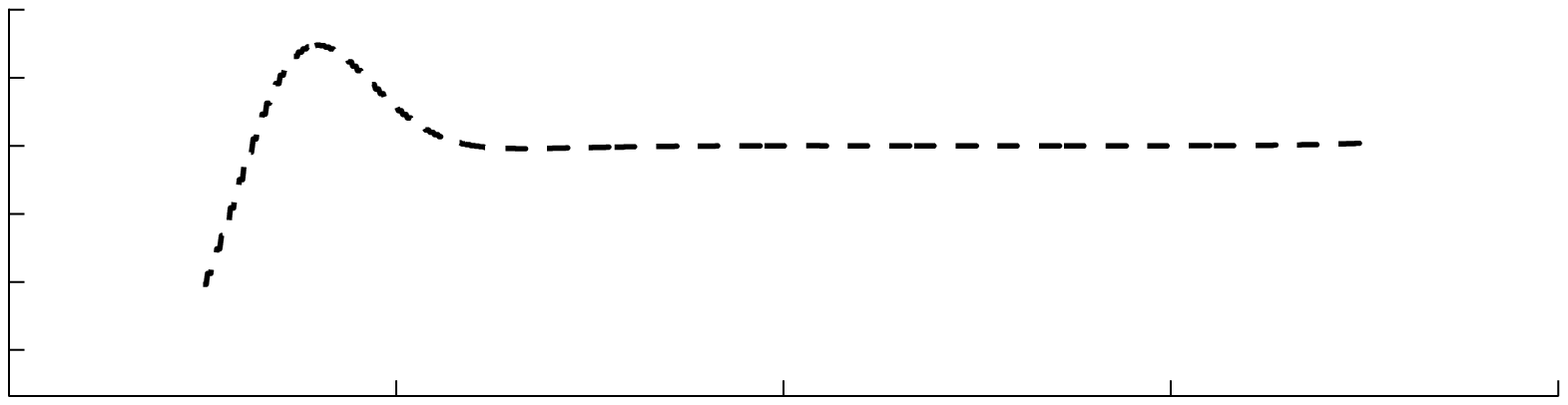}}
\put(170,200){\includegraphics[width=4.75cm]{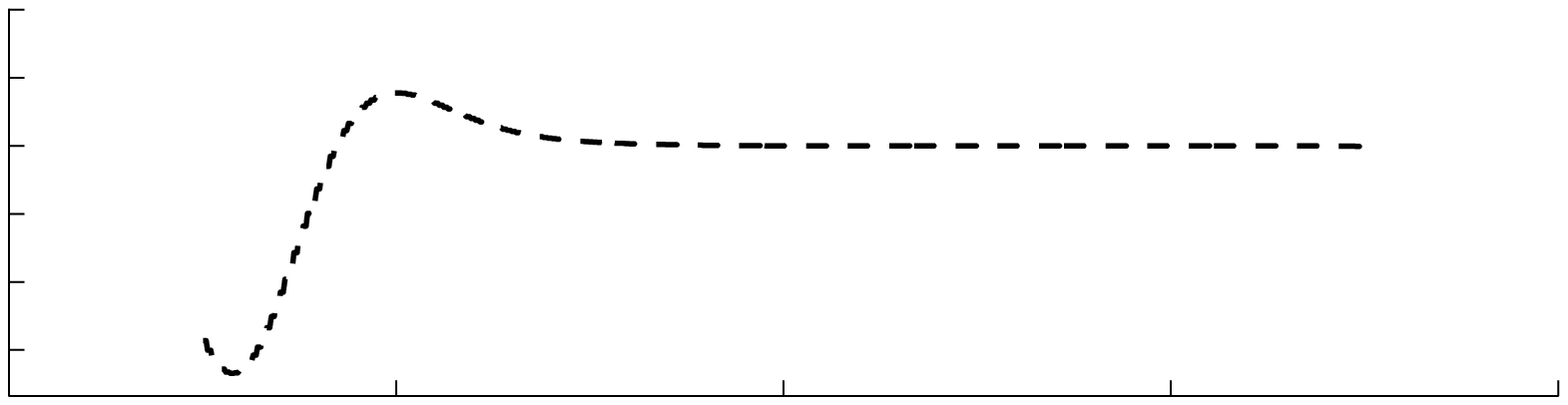}}
\put(-30,160){$\frac 12$ data}
\put(15,140){\includegraphics[width=4.75cm]{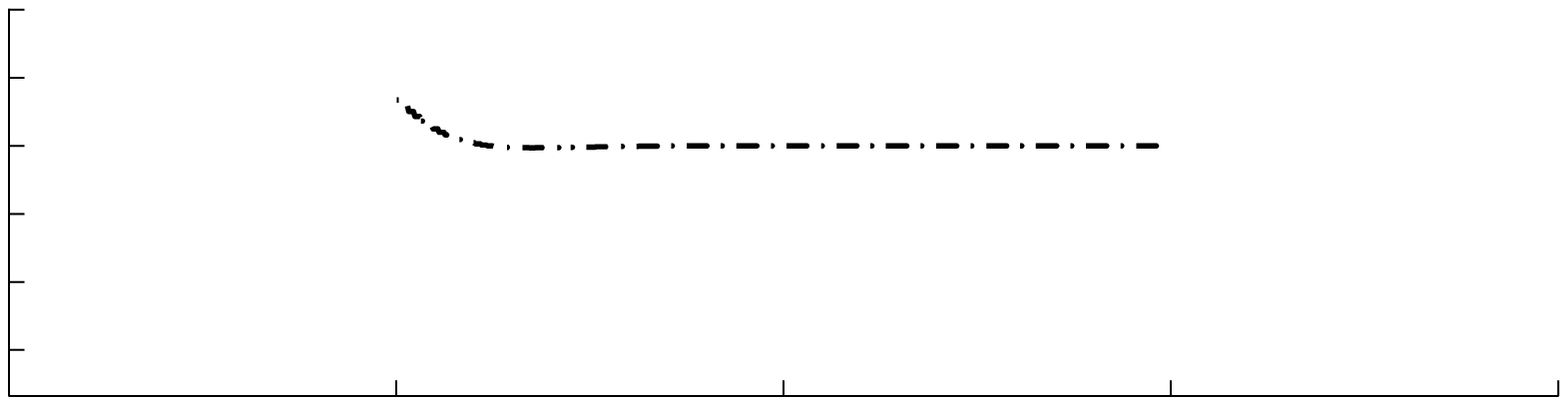}}
\put(170,140){\includegraphics[width=4.75cm]{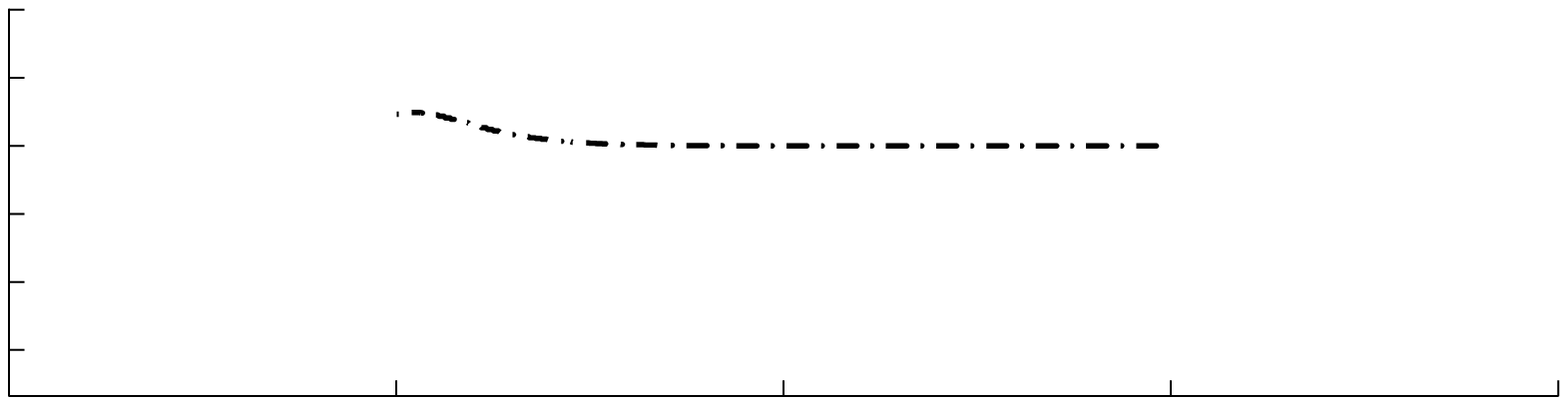}}
\put(-30,100){$\frac 14$ data}
\put(15,80){\includegraphics[width=4.75cm]{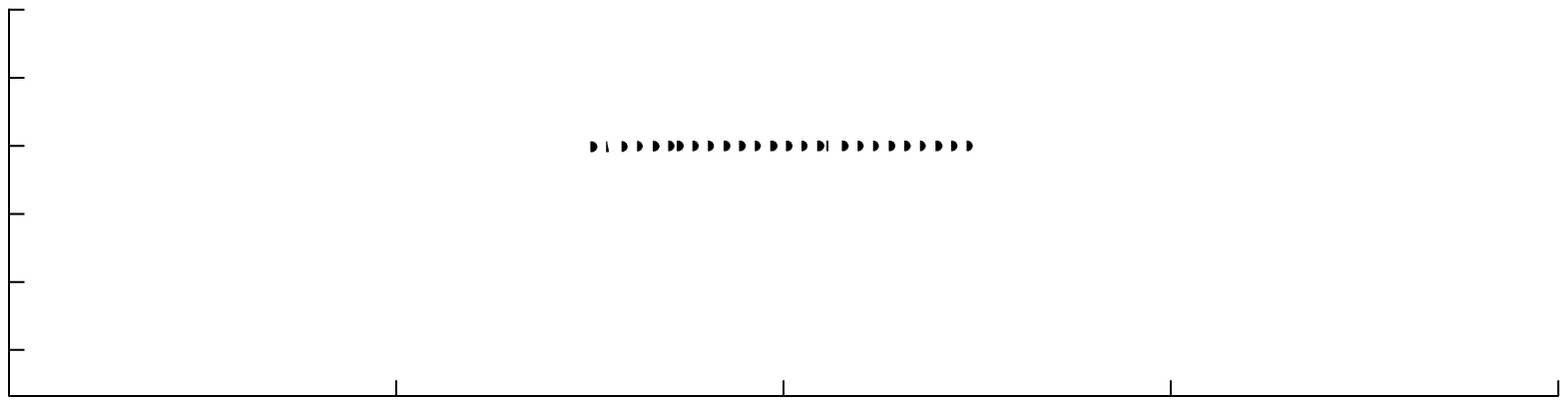}}
\put(170,80){\includegraphics[width=4.75cm]{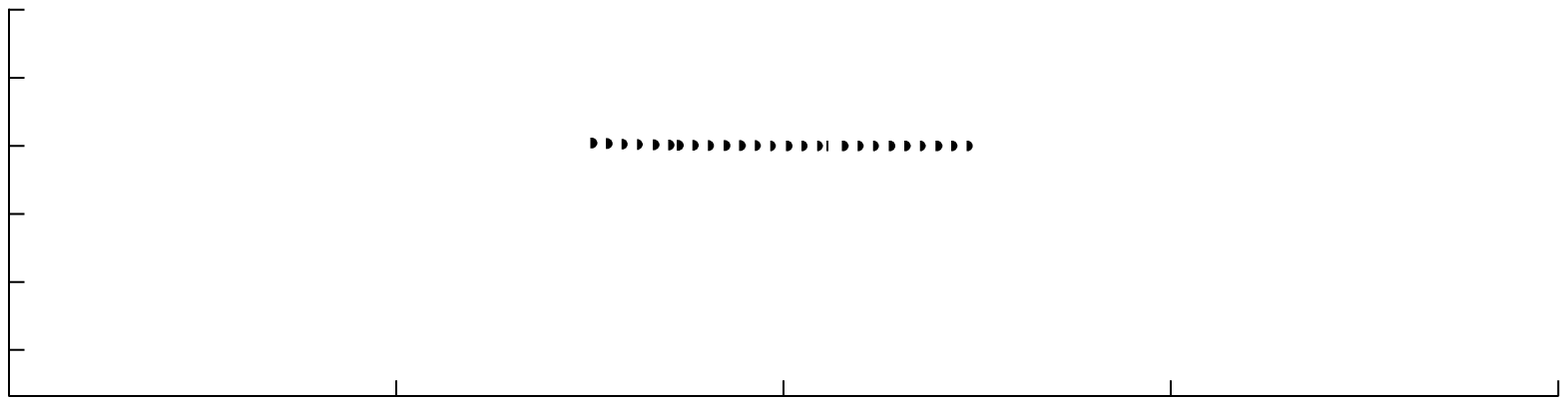}}
\put(-30,40){All}
\put(15,20){\includegraphics[width=4.75cm]{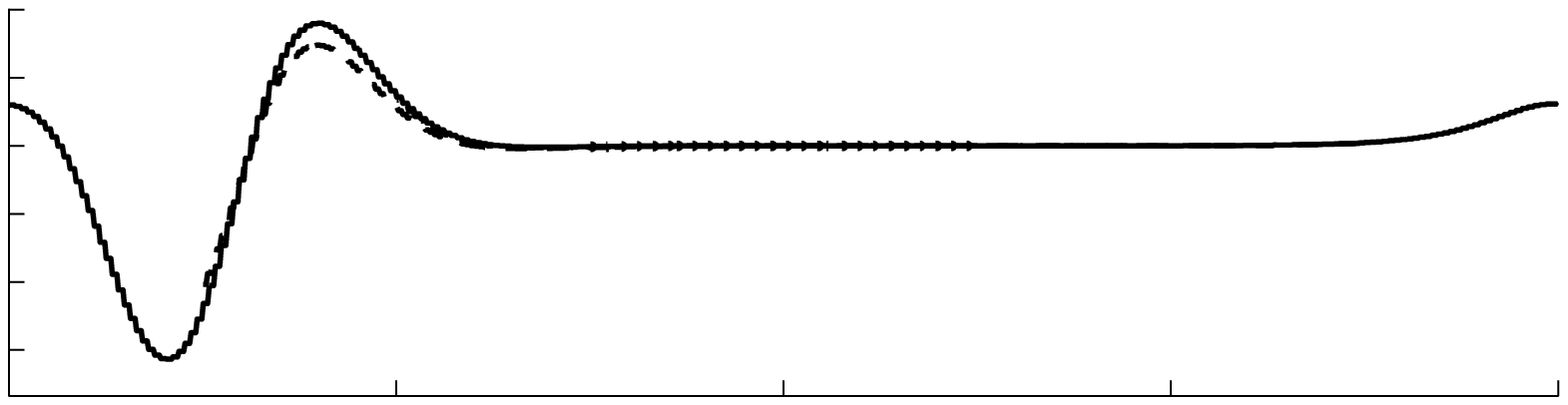}}
\put(170,20){\includegraphics[width=4.75cm]{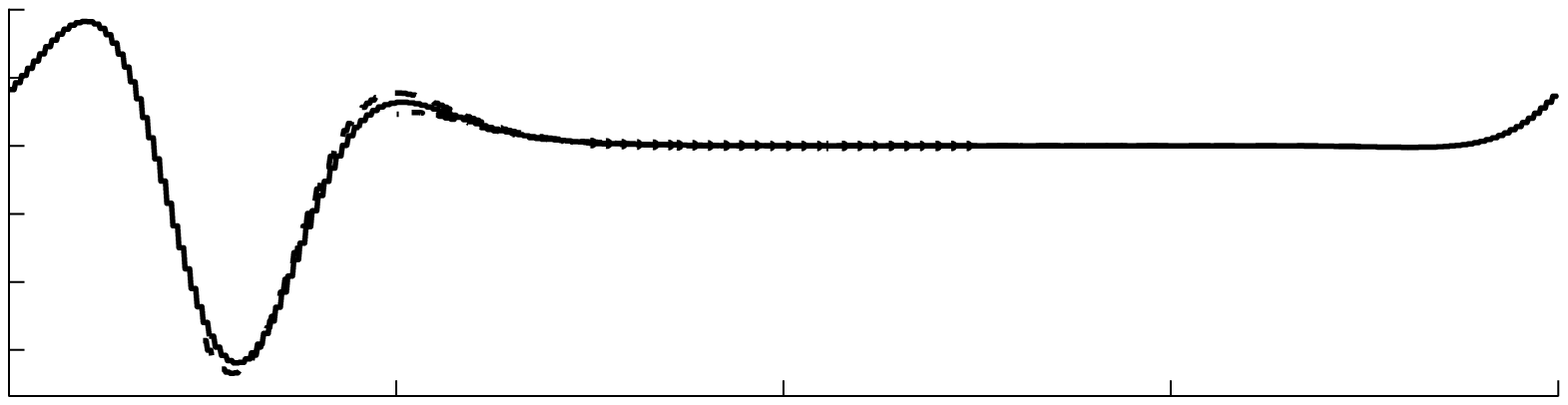}}
\put(10,10){\small$-\pi$}
\put(38,10){\small$-\pi/2$}
\put(81,10){\small$0$}
\put(110,10){\small$\pi/2$}
\put(147,10){\small$\pi$}

\put(164,10){\small$-\pi$}
\put(192,10){\small$-\pi/2$}
\put(236,10){\small$0$}
\put(265,10){\small$\pi/2$}
\put(302,10){\small$\pi$}
\end{picture}
\caption{\label{fig:traces-u1}Traces of the CGO solutions $u_1$ corresponding to the discontinuous conductivity in Figure~\ref{fig-Test-sig-2}. Here $k=3+3i$.}
\end{figure}

\begin{figure}
\hspace{2em}
\begin{picture}(300,360)
\put(15,320){\textbf{Real parts of $\mathbf{u_2}$}}
\put(170,320){\textbf{Imaginary parts $\mathbf{u_2}$}}
\put(-30,280){Full data}
\put(15,260){\includegraphics[width=4.75cm]{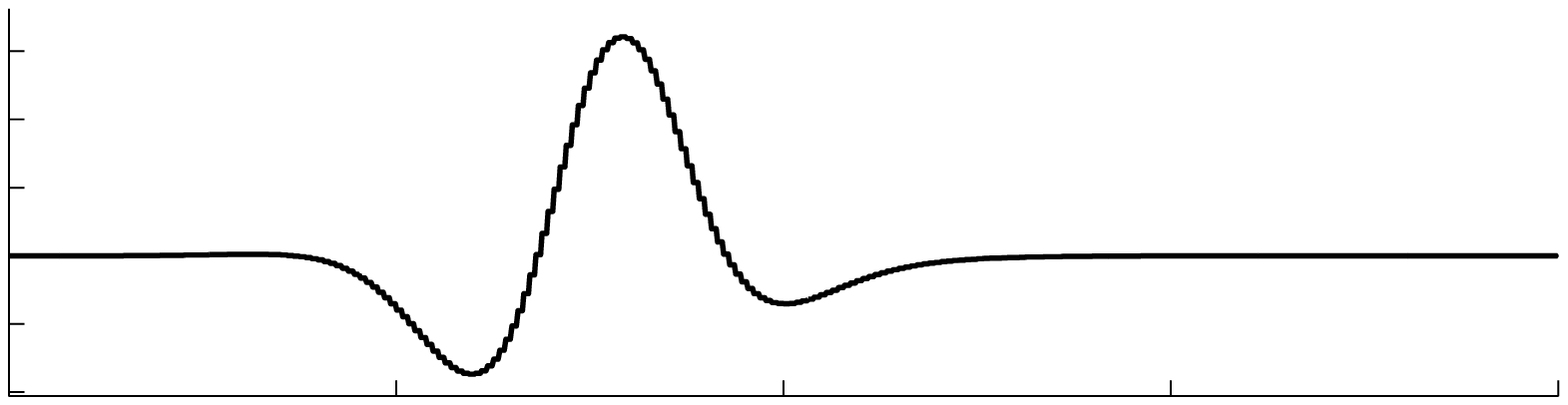}}
\put(170,260){\includegraphics[width=4.75cm]{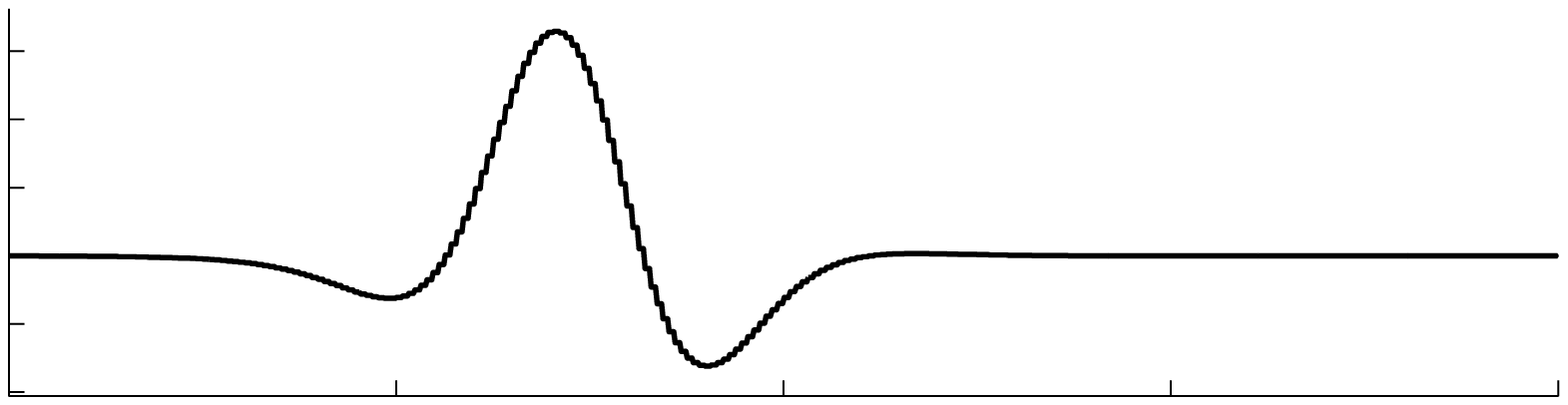}}
\put(-30,220){$\frac 34$ data}
\put(15,200){\includegraphics[width=4.75cm]{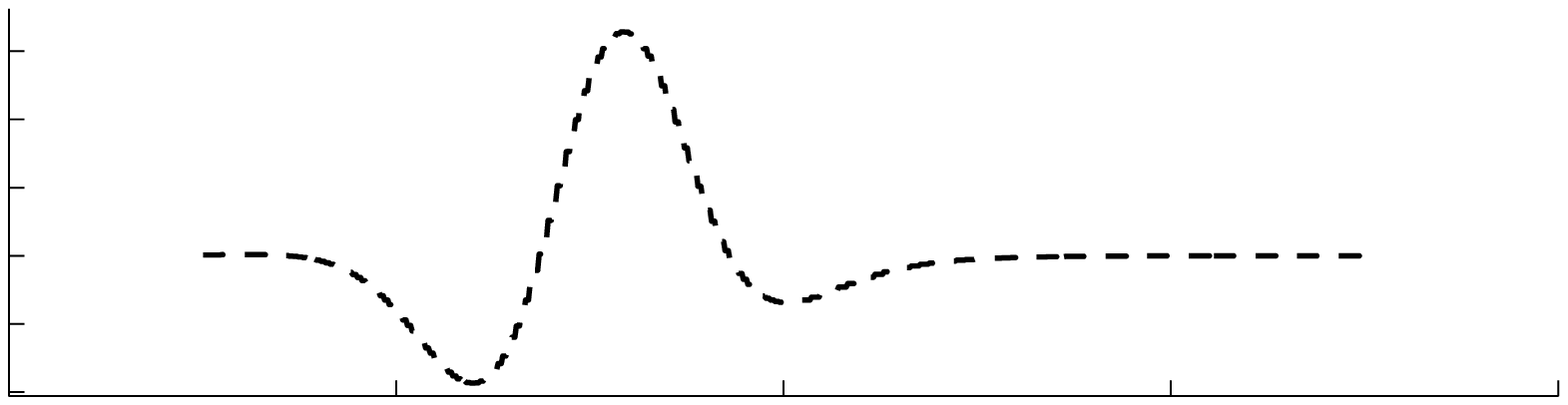}}
\put(170,200){\includegraphics[width=4.75cm]{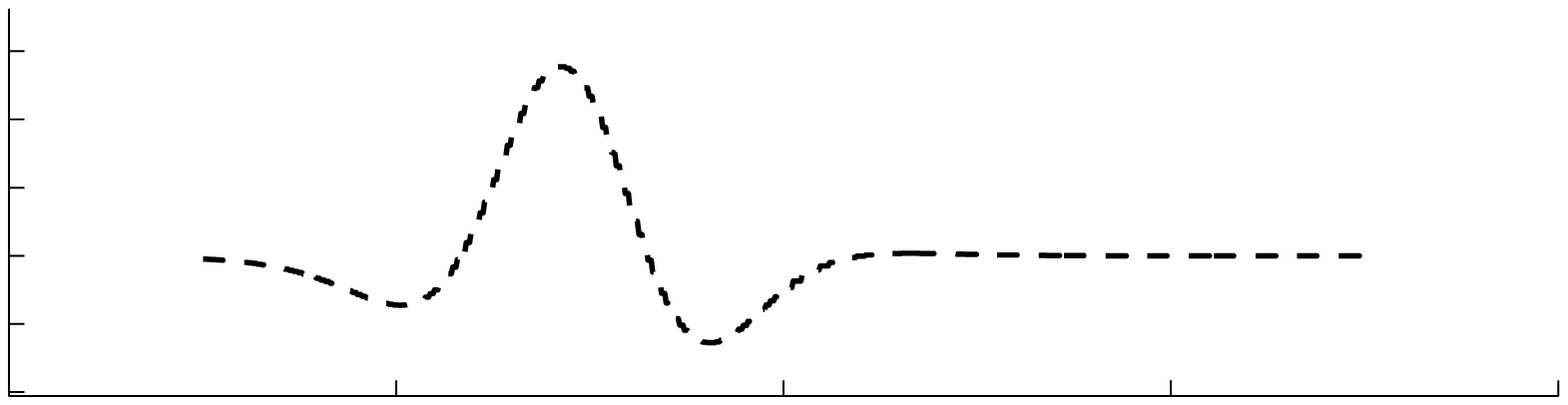}}
\put(-30,160){$\frac 12$ data}
\put(15,140){\includegraphics[width=4.75cm]{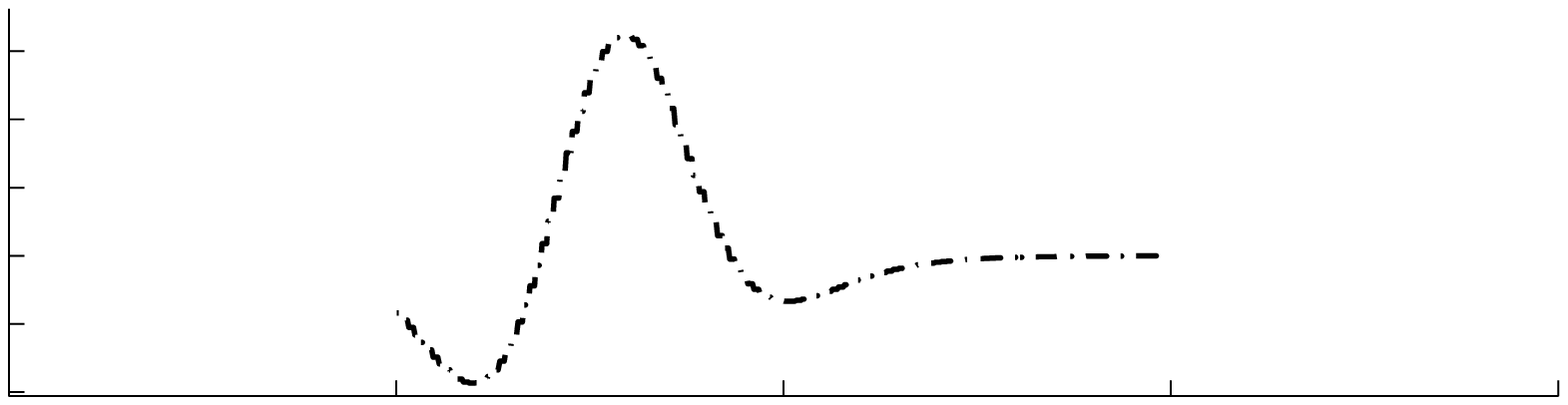}}
\put(170,140){\includegraphics[width=4.75cm]{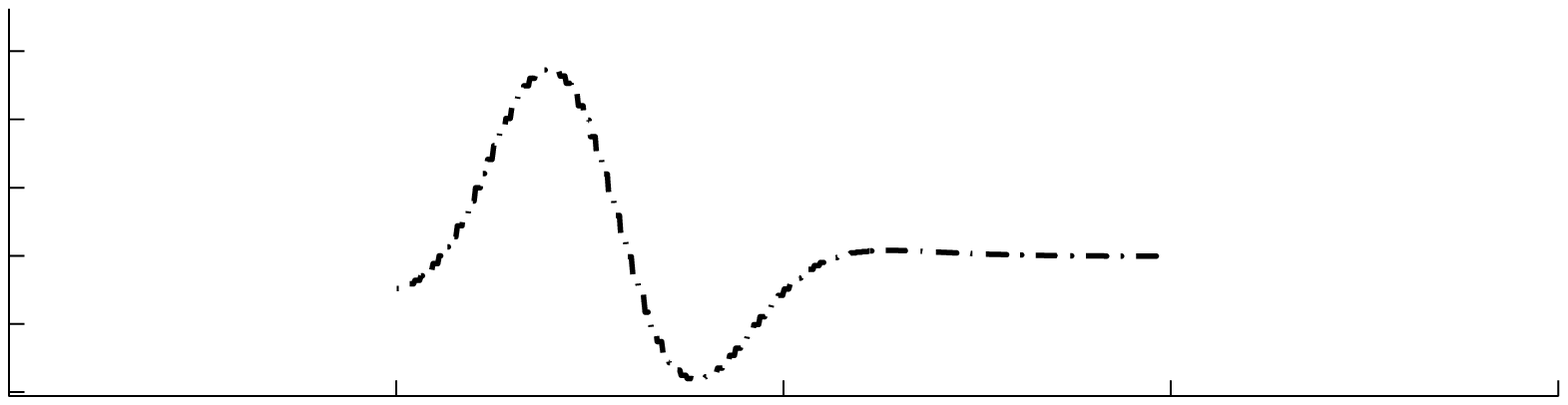}}
\put(-30,100){$\frac 14$ data}
\put(15,80){\includegraphics[width=4.75cm]{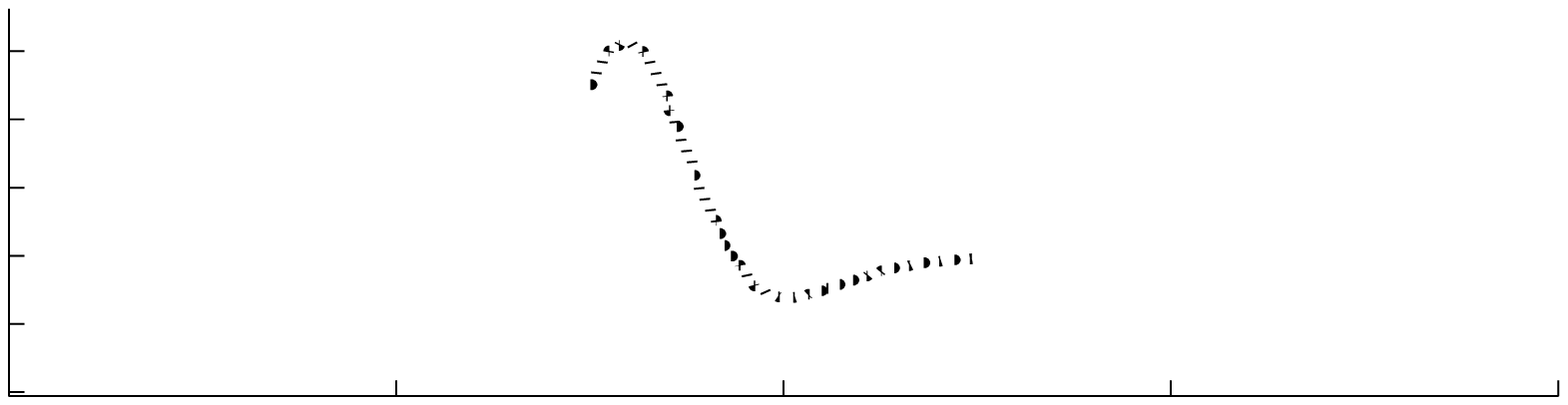}}
\put(170,80){\includegraphics[width=4.75cm]{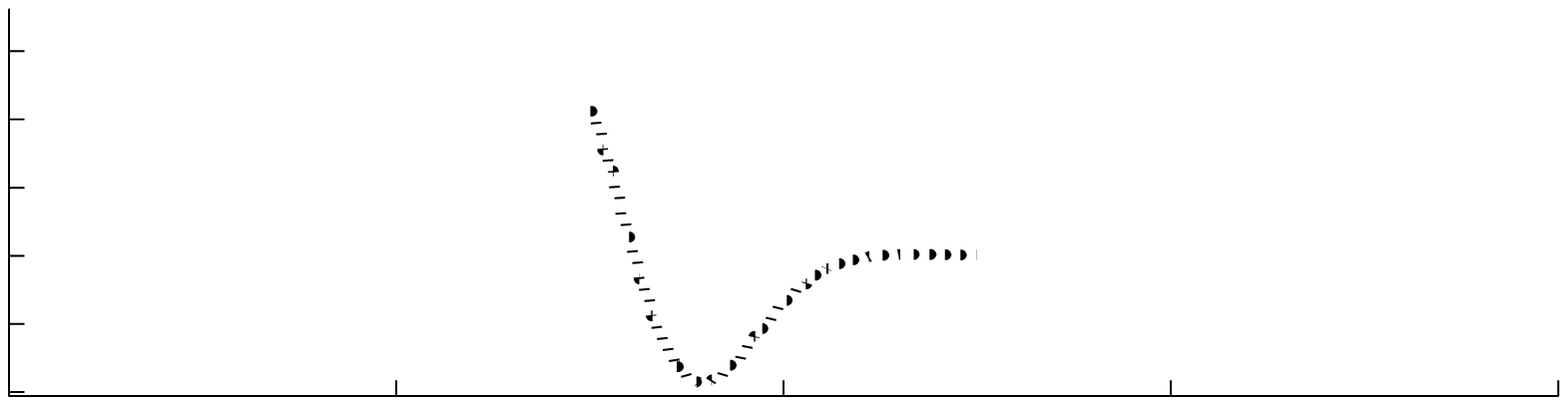}}
\put(-30,40){All}
\put(15,20){\includegraphics[width=4.75cm]{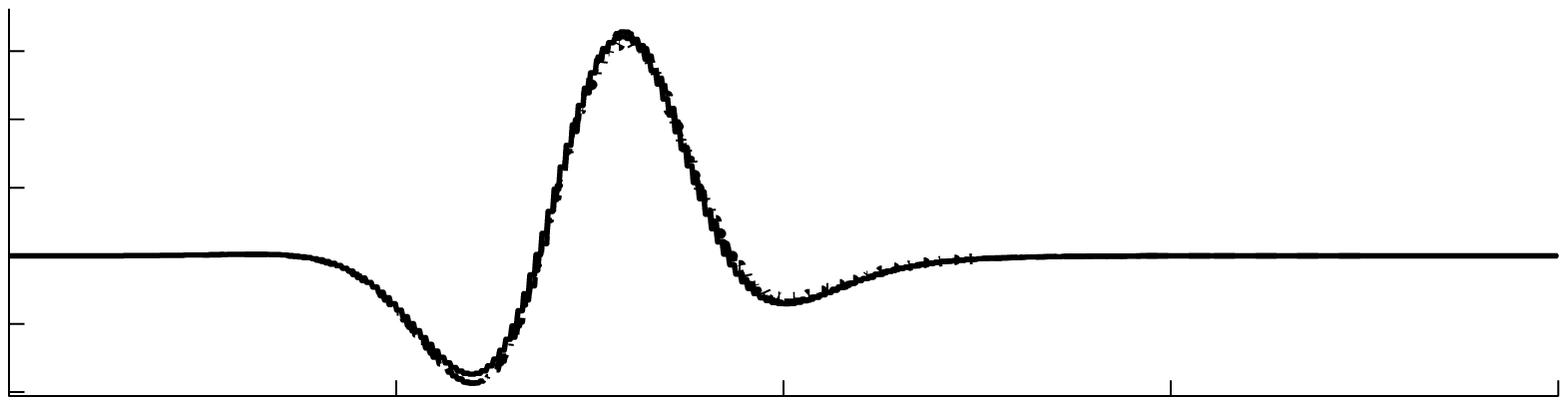}}
\put(170,20){\includegraphics[width=4.75cm]{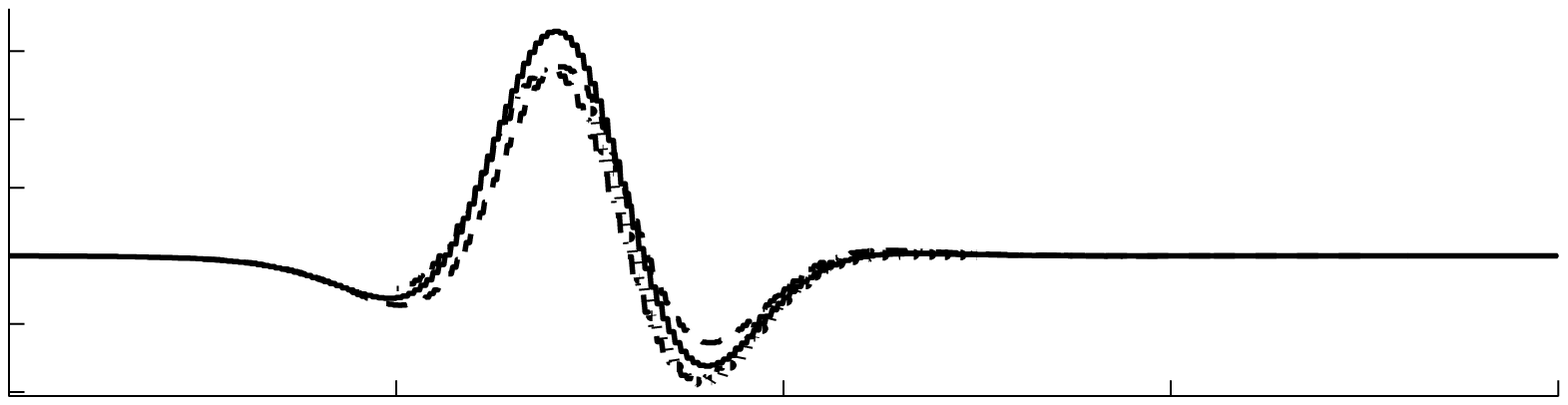}}
\put(10,10){\small$-\pi$}
\put(38,10){\small$-\pi/2$}
\put(81,10){\small$0$}
\put(110,10){\small$\pi/2$}
\put(147,10){\small$\pi$}

\put(164,10){\small$-\pi$}
\put(192,10){\small$-\pi/2$}
\put(236,10){\small$0$}
\put(265,10){\small$\pi/2$}
\put(302,10){\small$\pi$}
\end{picture}
\caption{\label{fig:traces-u2}Traces of the CGO solutions $u_2$ corresponding to the discontinuous conductivity in Figure~\ref{fig-Test-sig-2}. Here $k=3+3i$.}
\end{figure}

As the smoothness condition for Method~2 is violated (as well as that for Method~1), we cannot compare the full or partial data traces of the CGO solutions to their ``true'' traces.  Precisely, we cannot compute the potential $Q$ (or $q$) and the associated $\dez-\dbarz$ system (or Lippmann-Schwinger equation) for the ``true'' traces of the CGO solutions.  Instead, we compare the partial D-N data traces of the CGO solutions to the corresponding full D-N data traces.  Figures~\ref{fig:traces-u1} and \ref{fig:traces-u2} shows the reconstructed traces of the CGO solutions $\tilde{u}_1$ and $\tilde{u}_2$, respectively, for $k=3+3i$ resulting from the conductivity distribution in Figure~\ref{fig-Test-sig-2} plotted for full, $3/4$, $1/2$, and $1/4$ D-N data  with 256, 192, 128, and 64 Haar wavelets respectively.   Clearly the partial D-N data traces of the CGO solutions appear to approximate the full D-N data traces of the CGO solutions in the accessible region $\Gamma$ of the boundary.

%+++++++++++++++++++++++++++++++
%+++++++++++++++++++++++++++++++

%+++++++++++++++++++++++++++++++
%+++++++++++++++++++++++++++++++

\begin{figure}
\hspace{2em}
\begin{picture}(300,150)
\put(20,100){\textbf{Full Data}}
\put(10,20){\includegraphics[width=2.5cm]{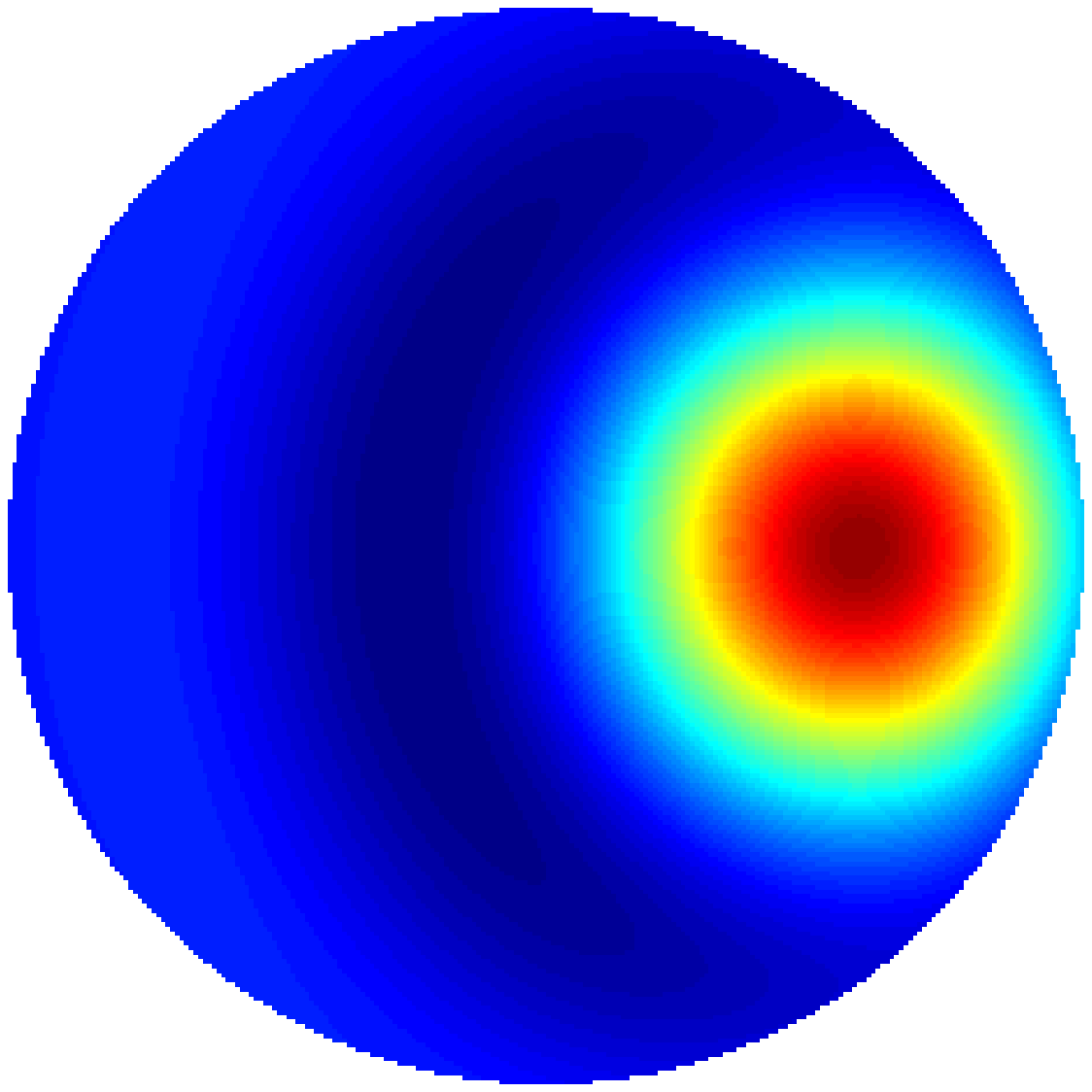}}
\put(17,8){\small$\begin{array}{rl}
\text{max}=& \hspace*{-0.5em}1.19\\
\text{min} =& \hspace*{-0.5em}0.97
\end{array}$}

\put(102,100){$\mathbf{\frac 34}$ \textbf{Data}}
\put(85,20){\includegraphics[width=2.5cm]{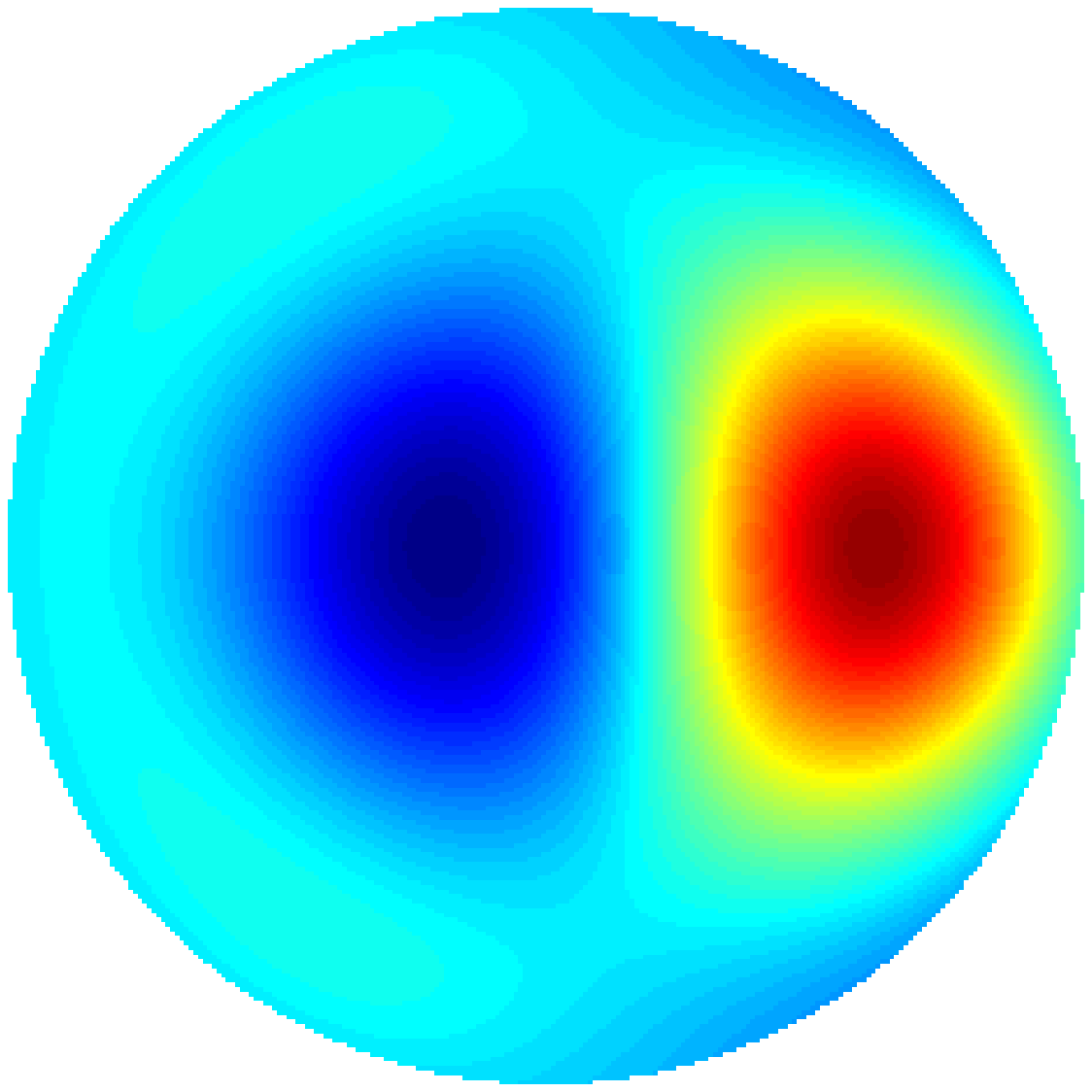}}
\put(92,8){\small$\begin{array}{rl}
\text{max}=& \hspace*{-0.5em}1.14\\
\text{min} =& \hspace*{-0.5em}0.93
\end{array}$}

\put(180,100){$\mathbf{\frac 12}$ \textbf{Data}}
\put(160,20){\includegraphics[width=2.5cm]{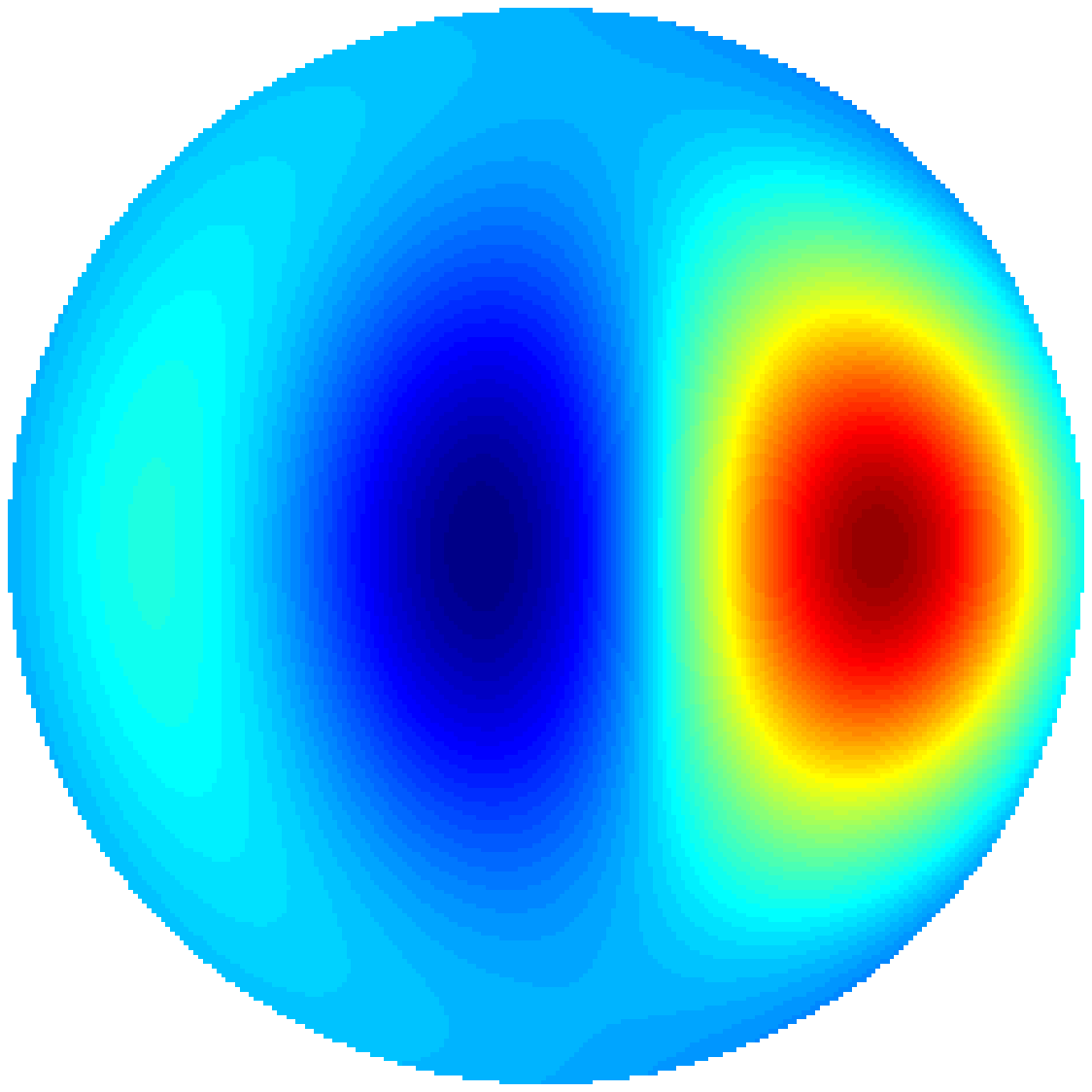}}
\put(167,8){\small$\begin{array}{rl}
\text{max}=& \hspace*{-0.5em}1.10\\
\text{min} =& \hspace*{-0.5em}0.96
\end{array}$}

\put(250,100){$\mathbf{\frac 14}$ \textbf{Data}}
\put(235,20){\includegraphics[width=2.5cm]{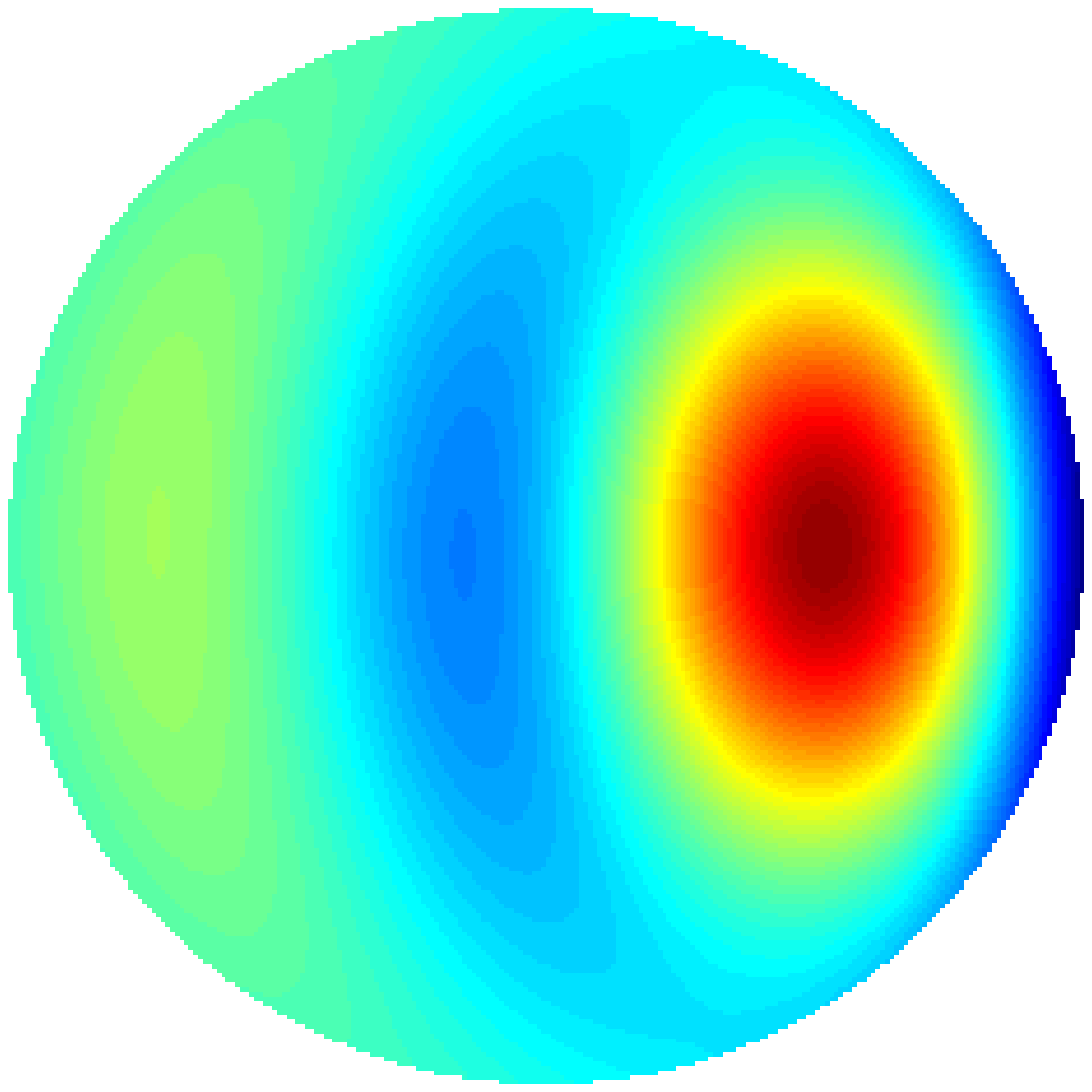}}
\put(242,8){\small$\begin{array}{rl}
\text{max}=& \hspace*{-0.5em}1.06\\
\text{min} =& \hspace*{-0.5em}0.97
\end{array}$}

\put(-30,50){$\mathbf{|k|\leq3}$}
\end{picture}
\caption{Reconstructions of the discontinuous conductivity in Figure~\ref{fig-Test-sig-2} produced using scattering data for $|k|\leq3$ using the method described in Section~\ref{sec-method2-pdata}.  From left to right, the reconstructions are for Full, 3/4, 1/2, and 1/4 Dirichlet-to-Neumann data.}\label{fig-Test-sig-2-Recons-kInt3o0}
\end{figure}

\begin{figure}
\hspace{2em}
\begin{picture}(300,150)
\put(20,100){\textbf{Full Data}}
\put(10,20){\includegraphics[width=2.5cm]{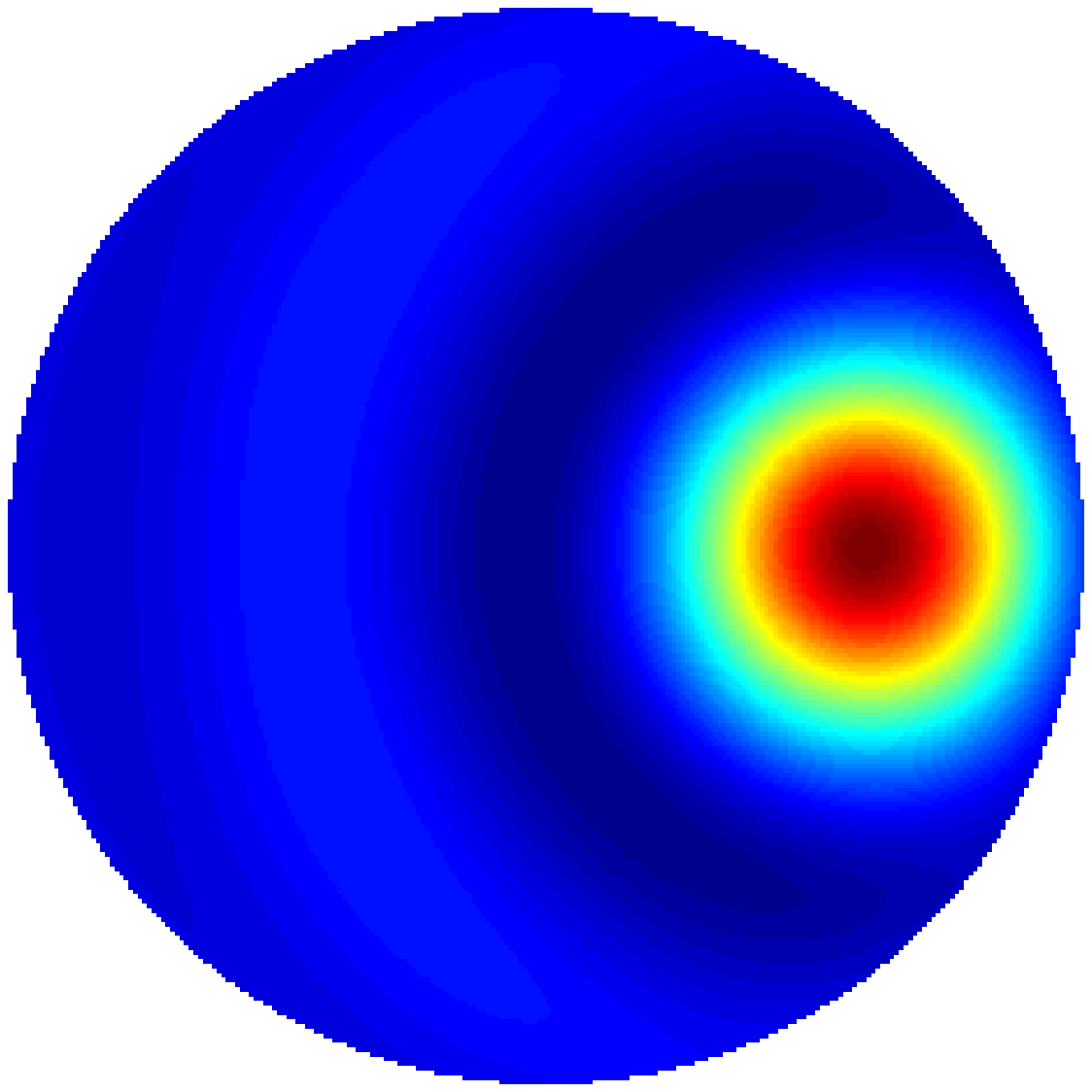}}
\put(17,8){\small$\begin{array}{rl}
\text{max}=& \hspace*{-0.5em}1.38\\
\text{min} =& \hspace*{-0.5em}0.96
\end{array}$}

\put(102,100){$\mathbf{\frac 34}$ \textbf{Data}}
\put(85,20){\includegraphics[width=2.5cm]{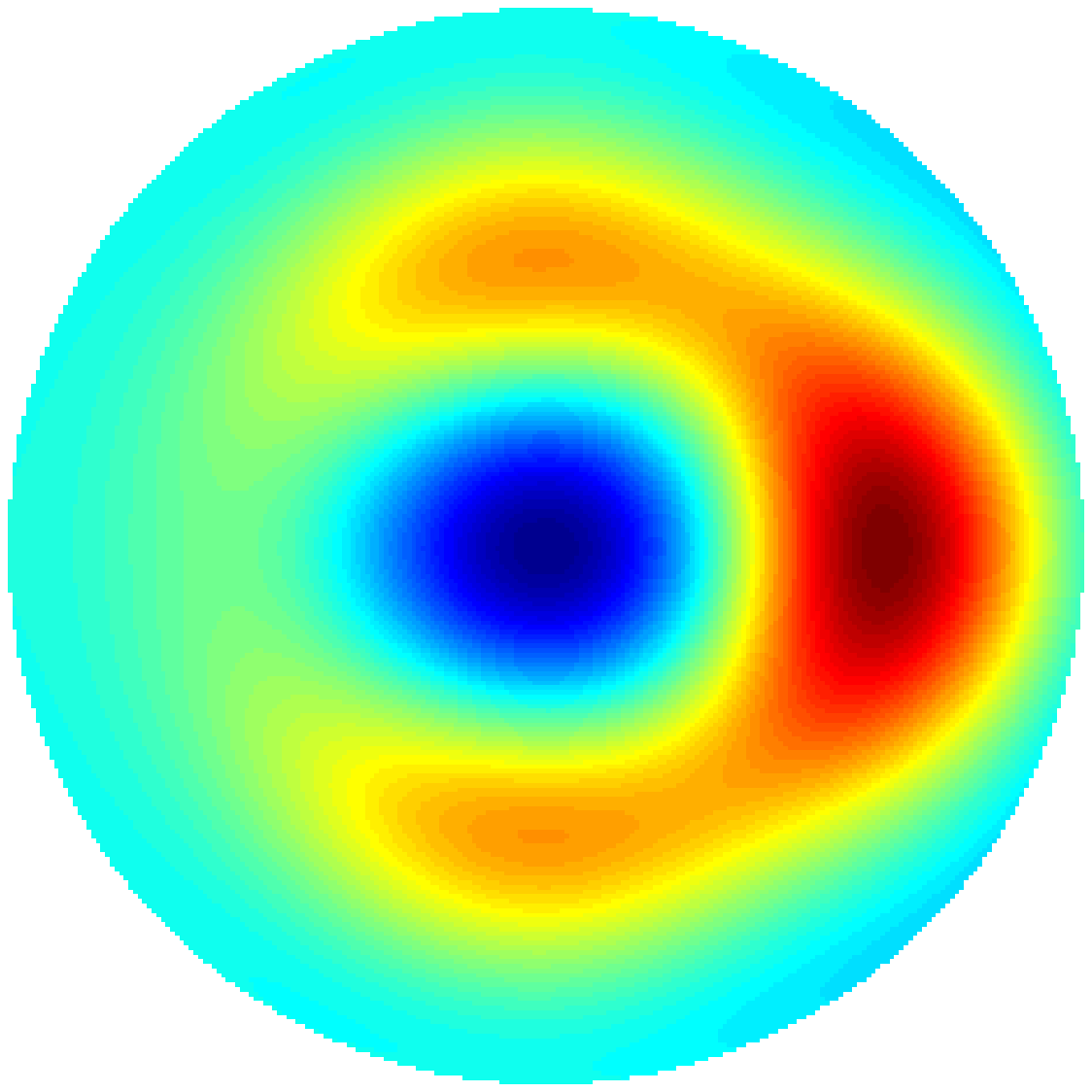}}
\put(92,8){\small$\begin{array}{rl}
\text{max}=& \hspace*{-0.5em}1.30\\
\text{min} =& \hspace*{-0.5em}0.82
\end{array}$}

\put(180,100){$\mathbf{\frac 12}$ \textbf{Data}}
\put(160,20){\includegraphics[width=2.5cm]{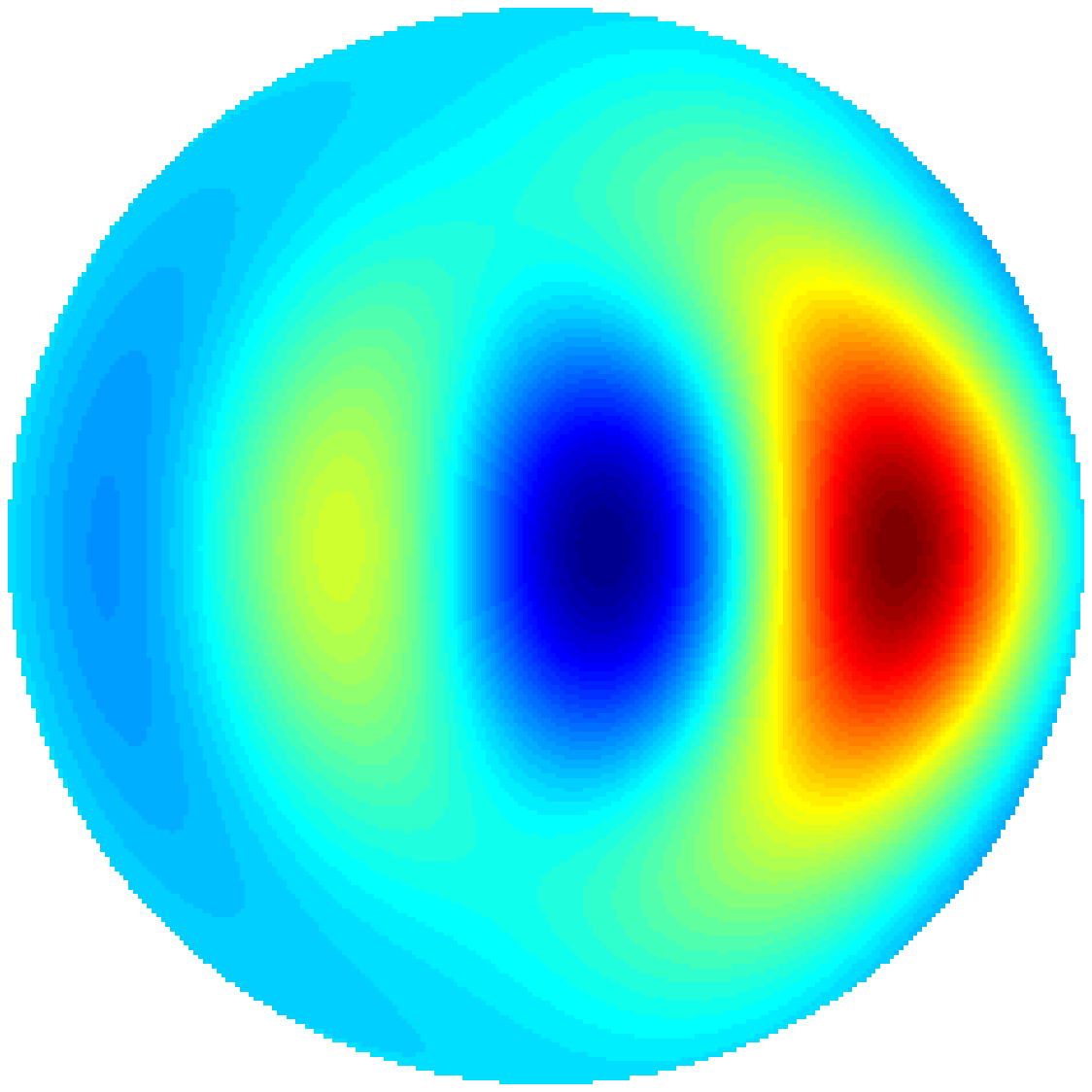}}
\put(167,8){\small$\begin{array}{rl}
\text{max}=& \hspace*{-0.5em}1.24\\
\text{min} =& \hspace*{-0.5em}0.88
\end{array}$}

\put(250,100){$\mathbf{\frac 14}$ \textbf{Data}}
\put(235,20){\includegraphics[width=2.5cm]{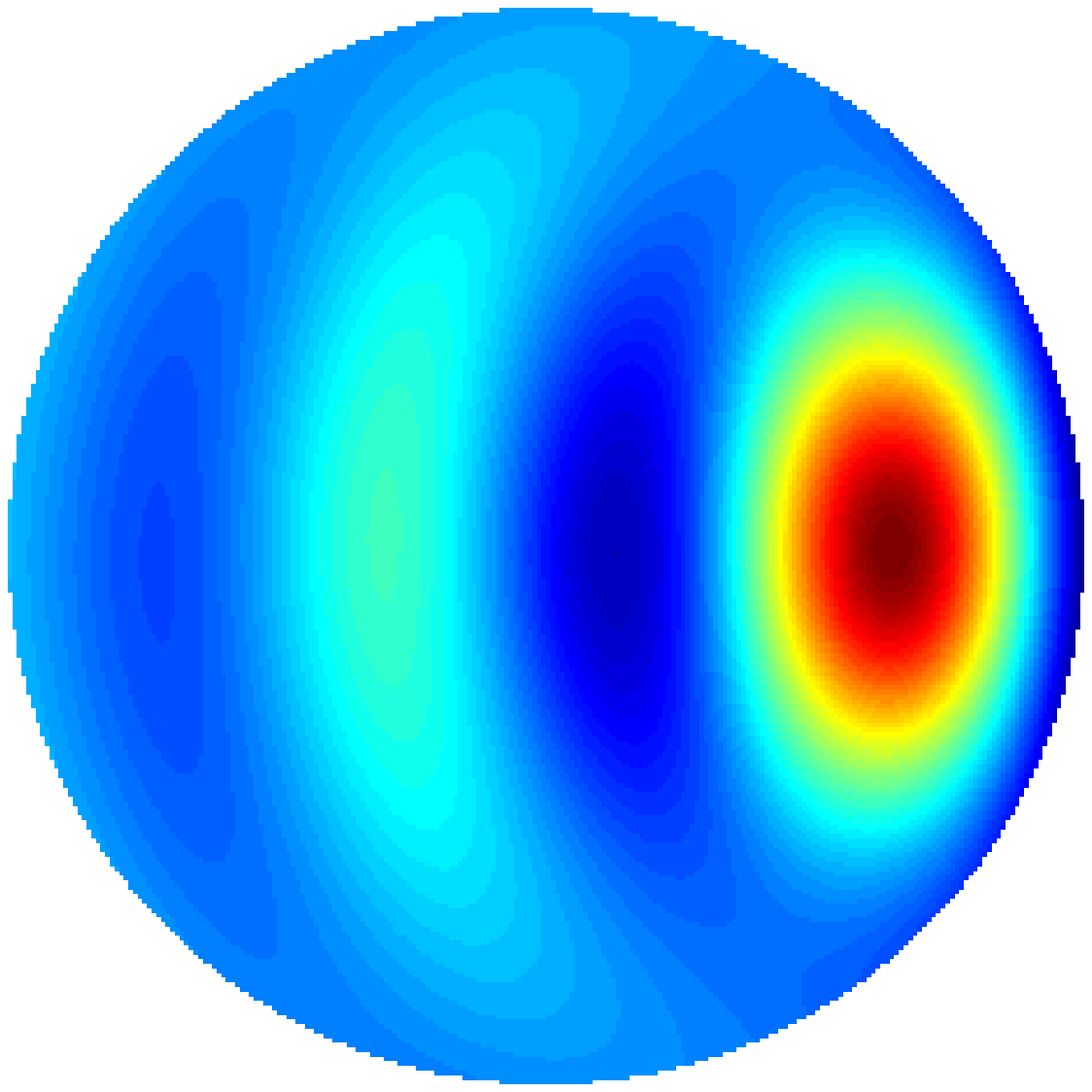}}
\put(242,8){\small$\begin{array}{rl}
\text{max}=& \hspace*{-0.5em}1.13\\
\text{min} =& \hspace*{-0.5em}0.96
\end{array}$}

\put(-30,50){$\mathbf{|k|\leq4}$}
\end{picture}
\caption{Reconstructions of the discontinuous conductivity in Figure~\ref{fig-Test-sig-2} produced using scattering data for $|k|\leq4$ using the method described in Section~\ref{sec-method2-pdata}.  From left to right, the reconstructions are for Full, 3/4, 1/2, and 1/4 Dirichlet-to-Neumann data.}\label{fig-Test-sig-2-Recons-kInt4o0}
\end{figure}

Next we used the partial D-N data traces of the CGO solutions $\tilde{u}_1$ and $\tilde{u}_2$ in the modified D-bar algorithm for Method~2, described above in Section~\ref{sec-method2-pdata}, to determine their effect on the algorithm and thus the reconstructed conductivity distribution.  Figures~\ref{fig-Test-sig-2-Recons-kInt3o0} and \ref{fig-Test-sig-2-Recons-kInt4o0} show the reconstructed conductivity from full, $3/4$, $1/2$, and $1/4$ D-N data using scattering data satisfying $|k|\leq3$ and $|k|\leq4$, respectively.  The range of the reconstructed values decreases with the size of the accessible region of the boundary $\Gamma$, however an object of higher conductivity is clearly visible in all cases.  As the magnitude of $k$ increases the reconstructed values of the conductivity improve.  However, as in the full data D-N case, increasing the scattering radius too much can introduce artifacts into the reconstruction. 

Note that the reconstructions of the conductivity shown in Figures~\ref{fig-Test-sig-2-Recons-kInt3o0} and \ref{fig-Test-sig-2-Recons-kInt4o0} are all plotted on their own scales.  In both figures we are clearly able to determine whether the inclusion is more or less conductive than the background, as well as its approximate location, even from as little as $25\%$ D-N data.

\begin{figure} % adjust colors on this one?
\hspace{1em}
\begin{picture}(300,230)
\put(15,170){\textbf{Full Data}}
\put(0,90){\includegraphics[width=2.7cm]{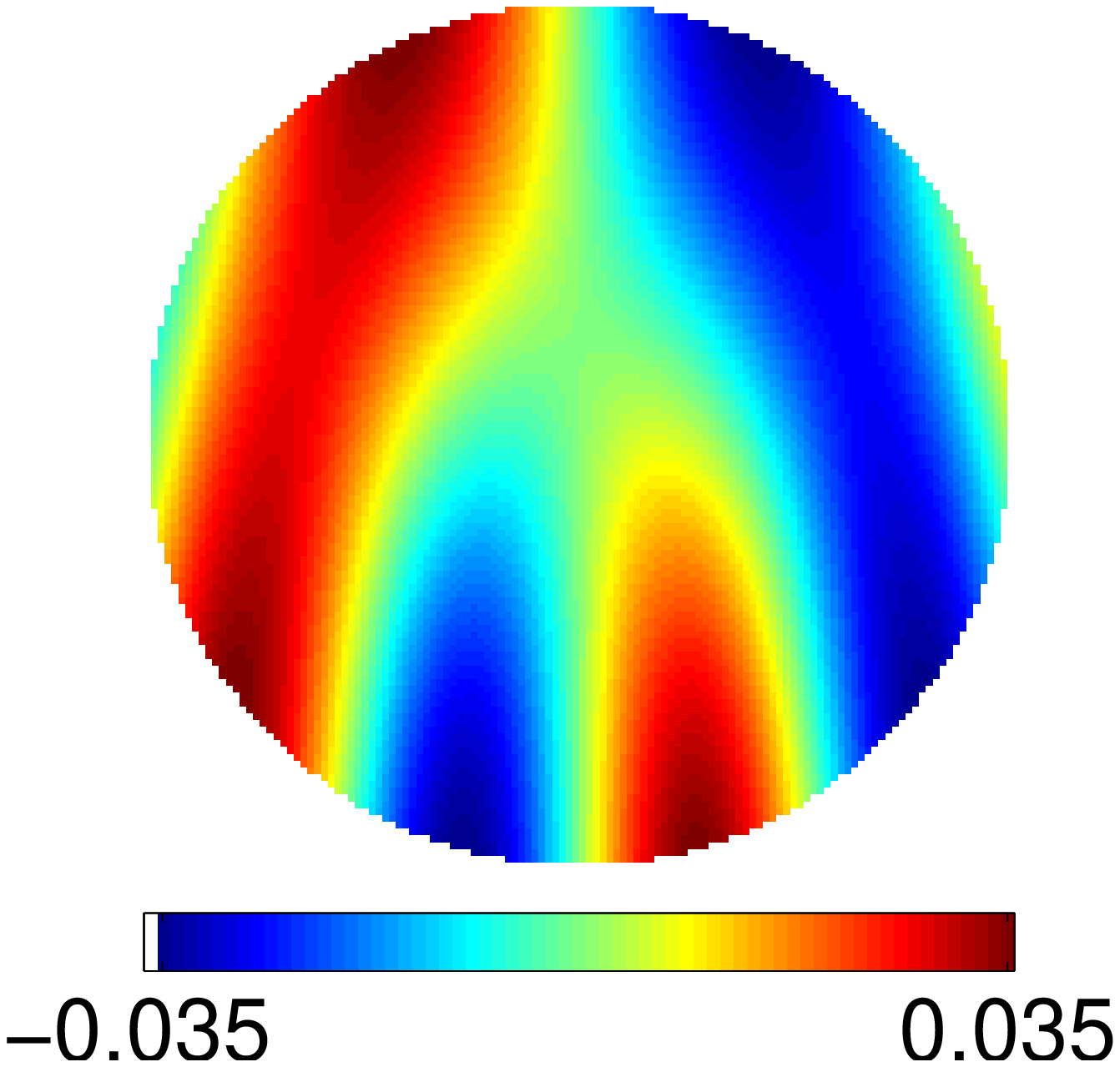}}
\put(0,0){\includegraphics[width=2.7cm]{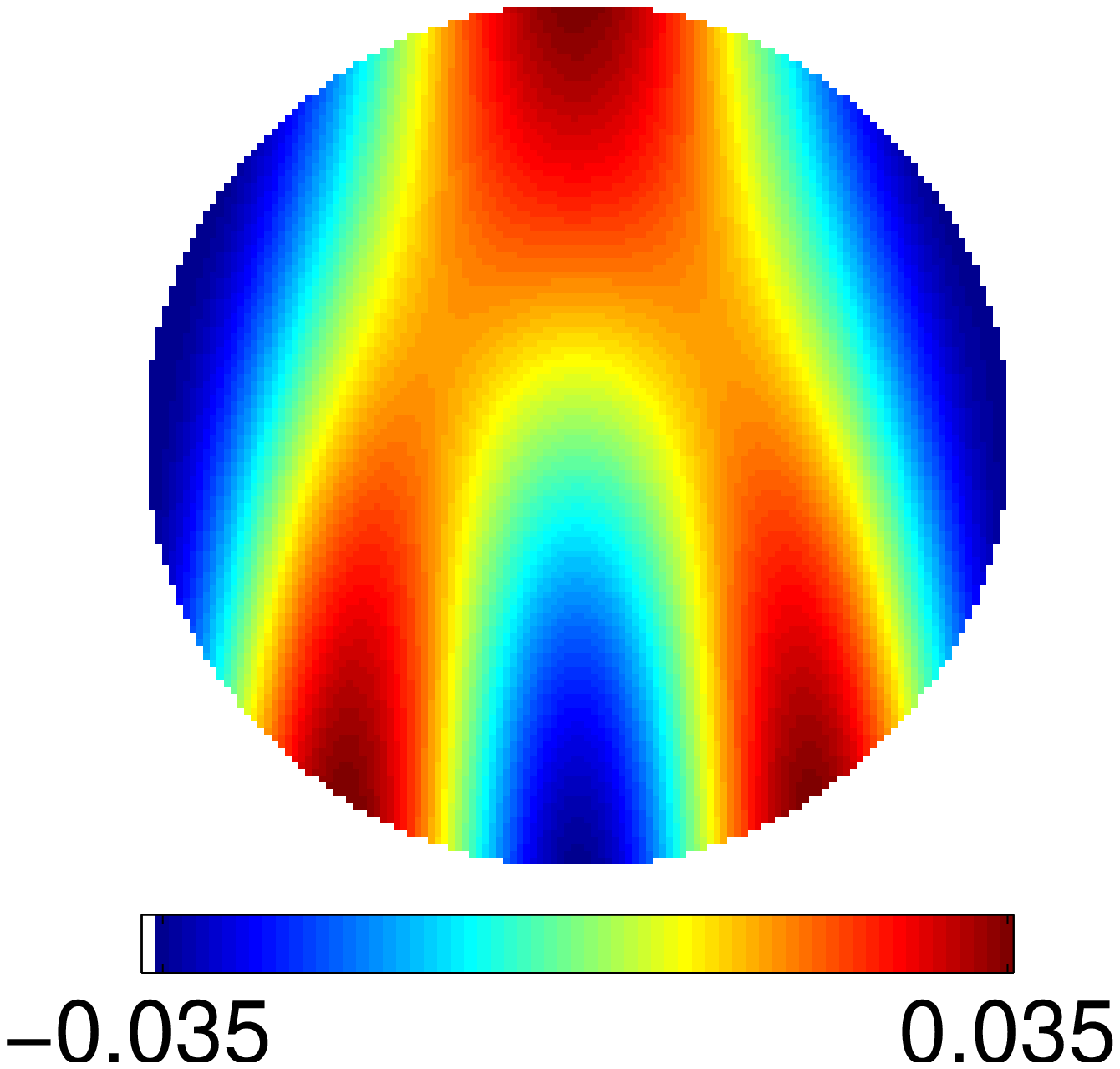}}

\put(102,170){$\mathbf{\frac 34}$ \textbf{Data}}
\put(80,90){\includegraphics[width=2.7cm]{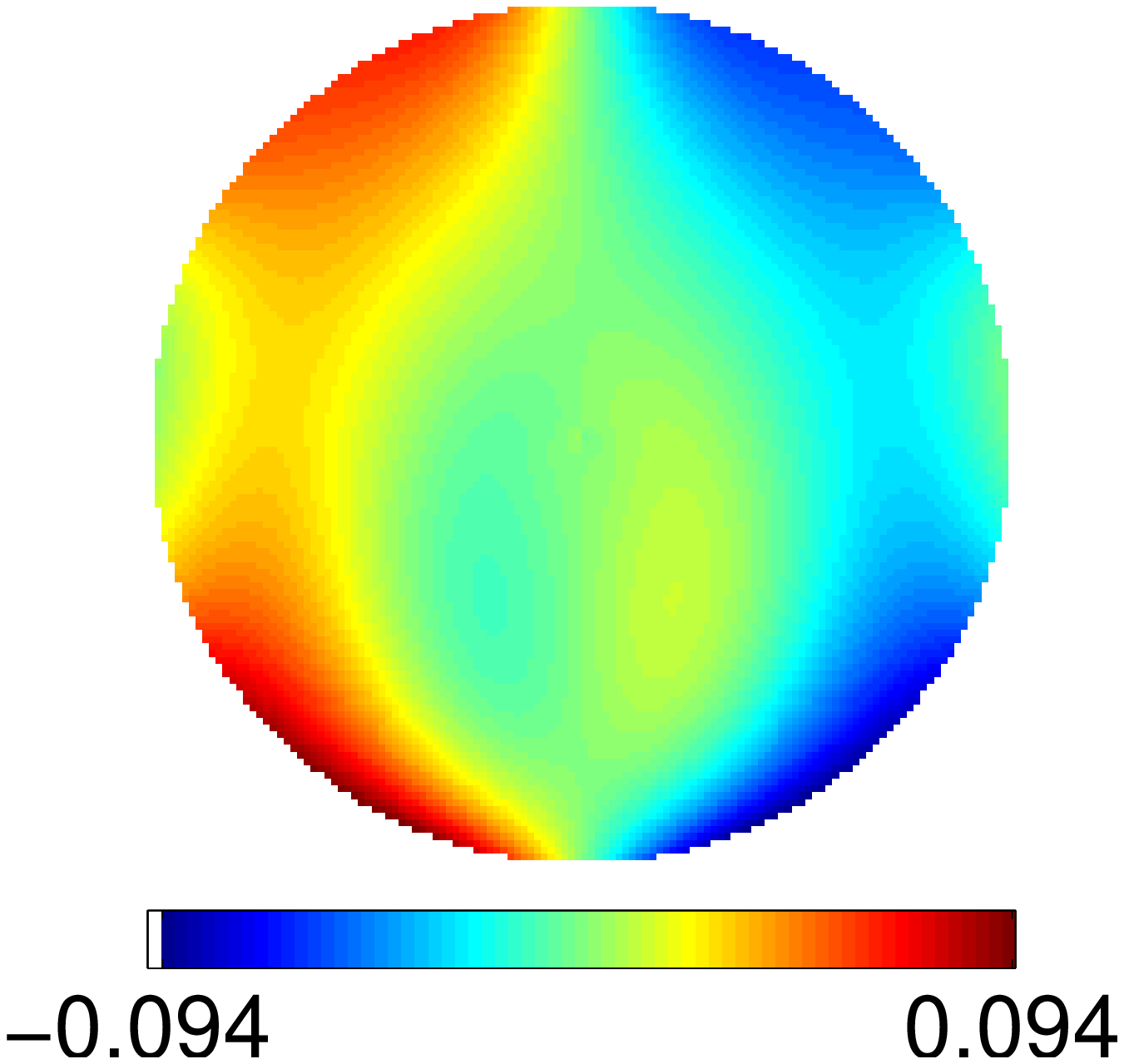}}
\put(80,0){\includegraphics[width=2.7cm]{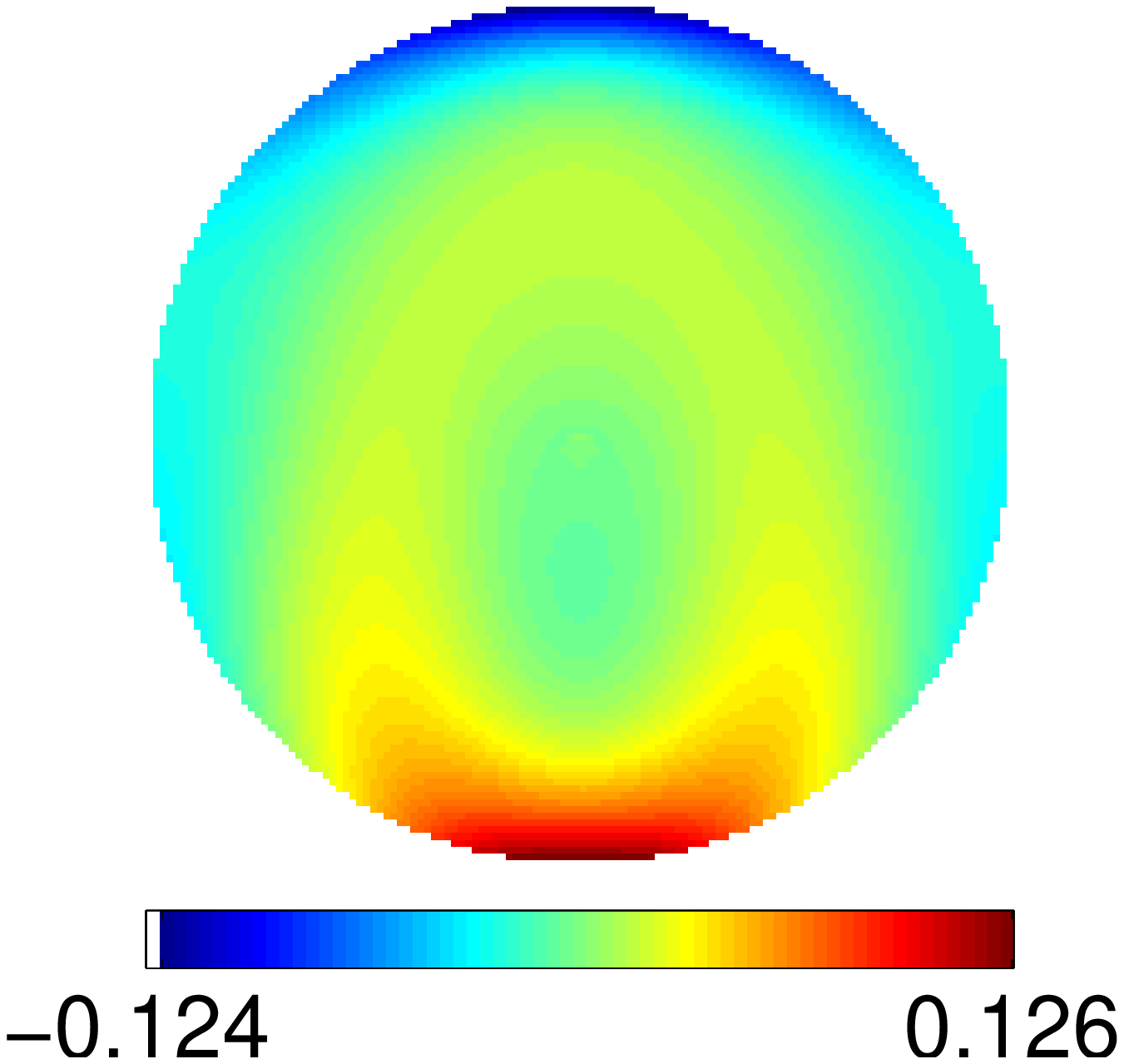}}

\put(180,170){$\mathbf{\frac 12}$ \textbf{Data}}
\put(160,90){\includegraphics[width=2.7cm]{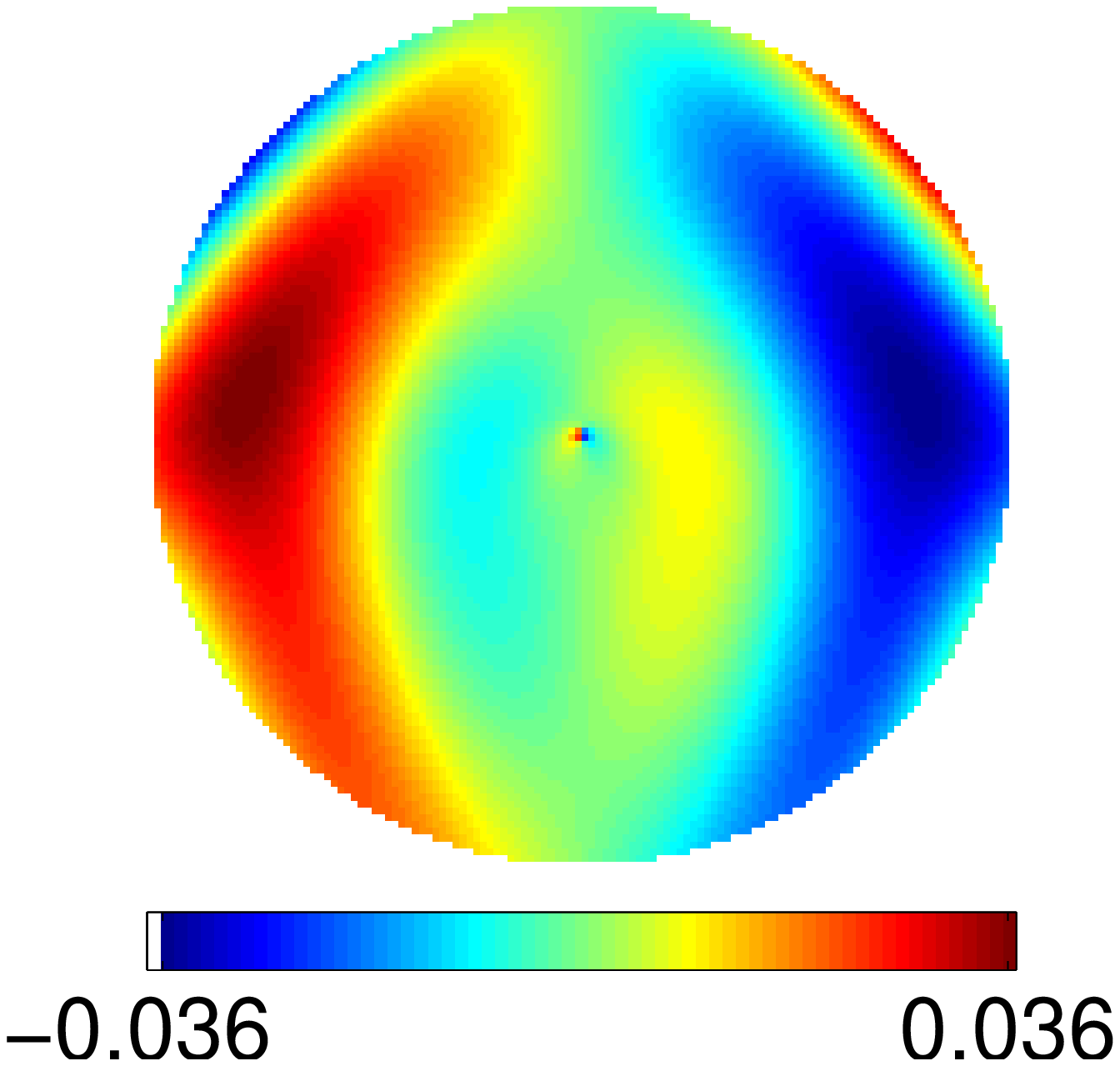}}
\put(160,0){\includegraphics[width=2.7cm]{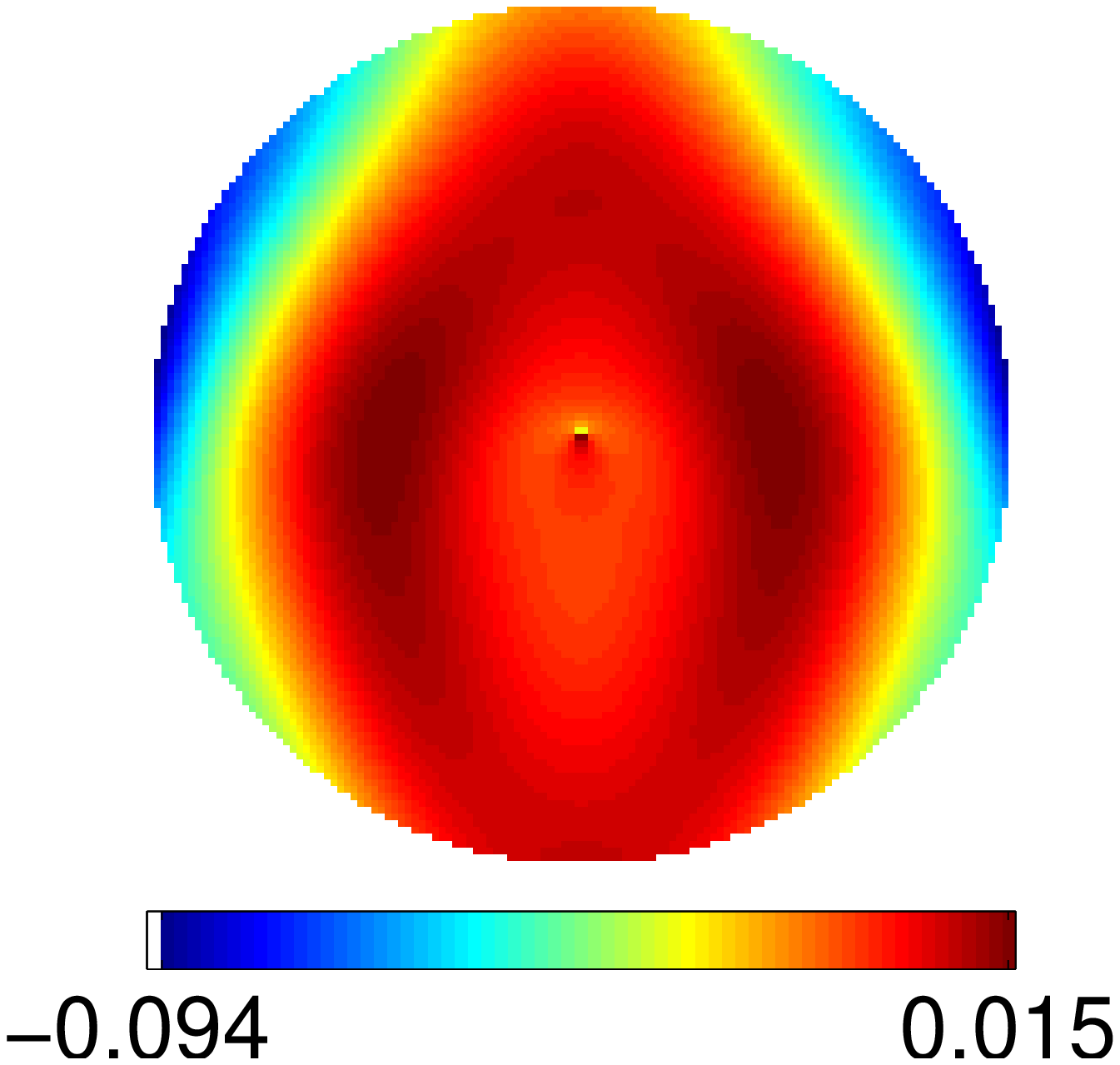}}

\put(263,170){$\mathbf{\frac 14}$ \textbf{Data}}
\put(240,90){\includegraphics[width=2.7cm]{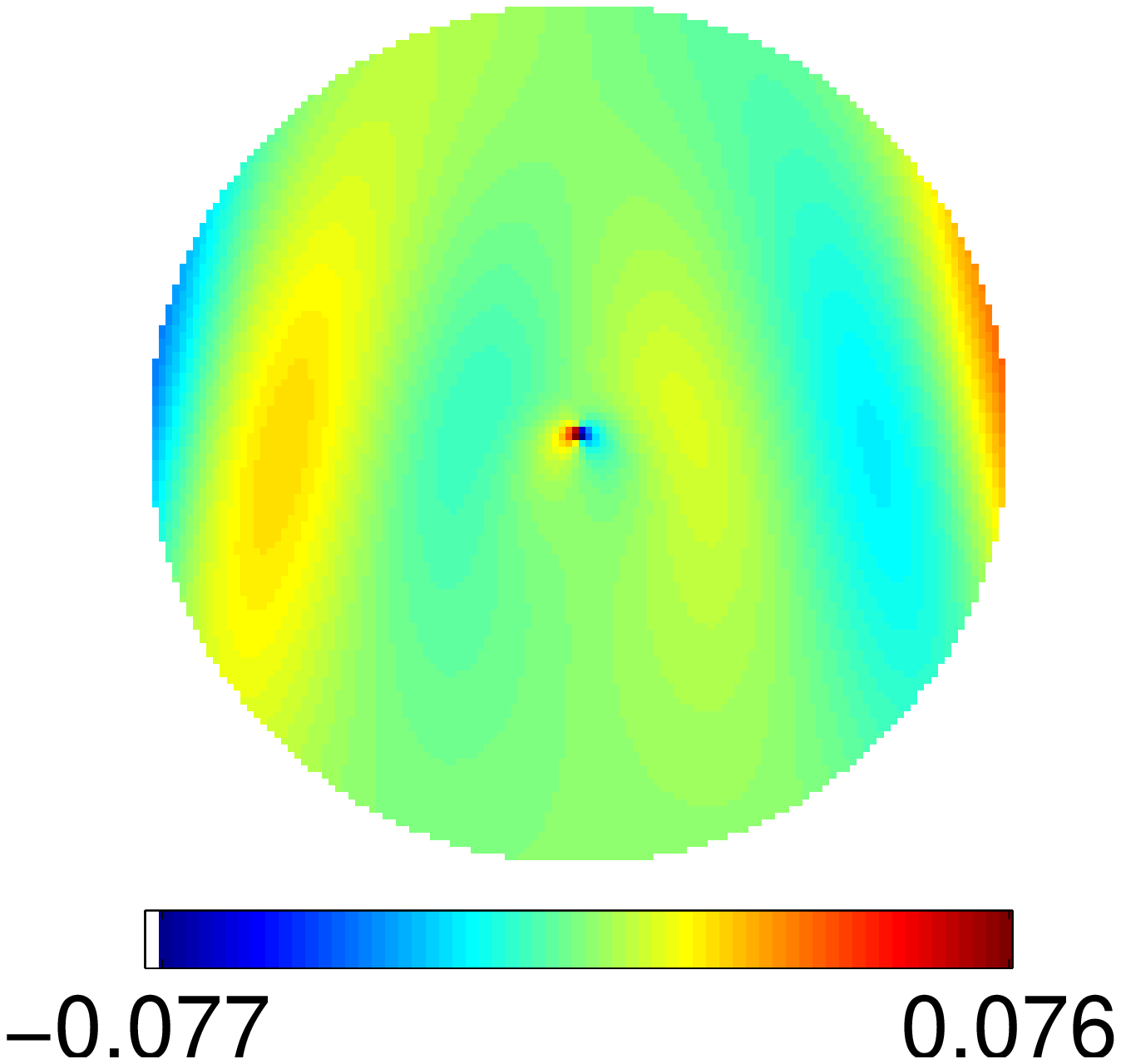}}
\put(240,0){\includegraphics[width=2.7cm]{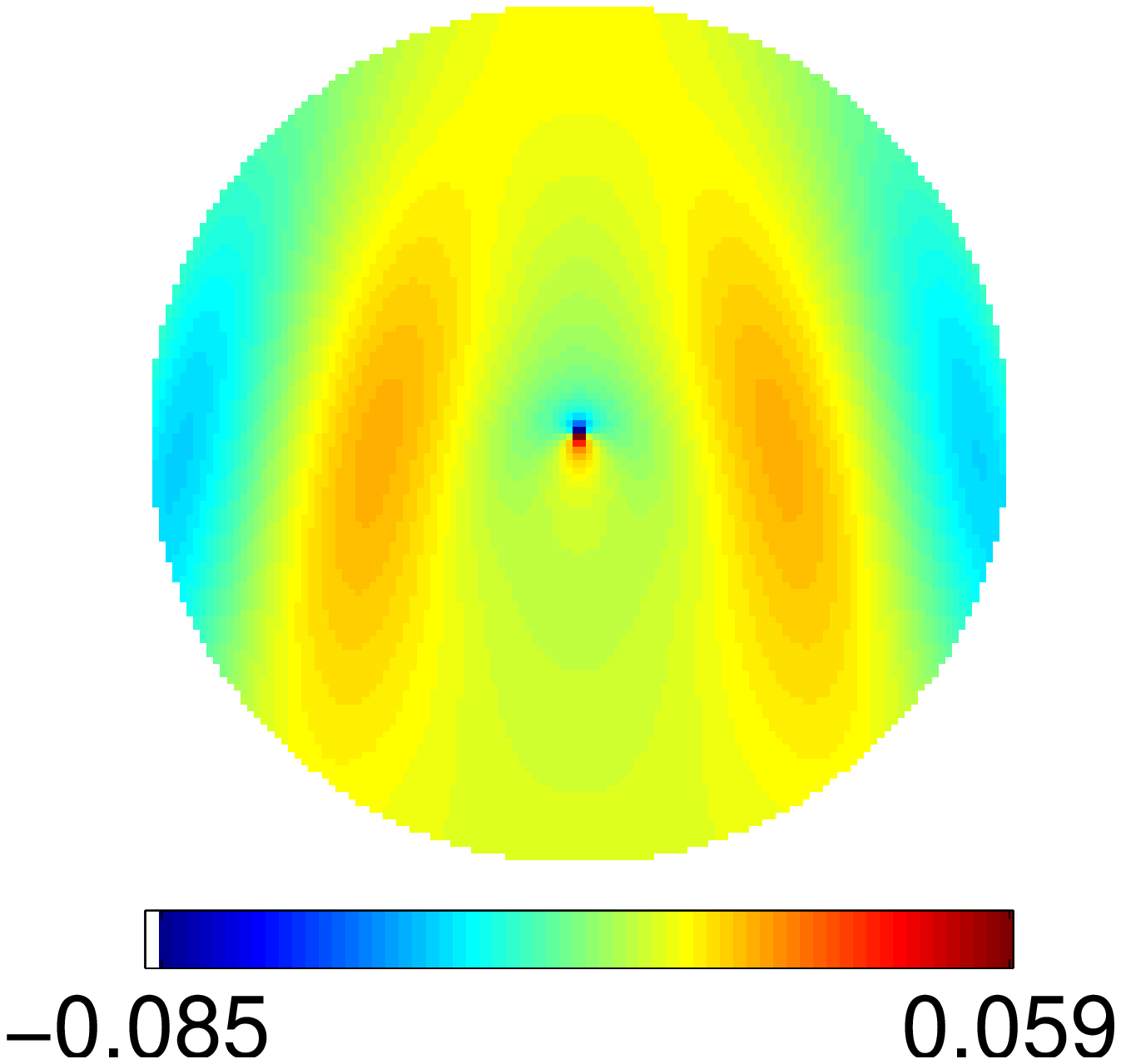}}

\put(-30,140){\textbf{Re} $\mathbf{S_{21}}$}
\put(-30,40){\textbf{Im} $\mathbf{S_{21}}$}
\end{picture}
\caption{Scattering transforms $S_{21}$ for the discontinuous conductivity in Figure~\ref{fig-Test-sig-2}  for $|k|\leq4$ using the method described in Section~\ref{sec-method2-pdata}.  From left to right, the plots are for Full, 3/4, 1/2, and 1/4 Dirichlet-to-Neumann data.}\label{fig-Test-scat-2-kInt4o0}
\end{figure}

Figure \ref{fig-Test-scat-2-kInt4o0} shows the real and imaginary parts of the scattering transform $S_{21}(k)$ for Full, 3/4, 1/2, and 1/4 data with $|k|\leq 4$.  Note that the scattering data is clearly affected by the loss of information in the D-N map, yet the reconstructions of the conductivity (seen in Figures~\ref{fig-Test-sig-2-Recons-kInt3o0} and \ref{fig-Test-sig-2-Recons-kInt4o0}) continue to contain valuable information.  If one were to continue in this direction (using the partial data D-bar algorithm described in Section~\ref{sec-method2-pdata}) a more in-depth study to determine which values of $k$ are admissible in the scattering data is recommended.
%+++++++++++++++++++++++++++++++
%+++++++++++++++++++++++++++++++

%+++++++++++++++++++++++++++++++
%+++++++++++++++++++++++++++++++

%+++++++++++++++++++++++++++++++
%%+++++++++++++++++++++++++++++++

\clearpage
%--------------------------------------------------------------------
\section{Conclusions}
%--------------------------------------------------------------------
The first step in nonlinear EIT imaging uses the voltage and current boundary data to determine the traces of the CGO solutions at the boundary. This is done by solving a boundary integral equation which is a Fredholm equation of the second kind, e.g., \eqref{eq-Psi-INT}.

In this work, we used simulated partial boundary data and a wavelet-based integral equation solver to demonstrate that CGO solutions can be approximately recovered from partial data, on the part of the boundary where the data was acquired. This result is clearly seen in Figures~\ref{fig:tracesA} and \ref{fig:tracesC} for a $C^2$ conductivity, and in Figures~\ref{fig:traces-u1} and \ref{fig:traces-u2} for a discontinuous conductivity.  In addition, we have demonstrated that such partial data CGO solutions can be used in existing full data D-bar methods to provide useful and informative reconstructions, even in the case of discontinuous conductivities, see Figures~\ref{fig-Test-sig-2-Recons-kInt3o0} and \ref{fig-Test-sig-2-Recons-kInt4o0}.
%--------------------------------------------------------------------
\section*{Acknowledgments}
%--------------------------------------------------------------------
\noindent
The study was supported by the SalWe Research Program for Mind and Body (Tekes - the Finnish Funding Agency for Technology and Innovation grant 1104/10) and by the Academy of Finland (Finnish Centre of Excellence in Inverse Problems Research 2012--2017, decision number 250215).
%--------------------------------------------------------------------
% Bibliography
%--------------------------------------------------------------------
\bibliographystyle{amsalpha}
\bibliography{InverseProblemsRefs}
%--------------------------------------------------------------------
%\affiliation{University of Helsinki: Department of Mathematics and Statistics} 
\end{document}